\documentclass[leqno,draft]{article}

\usepackage{makeidx}
\makeindex

\def\diam{\mathop{\rm diam}}

\newtheorem{theorem}{Theorem}
\newtheorem{lemma}[theorem]{Lemma}
\newtheorem{proposition}[theorem]{Proposition}
\newtheorem{definition}[theorem]{Definition}
\newtheorem{corollary}[theorem]{Corollary}

\newcommand{\begintheorem}{\addtocounter{equation}{1}\begin{theorem}}
\newcommand{\beginlemma}{\addtocounter{equation}{1}\begin{lemma}}
\newcommand{\beginproposition}{\addtocounter{equation}{1}\begin{proposition}}
\newcommand{\begindefinition}{\addtocounter{equation}{1}\begin{definition}}
\newcommand{\begincorollary}{\addtocounter{equation}{1}\begin{corollary}}

\begin{document}

\title{Some topics related to metrics and norms, including ultrametrics
and ultranorms, 2}

\author{Stephen Semmes \\
        Rice University}

\date{}

\maketitle

\begin{abstract}
Here we look at (collections of) semimetrics and seminorms, including
their ultrametric versions.  In particular, we are concerned with
geometric properties related to connectedness and topological
dimension $0$.
\end{abstract}

\tableofcontents

\part{Semimetrics and seminorms}
\label{semimetrics, seminorms}

\section{Semimetrics}
\label{semimetrics}
\setcounter{equation}{0}

        Let $X$ be a set, and let $d(x, y)$ be a nonnegative real-valued
function defined for $x, y \in X$.  We say that $d(x, y)$ is a
\emph{semimetric}\index{semimetrics} on $X$ if
\begin{equation}
\label{d(x, x) = 0}
        d(x, x) = 0
\end{equation}
for every $x \in X$,
\begin{equation}
\label{d(x, y) = d(y, x)}
        d(x, y) = d(y, x)
\end{equation}
for every $x, y \in X$, and
\begin{equation}
\label{d(x, z) le d(x, y) + d(y, z)}
        d(x, z) \le d(x, y) + d(y, z)
\end{equation}
for every $x, y, z \in X$.  Of course, (\ref{d(x, z) le d(x, y) + d(y,
  z)}) is known as the \emph{triangle inequality}.\index{triangle inequality}
If we also have that
\begin{equation}
\label{d(x, y) > 0}
        d(x, y) > 0
\end{equation}
for every $x, y \in X$ with $x \ne y$, then $d(x, y)$ is said to be a
\emph{metric}\index{metrics} on $X$.

        Let $d(\cdot, \cdot)$ be a semimetric on $X$, let $x$ be an 
element of $X$, and let $r$ be a positive real number.  The \emph{open
  ball}\index{open balls} in $X$ centered at $x$ with radius $r$
associated to $d(\cdot, \cdot)$ is defined as usual by
\begin{equation}
\label{B(x, r) = B_d(x, r) = {y in X : d(x, y) < r}}
        B(x, r) = B_d(x, r) = \{y \in X : d(x, y) < r\}.
\end{equation}
If $y \in B(x, r)$, then $t = r - d(x, y) > 0$, and it is easy to see
that
\begin{equation}
\label{B(y, t) subseteq B(x, r)}
        B(y, t) \subseteq B(x, r),
\end{equation}
using the triangle inequality.

        A set $U \subseteq X$ is said to be an \emph{open set}\index{open sets}
with respect to the semimetric $d(\cdot, \cdot)$ if for each $x \in U$
there is an $r > 0$ such that
\begin{equation}
\label{B(x, r) subseteq U}
        B(x, r) \subseteq U.
\end{equation}
It is well known and easy to check that this defines a topology on
$X$.  Note that open balls in $X$ with respect to $d(\cdot, \cdot)$
are open sets, by (\ref{B(y, t) subseteq B(x, r)}).  The collection of
open balls with respect to $d(\cdot, \cdot)$ centered at a point $x
\in X$ is a local base for the topology of $X$ determined by $d(\cdot,
\cdot)$ at $x$, and the collection of all open balls in $X$ with
respect to $d(\cdot, \cdot)$ is a base for the topology on $X$
determined by $d(\cdot, \cdot)$.  If $d(\cdot, \cdot)$ is a metric on
$X$, then $X$ is Hausdorff with respect to the corresponding topology.

        Similarly, the \emph{closed ball}\index{closed balls} in $X$
centered at a point $x \in X$ with radius $r \ge 0$ with respect to
a semimetric $d(\cdot, \cdot)$ is defined by
\begin{equation}
\label{overline{B}(x, r) = overline{B}_d(x, r) = {y in X : d(x, y) le r}}
 \overline{B}(x, r) = \overline{B}_d(x, r) = \{y \in X : d(x, y) \le r\}.
\end{equation}
Put
\begin{equation}
\label{V(x, r) = X setminus overline{B}(x, r) = {y in X : d(x, y) > r}}
 V(x, r) = X \setminus \overline{B}(x, r) = \{y \in X : d(x, y) > r\}.
\end{equation}
If $y \in V(x, r)$, then $t = d(x, y) - r > 0$, and one can check that
\begin{equation}
\label{B(y, t) subseteq V(x, r)}
        B(y, t) \subseteq V(x, r),
\end{equation}
using the triangle inequality.  This implies that $V(x, r)$ is an open
set with respect to the topology determined by $d(\cdot, \cdot)$ for
every $x \in X$ and $r \ge 0$, so that $\overline{B}(x, r)$ is a
closed set with respect to this topology.

        Let $X$ be an arbitrary topological space for the moment.
Strictly speaking, one often says that $X$ is \emph{regular}\index{regular
topological spaces} if for each $x \in X$ and closed set $E \subseteq X$
with $x \not\in E$ there are disjoint open sets $U, V \subseteq X$
such that $x \in U$ and $E \subseteq V$.  This is equivalent to asking
that for each $x \in X$ and open set $W \subseteq X$ with $x \in W$
there is an open set $U \subseteq X$ such that $x \in U$ and the closure
$\overline{U}$ of $U$ in $X$ is contained in $W$.  If the topology
on $X$ is determined by a semimetric, then it is easy to see that
$X$ is regular in this sense, by standard arguments.  In particular,
a regular topological space in this sense need not be Hausdorff.

        If $X$ satisfies the first separation condition, then subsets
of $X$ with exactly one element are closed sets, and so regularity of
$X$ as in the preceding paragraph implies that $X$ is Hausdorff.  If
fact, it would be enough to ask that $X$ satisfy the $0$th separation
separation condition for this to work.  Sometimes this may be included
in the definition of regularity, and otherwise one may say that a
topological space $X$ satisfies the third separation condition when
$X$ satisfies the first or $0$th separation condition and $X$ is
regular.

        A topological space $X$ is \emph{normal}\index{normal topological
spaces} in the strict sense if for every pair $A$, $B$ of disjoint closed
subsets of $X$ there are disjoint open sets $U, V \subseteq X$ such that
$A \subseteq U$ and $B \subseteq V$.  If $X$ satisfies the first separation
condition and is normal in this strict sense, then $X$ is Hausdorff and
regular.  As before, the first separation condition can also be included
in the definition of normality, or one may say that $X$ satisfies the fourth
separation condition when $X$ satisfies the first separation condition
and $X$ is normal.  If the topology on $X$ is determined by a semimetric,
then one can check that $X$ is normal in the strict sense, in the same way
as for metric spaces.

        If $d(x, y)$ is a semimetric on $X$ and $Y \subseteq X$, then
the restriction of $d(x, y)$ to $x, y \in Y$ is a semimetric on $Y$.
Let $B_{d, Y}(x, r)$ be the open ball in $Y$ centered at $x \in Y$
with radius $r > 0$ with respect to the restriction of $d(\cdot, \cdot)$
to $Y$, so that
\begin{equation}
\label{B_{d, Y}(x, r) = B_{d, X}(x, r) cap Y}
        B_{d, Y}(x, r) = B_{d, X}(x, r) \cap Y.
\end{equation}
Thus (\ref{B_{d, Y}(x, r) = B_{d, X}(x, r) cap Y}) is an open set in
$Y$ with respect to the topology induced on $Y$ by the topology on $X$
determined by $d(\cdot, \cdot)$.  Every open set in $Y$ with respect
to the topology determined by the restriction of $d(\cdot, \cdot)$ to
$Y$ can be expressed as a union of open balls in $Y$, and hence is an
open set with respect to the topology induced on $Y$ by the topology
on $X$ determined by $d(\cdot, \cdot)$.  Conversely, it is easy to see
that every open set in $Y$ with respect to the topology induced on $Y$
by the topology on $X$ determined by $d(\cdot, \cdot)$ is an open set
with respect to the topology on $Y$ determined by the restriction of
$d(\cdot, \cdot)$ to $Y$, directly from the definitions.

\section{Collections of semimetrics}
\label{collections of semimetrics}
\setcounter{equation}{0}

        Let $X$ be a set, and let $l$ be a positive integer.  If $d_j(x, y)$
is a semimetric on $X$ for $j = 1, \ldots, l$, then it is easy to see that
\begin{equation}
\label{d(x, y) = max_{1 le j le l} d_j(x, y)}
        d(x, y) = \max_{1 \le j \le l} d_j(x, y)
\end{equation}
is a semimetric on $X$ as well.  Observe that
\begin{equation}
\label{B_d(x, r) = bigcap_{j = 1}^l B_{d_j}(x, r)}
        B_d(x, r) = \bigcap_{j = 1}^l B_{d_j}(x, r)
\end{equation}
for every $x \in X$ and $r > 0$, where $B_d(x, r)$ and $B_{d_j}(x, r)$
are as in (\ref{B(x, r) = B_d(x, r) = {y in X : d(x, y) < r}}).
Similarly,
\begin{equation}
\label{widetilde{d}(x, y) = sum_{j = 1}^l d_j(x, y)}
        \widetilde{d}(x, y) = \sum_{j = 1}^l d_j(x, y)
\end{equation}
is a semimetric on $X$, and
\begin{equation}
\label{d(x, y) le widetilde{d}(x, y) le l d(x, y)}
        d(x, y) \le \widetilde{d}(x, y) \le l \, d(x, y)
\end{equation}
for every $x, y \in X$.  This implies that $d(x, y)$ and
$\widetilde{d}(x, y)$ determine the same topology on $X$.

        Now let $\mathcal{M}$ be a nonempty collection of semimetrics
on $X$.  Let us say that a set $U \subseteq X$ is an \emph{open
  set}\index{open sets} with respect to $\mathcal{M}$ if for each $x
\in U$ there are finitely many elements $d_1, \ldots, d_l$ of
$\mathcal{M}$ and finitely many positive real numbers $r_1, \ldots,
r_l$ such that
\begin{equation}
\label{bigcap_{j = 1}^l B_{d_j}(x, r_j) subseteq U}
        \bigcap_{j = 1}^l B_{d_j}(x, r_j) \subseteq U.
\end{equation}
It is easy to see that this defines a topology on $X$.  If
$\mathcal{M} = \emptyset$, then we interpret the corresponding
topology on $X$ as being the indiscrete topology.  If $\mathcal{M}$
consists of a single semimetric on $X$, then this is the same as the
topology determined on $X$ by that semimetric, as in the previous
section.  If $\mathcal{M}$ consists of finitely many semimetrics on
$X$, then the topology on $X$ associated to $\mathcal{M}$ is the same
as the topology determined on $X$ by the maximum or sum of the
elements of $\mathcal{M}$.  This follows from the remarks in the
preceding paragraph, since one can take the $r_j$'s in (\ref{bigcap_{j
    = 1}^l B_{d_j}(x, r_j) subseteq U}) to be equal to each other.  If
$\mathcal{M}$ is any collection of semimetrics on $X$, then every open
set in $X$ with respect to the topology determined by any $d \in
\mathcal{M}$ is also an open set in $X$ with respect to the topology
associated to $\mathcal{M}$.  In particular, open balls in $X$ with
respect to the elements of $\mathcal{M}$ are open sets with respect to
the topology associated to $\mathcal{M}$.

        Let $\mathcal{M}$ be a nonempty collection of semimetrics
on $X$ again, and let $x \in X$ be given, along with finitely many
elements $d_1, \ldots, d_l$ of $\mathcal{M}$ and positive real numbers
$r_1, \ldots, r_l$.  Thus $B_{d_j}(x, r_j)$ is an open set with
respect to $d_j$ for $j = 1, \ldots, l$, and hence with respect to the
topology associated to $\mathcal{M}$, as before.  This implies that
\begin{equation}
\label{bigcap_{j = 1}^l B_{d_j}(x, r_j)}
        \bigcap_{j = 1}^l B_{d_j}(x, r_j)
\end{equation}
is an open set in $X$ with respect to the topology associated to
$\mathcal{M}$ too.  If $x \in X$ is fixed, then the collection of open
sets of the form (\ref{bigcap_{j = 1}^l B_{d_j}(x, r_j)}) is a local
base for the topology on $X$ associated to $\mathcal{M}$ at $x$.  The
collection of open sets of the form (\ref{bigcap_{j = 1}^l B_{d_j}(x,
  r_j)}) for any $x \in X$ is a base for the topology on $X$
associated to $\mathcal{M}$.

        Equivalently, the collection of open balls $B_d(x, r)$
with $d \in \mathcal{M}$ and $r > 0$ is a local sub-base for the
topology on $X$ associated to $\mathcal{M}$ at $x$.  Similarly, the
collection of open balls $B_d(x, r)$ with $x \in X$, $d \in
\mathcal{M}$, and $r > 0$ is a sub-base for the topology on $X$
associated to $\mathcal{M}$.  This is the same as saying that the
collection of finite intersections of open balls in $X$ corresponding
to elements of $\mathcal{M}$ is a base for the topology on $X$
associated to $\mathcal{M}$.  This is slightly less precise than
using intersections of the form (\ref{bigcap_{j = 1}^l B_{d_j}(x, r_j)}),
where the balls are centered at the same point in $X$.

        Every closed set in $X$ with respect to any $d \in \mathcal{M}$
is a closed set in $X$ with respect to the toplogy on $X$ associated
to $\mathcal{M}$, because of the analogous statement for open sets.
This includes closed balls in $X$ with respect to any $d \in
\mathcal{M}$, as in the previous section.  It follows that the
intersection of any family of closed balls in $X$ with respect to
elements of $\mathcal{M}$ is a closed set in $X$ with respect to the
topology associated to $\mathcal{M}$ as well.  Using this, one can
check that $X$ is regular in the strict sense discussed in the
previous section, with respect to the topology associated to
$\mathcal{M}$.

        Let us say that $\mathcal{M}$ is
\emph{nondegenerate}\index{nondegenerate collections of semimetrics}
on $X$ if for each pair of distinct elements $x$, $y$ of $X$
there is a $d \in \mathcal{M}$ such that
\begin{equation}
\label{d(x, y) > 0, 2}
        d(x, y) > 0.
\end{equation}
This implies that $X$ is Hausdorff with respect to the topology
associated to $\mathcal{M}$, by essentially the same argument as for
metric spaces.  If $\mathcal{M}$ consists of finitely many semimetrics
on $X$, then $\mathcal{M}$ is nondegenerate exactly when the sum or
maximum of these semimetrics is a metric on $X$.

        If $\mathcal{M}$ is any collection of semimetrics on $X$
and $Y \subseteq X$, then let
\begin{equation}
\label{mathcal{M}_Y}
        \mathcal{M}_Y
\end{equation}
be the collection of semimetrics on $Y$ obtained by restricting the
elements of $\mathcal{M}$ to $Y$.  One can check that the topology on
$Y$ associated to $\mathcal{M}_Y$ as before is the same as the
topology induced on $Y$ by the topology on $X$ associated to
$\mathcal{M}$.  Of course, this is trivial when $\mathcal{M} =
\emptyset$, and so we suppose that $\mathcal{M} \ne \emptyset$.  The
argument is analogous to the one for a single semimetric, as in the
previous section.  More precisely, if $E \subseteq Y$ is an open set
with respect to the topology induced on $Y$ by the topology on $X$
associated to $\mathcal{M}$, then it is easy to see that $E$ is an
open set with respect to the topology on $Y$ associated to
$\mathcal{M}_Y$, directly from the definitions.  Conversely, if $E
\subseteq Y$ is an open set with respect to the topology associated to
$\mathcal{M}_Y$, then one would like to verify that $E$ is an open set
with respect to the topology induced on $Y$ by the topology on $X$
associated to $\mathcal{M}$.  If $E$ is an open ball in $Y$ centered
at a point in $Y$ with respect to an element of $\mathcal{M}_Y$, then
this follows from (\ref{B_{d, Y}(x, r) = B_{d, X}(x, r) cap Y}).
Similarly, if $E$ is the intersection of finitely many open balls in
$Y$ with respect to elements of $\mathcal{M}_Y$, then $E$ can be
expressed as the intersection of $Y$ with finitely many open balls in
$X$ with respect to elements of $\mathcal{M}$, which implies that $E$
is an open set in $Y$ with respect to the topology induced on $Y$ by
the topology on $X$ associated to $\mathcal{M}$.  As before, the
collection of finite intersection of open balls in $Y$ with respect to
elements of $\mathcal{M}_Y$ is a base for the topology on $Y$
associated to $\mathcal{M}_Y$, so that every open set $E$ in $Y$ with
respect to the topology associated to $\mathcal{M}_Y$ can be expressed
as a union of finite intersections of open balls in $Y$ with respect
to elements of $\mathcal{M}_Y$.  It follows that $E$ is an open set
with respect to the topology induced on $Y$ by the topology on $X$
associated to $\mathcal{M}$, as desired, because $E$ can be expressed
as a union of open sets with respect to the induced topology.

\section{Seminorms}
\label{seminorms}
\setcounter{equation}{0}

        Let $V$ be a vector space over either the real numbers
${\bf R}$\index{R@${\bf R}$} or the complex numbers
${\bf C}$.\index{C@${\bf C}$}  A nonnegative real-valued function
$N$ on $V$ is said to be a \emph{seminorm}\index{seminorms} on $V$
if it satisfies the following two conditions.  First,
\begin{equation}
\label{N(t v) = |t| N(v)}
        N(t \, v) = |t| \, N(v)
\end{equation}
for every $v \in V$ and real or complex number $t$, as appropriate,
where $|t|$ is the usual absolute value of $t$.  Second,
\begin{equation}
\label{N(v + w) le N(v) + N(w)}
        N(v + w) \le N(v) + N(w)
\end{equation}
for every $v, w \in V$, which is another version of the triangle
inequality.\index{triangle inequality} Note that (\ref{N(t v) = |t|
  N(v)}) implies that $N(0) = 0$, by taking $t = 0$.  If we also have
that
\begin{equation}
\label{N(v) > 0}
        N(v) > 0
\end{equation}
for every $v \in V$ with $v \ne 0$, then $N$ is said to be a
\emph{norm}\index{norms} on $V$.  If $N$ is a seminorm on $V$, then it
is easy to see that
\begin{equation}
\label{d(v, w) = N(v - w)}
        d(v, w) = N(v - w)
\end{equation}
defines a semimetric on $V$.  Similarly, if $N$ is a norm on $V$, then
(\ref{d(v, w) = N(v - w)}) defines a metric on $V$.

        If $N_1, \ldots, N_l$ are finitely many seminorms on $V$, then
it is easy to see that
\begin{equation}
\label{N(v) = max_{1 le j le n} N_j(v)}
        N(v) = \max_{1 \le j \le n} N_j(v)
\end{equation}
and
\begin{equation}
\label{widetilde{N}(v) = sum_{j = 1}^l N_j(v)}
        \widetilde{N}(v) = \sum_{j = 1}^l N_j(v)
\end{equation}
are seminorms on $V$ too.  Let $d(v, w)$ be the semimetric on $V$
corresponding to $N$ as in (\ref{d(v, w) = N(v - w)}), let
\begin{equation}
\label{d_j(v, w) = N_j(v - w)}
        d_j(v, w) = N_j(v - w)
\end{equation}
be the semimetric corresponding to $N_j$ for each $j = 1, \ldots, l$,
and let
\begin{equation}
\label{widetilde{d}(v, w) = widetilde{N}(v - w)}
        \widetilde{d}(v, w) = \widetilde{N}(v - w)
\end{equation}
be the semimetric corresponding to $\widetilde{N}$.  Thus $d(v, w)$
and $\widetilde{d}(v, w)$ can be given in terms of the $d_j(v, w)$'s
as in (\ref{d(x, y) = max_{1 le j le l} d_j(x, y)}) and
(\ref{widetilde{d}(x, y) = sum_{j = 1}^l d_j(x, y)}).  As before, we
also have that
\begin{equation}
\label{N(v) le widetilde{N}(v) le l N(v)}
        N(v) \le \widetilde{N}(v) \le l \, N(v)
\end{equation}
for every $v \in V$.

        Now let $\mathcal{N}$ be a collection of seminorms on $V$, and
let $\mathcal{M}$ be the corresponding collection of semimetrics on
$V$, as in (\ref{d(v, w) = N(v - w)}).  This leads to a topology on
$V$, as in the previous section.  If $\mathcal{N}$ consists of
finitely many seminorms on $V$, then one could get the same topology
on $V$ using a single seminorm, as in (\ref{N(v) = max_{1 le j le n}
  N_j(v)}) or (\ref{widetilde{N}(v) = sum_{j = 1}^l N_j(v)}).  Let us
say that $\mathcal{N}$ is
\emph{nondegenerate}\index{nondegenerate collections of seminorms} on
$V$ if for each $v \in V$ with $v \ne 0$ there is an $N \in
\mathcal{N}$ such that (\ref{N(v) > 0}) holds.  This implies that
$\mathcal{M}$ is nondegenerate as a collection of semimetrics on $V$,
as in the previous section, and hence that the associated topology on
$V$ is Hausdorff.

        If $W$ is a linear subspace of $V$ and $N$ is a seminorm on $V$,
then the restriction of $N$ to $W$ defines a seminorm on $W$ too.  Let
$\mathcal{N}$ be a collection of seminorms on $V$ again, and let
$\mathcal{N}_W$ be the collection of seminorms on $W$ obtained by
restricting the elements of $\mathcal{N}$ to $W$.  Also let
$\mathcal{M}$ be the collection of semimetrics on $V$ corresponding to
$\mathcal{N}$ as before, and let $\mathcal{M}_W$ be the collection of
semimetrics on $W$ that correspond to elements of $\mathcal{N}_W$ in
the same way.  Equivalently, $\mathcal{M}_W$ is the same as the
collection of semimetrics on $W$ obtained by restricting the elements
of $\mathcal{M}$ to $W$, as in the previous section.  Thus the
topology on $W$ associated to $\mathcal{N}_W$ is the same as the
topology induced on $W$ by the topology on $V$ associated to
$\mathcal{N}$, as discussed in the previous section again.

\section{Semi-ultrametrics}
\label{semi-ultrametrics}
\setcounter{equation}{0}

        A semimetric $d(x, y)$ on a set $X$ is said to be a
\emph{semi-ultrametric}\index{semi-ultrametrics} on $X$ if
\begin{equation}
\label{d(x, z) le max(d(x, y), d(y, z))}
        d(x, z) \le \max(d(x, y), d(y, z))
\end{equation}
for every $x, y, z \in X$.  Note that (\ref{d(x, z) le max(d(x, y),
  d(y, z))}) implies the ordinary triangle inequality (\ref{d(x, z) le
  d(x, y) + d(y, z)}).  An \emph{ultrametric}\index{ultrametrics} on
$X$ is a metric that satisfies (\ref{d(x, z) le max(d(x, y), d(y,
  z))}).  Remember that the \emph{discrete metric}\index{discrete metric}
is defined on $X$ by putting $d(x, y)$ equal to $1$ when $x \ne y$,
and to $0$ when $x = y$.  It is easy to see that this defines an
ultrametric on $X$, for which the corresponding topology is the
discrete topology.  One can check that the maximum of finitely many
semi-ultrametrics on $X$ is also a semi-ultrametric on $X$, as in
Section \ref{collections of semimetrics}.  However, the sum of
finitely many semi-ultrametrics on $X$ is not necessarily a
semi-ultrametric on $X$.

        Let $d(\cdot, \cdot)$ be a semi-ultrametric on a set $X$,
and let $x \in X$ and $r > 0$ be given.  If $y \in B(x, r)$, then
it is easy to see that
\begin{equation}
\label{B(y, r) subseteq B(x, r)}
        B(y, r) \subseteq B(x, r),
\end{equation}
using the ultrametric version of the triangle inequality (\ref{d(x, z)
  le max(d(x, y), d(y, z))}).  More precisely, if $d(x, y) < r$, then
\begin{equation}
\label{B(x, r) = B(y, r)}
        B(x, r) = B(y, r),
\end{equation}
since both (\ref{B(y, r) subseteq B(x, r)}) and the opposite inclusion
hold, for the same reasons.  Similarly, if $y \in \overline{B}(x, r)$
for some $r \ge 0$, then
\begin{equation}
\label{overline{B}(y, r) subseteq overline{B}(x, r)}
        \overline{B}(y, r) \subseteq \overline{B}(x, r),
\end{equation}
by the ultrametric version of the triangle inequality.  As before, we
get that
\begin{equation}
\label{overline{B}(x, r) = overline{B}(y, r)}
        \overline{B}(x, r) = \overline{B}(y, r)
\end{equation}
when $d(x, y) \le r$ for some $r \ge 0$, since both
(\ref{overline{B}(y, r) subseteq overline{B}(x, r)}) and the opposite
inclusion hold.  Note that (\ref{overline{B}(y, r) subseteq
  overline{B}(x, r)}) implies that $\overline{B}(x, r)$ is an open set
when $r > 0$, with respect to the topology determined by $d(\cdot,
\cdot)$.  One can also check that $B(x, r)$ is a closed set in $X$ for
every $x \in X$ and $r > 0$, and we shall return to this in a moment.

        If $x, y, z \in X$ satisfy $d(y, z) \le d(x, y)$, then
(\ref{d(x, z) le max(d(x, y), d(y, z))}) implies that
\begin{equation}
\label{d(x, z) le d(x, y)}
        d(x, z) \le d(x, y).
\end{equation}
Of course, we also have that
\begin{equation}
\label{d(x, y) le max(d(x, z), d(z, y))}
        d(x, y) \le \max(d(x, z), d(z, y)),
\end{equation}
by interchanging the roles of $y$ and $z$ in (\ref{d(x, z) le max(d(x,
  y), d(y, z))}).  If
\begin{equation}
\label{d(y, z) < d(x, y)}
        d(y, z) < d(x, y),
\end{equation}
then (\ref{d(x, y) le max(d(x, z), d(z, y))}) implies that
\begin{equation}
\label{d(x, y) le d(x, z)}
        d(x, y) \le d(x, z).
\end{equation}
Combining (\ref{d(x, z) le d(x, y)}) and (\ref{d(x, y) le d(x, z)}),
we get that
\begin{equation}
\label{d(x, y) = d(x, z)}
        d(x, y) = d(x, z)
\end{equation}
when $x, y, z \in X$ satisfy (\ref{d(y, z) < d(x, y)}).

        Let $x \in X$ and $r \ge 0$ be given, and let $V(x, r)$ be as in
(\ref{V(x, r) = X setminus overline{B}(x, r) = {y in X : d(x, y) > r}}).
If $y \in V(x, r)$, then (\ref{B(y, t) subseteq V(x, r)}) holds with
$t = d(x, y)$, by (\ref{d(x, y) = d(x, z)}).  Similarly, put
\begin{equation}
\label{W(x, r) = X setminus B(x, r) = {y in X : d(x, y) ge r}}
        W(x, r) = X \setminus B(x, r) = \{y \in X : d(x, y) \ge r\}
\end{equation}
for each $r > 0$.  If $y \in W(x, r)$, then it is easy to see that
\begin{equation}
\label{B(y, d(x, y)) subseteq W(x, r)}
        B(y, d(x, y)) \subseteq W(x, r),
\end{equation}
using (\ref{d(x, y) = d(x, z)}) again.  This implies that $W(x, r)$ is
an open set in $X$ with respect to the topology determined by
$d(\cdot, \cdot)$, so that $B(x, r)$ is a closed set in $X$.

\section{Absolute value functions}
\label{absolute value functions}
\setcounter{equation}{0}

        Let $k$ be a field.  A nonnegative real-valued function $|x|$
defined on $k$ is said to be an \emph{absolute value function}\index{absolute
value functions} on $k$ if
\begin{equation}
\label{|x| = 0 if and only if x = 0}
        |x| = 0 \quad\hbox{if and only if}\quad x = 0,
\end{equation}
and
\begin{eqnarray}
\label{|x y| = |x| |y|}
        |x \, y| & = & |x| \, |y|, \\
\label{|x + y| le |x| + |y|}
        |x + y| & \le & |x| + |y|
\end{eqnarray}
for every $x, y \in k$.  It is well known that the standard absolute
value functions on the real and complex numbers satisfy these
conditions.  If $k$ is any field, then the \emph{trivial absolute
  value function}\index{trivial absolute value function} is defined by
putting $|x| = 1$ for every $x \in k$ with $x \ne 0$ and $|0| = 0$.
It is easy to see that this satisfies the requirements of an absolute
value function just mentioned.

        Let $|\cdot|$ be an absolute value function on a field $k$.
We have already used $0$ to refer to the additive identity elements in
$k$ or ${\bf R}$ in the preceding paragraph, and we shall use
$1$ to refer to the multiplicative identity elements in $k$ or ${\bf
  R}$, depending on the context.  If $1$ is the multiplicative
identity element in $k$, then $1 \ne 0$ in $k$, by definition of a
field, and hence $|1| > 0$.  We also have that $1^2 = 1$ in $k$,
so that $|1| = |1^2| = |1|^2$, which implies that
\begin{equation}
\label{|1| = 1}
        |1| = 1.
\end{equation}
Similarly, if $x \in k$ satisfies $x^n = 1$ for some positive integer
$n$, then
\begin{equation}
\label{|x|^n = |x^n| = |1| = 1}
        |x|^n = |x^n| = |1| = 1,
\end{equation}
so that $|x| = 1$.

        The additive inverse of $x \in k$ is denoted $-x$, as usual,
so that
\begin{equation}
\label{(-1) x = -x}
        (-1) \, x = -x
\end{equation}
for every $x \in k$.  In particular, $(-1)^2 = 1$, which implies that
\begin{equation}
\label{|-1| = 1}
        |-1| = 1,
\end{equation}
as in the previous paragraph.  Combining this with (\ref{(-1) x = -x}),
we get that
\begin{equation}
\label{|-x| = |x|}
        |-x| = |x|
\end{equation}
for every $x \in k$.  It follows that
\begin{equation}
\label{d(x, y) = |x - y|}
        d(x, y) = |x - y|
\end{equation}
defines a metric on $k$, using (\ref{|-x| = |x|}) to get that
(\ref{d(x, y) = |x - y|}) is symmetric in $x$ and $y$.

        If
\begin{equation}
\label{|x + y| le max(|x|, |y|)}
        |x + y| \le \max(|x|, |y|)
\end{equation}
for every $x, y \in k$, then we say that $|\cdot|$ defines an
\emph{ultrametric absolute value function}\index{ultrametric absolute
  value functions} on $k$.  This condition implies the ordinary
triangle inequality (\ref{|x + y| le |x| + |y|}), and that the
associated metric (\ref{|x + y| le max(|x|, |y|)}) is an ultrametric
on $k$.  It is easy to see that the trivial absolute value function on
$k$ is an ultrametric absolute value function, which corresponds to
the discrete metric on $k$.  It is well known that the \emph{$p$-adic
  absolute value function}\index{p-adic absolute value
  function@$p$-adic absolute value function} defines an ultrametric
absolute value function on the field ${\bf Q}$\index{Q@${\bf Q}$} of
rational numbers for every prime number $p$, for which the
corresponding ultrametric is known as the \emph{$p$-adic
  metric}.\index{p-adic metric@$p$-adic metric} The field ${\bf
  Q}_p$\index{Q_p@${\bf Q}_p$} of \emph{$p$-adic numbers}\index{p-adic
  numbers@$p$-adic numbers} is obtained by completing ${\bf Q}$ with
respect to the $p$-adic metric, and the $p$-adic absolute value
function can be extended to an ultrametric absolute value function on
${\bf Q}_p$ in a natural way.

        Let $|\cdot|$ be an ultrametric absolute value function on 
a field $k$.  If $y, z \in k$ satisfy
\begin{equation}
\label{|y - z| < |y|}
        |y - z| < |y|,
\end{equation}
then
\begin{equation}
\label{|y| = |z|}
        |y| = |z|.
\end{equation}
This follows from (\ref{d(x, y) = d(x, z)}) in the previous section,
with $x = 0$ and $d(\cdot, \cdot)$ as in (\ref{d(x, y) = |x - y|}).
Of course, one can also verify (\ref{|y| = |z|}) more directly in this
situation.

        Let ${\bf Z}_+$\index{Z_+@${\bf Z}_+$} be the set of positive
integers, and let $n \cdot x$ be the sum of $n$ $x$'s in a field $k$
for each $n \in {\bf Z}_+$ and $x \in k$.  An absolute value function
$|\cdot|$ on $k$ is said to be \emph{archimedian}\index{archimedian
absolute value functions} if the set of nonnegative real numbers of
the form $|n \cdot 1|$ with $n \in {\bf Z}_+$ has a finite upper
bound, and otherwise $|\cdot|$ is said to be
\emph{non-archimedian}\index{non-archimedian absolute value functions}
on $k$.  If $|\cdot|$ is an ultrametric absolute value function on $k$,
then
\begin{equation}
\label{|n cdot 1| le 1}
        |n \cdot 1| \le 1
\end{equation}
for every $n \in {\bf Z}_+$, so that $|\cdot|$ is non-archimedian on
$k$.  One can check that every non-archimedian absolute value function
on $k$ satisfies (\ref{|n cdot 1| le 1}) for every $n \in {\bf Z}_+$,
using (\ref{|x y| = |x| |y|}).  Equivalently, if $|n_0 \cdot 1| > 1$
for some $n_0 \in {\bf Z}_+$, then one can get positive integers $n$
such that $|n \cdot 1|$ is arbitrarily large, by taking powers of
$n_0$.  If an absolute value function $|\cdot|$ on $k$ satisfies
(\ref{|n cdot 1| le 1}) for every $n \in {\bf Z}_+$, then $|\cdot|$ is
an ultrametric absolute value function on $k$, as in Lemma 1.5 on p16
of \cite{c}, or Theorem 2.2.2 on p28 of \cite{fg}.  Thus
non-archimedian absolute value functions are in fact ultrametric
absolute value functions.

        If $k$ has positive characteristic, then there are only
finitely many elements of $k$ of the form $n \cdot 1$ for some $n \in
{\bf Z}_+$.  This implies that every absolute value function on $k$ is
non-archimedian.  If $k$ has only finitely many elements, then every
$x \in k$ with $x \ne 0$ satisfies $x^n = 1$ for some $n \in {\bf
  Z}_+$.  In this case, the only absolute value function on $k$ is the
trivial absolute value function.

\section{Seminorms, continued}
\label{seminorms, continued}
\setcounter{equation}{0}

        Let $k$ be a field with an absolute value function $|\cdot|$,
and let $V$ be a vector space over $k$.  A nonnegative real-valued function
$N$ on $V$ is said to be a \emph{seminorm}\index{seminorms} on $V$
with respect to $|\cdot|$ on $k$ if
\begin{equation}
\label{N(t v) = |t| N(v), 2}
        N(t \, v) = |t| \, N(v)
\end{equation}
for every $v \in V$ and $t \in k$, and
\begin{equation}
\label{N(v + w) le N(v) + N(w), 2}
        N(v + w) \le N(v) + N(w)
\end{equation}
for every $v, w \in V$.  Of course, this is the same as the definition
in Section \ref{seminorms} when $k = {\bf R}$ or ${\bf C}$, with the
standard absolute value function.  As before, (\ref{N(t v) = |t| N(v),
  2}) implies that $N(0) = 0$, by taking $t = 0$, and
\begin{equation}
\label{d(v, w) = N(v - w), 2}
        d(v, w) = N(v - w)
\end{equation}
defines a semimetric on $V$.  If
\begin{equation}
\label{N(v) > 0, 2}
        N(v) > 0
\end{equation}
for every $v \in V$ with $v \ne 0$, then $N$ is said to be a
\emph{norm}\index{norms} on $V$, in which case (\ref{d(v, w) = N(v -
  w), 2}) defines a metric on $V$.

        Suppose for the moment that $|\cdot|$ is an ultrametric
absolute value function on $k$.  If $N$ is a nonnegative real-valued
function on $V$ that satisfies (\ref{N(t v) = |t| N(v), 2}) and
\begin{equation}
\label{N(v + w) le max(N(v), N(w))}
        N(v + w) \le \max(N(v), N(w))
\end{equation}
for every $v, w \in V$, then $N$ is said to be a
\emph{semi-ultranorm}\index{semi-ultranorms} on $V$.  Note that
(\ref{N(v + w) le max(N(v), N(w))}) implies (\ref{N(v + w) le N(v) +
  N(w), 2}), so that a semi-ultranorm on $V$ is a seminorm.  If $N$ is
a semi-ultranorm on $V$, then (\ref{d(v, w) = N(v - w), 2}) is a
semi-ultrametric on $V$.  A semi-ultranorm $N$ on $V$ that satisfies
(\ref{N(v) > 0, 2}) is an \emph{ultranorm}\index{ultranorms} on $V$,
which corresponds to an ultrametric on $V$ as in (\ref{d(v, w) = N(v -
  w), 2}).  If $|\cdot|$ is the trivial absolute value function on
$k$, then the \emph{trivial ultranorm}\index{trivial ultranorm} is
defined on $V$ is defined on $V$ by putting $N(v) = 1$ when $v \ne 0$
and $N(0) = 0$.  It is easy to see that this is an ultranorm on $V$,
for which the corresponding ultrametric is the discrete metric.

        If $N_1, \ldots, N_l$ are finitely many seminorms on $V$
with respect to $|\cdot|$ on $k$, then their maximum and sum are
seminorms on $V$ too, as in Section \ref{seminorms}.  The maximum of
finitely many semi-ultranorms is a semi-ultranorm as well.  A
collection $\mathcal{N}$ of seminorms on $V$ leads to a collection
$\mathcal{M}$ of semimetrics on $V$ as in (\ref{d(v, w) = N(v - w),
  2}), and $\mathcal{M}$ determines a topology on $V$ as in Section
\ref{collections of semimetrics}.  Let us say that $\mathcal{N}$ is
\emph{nondegenerate}\index{nondegenerate collections of seminorms} on
$V$ if for each $v \in V$ with $v \ne 0$ there is an $N \in
\mathcal{N}$ such that (\ref{N(v) > 0, 2}) holds, as in Section
\ref{seminorms}.  As before, this implies that $\mathcal{M}$ is a
nondegenerate collection of semimetrics on $V$, so that the associated
topology on $V$ is Hausdorff.  The restriction of a seminorm $N$ on
$V$ to a linear subspace $W$ of $V$ is a seminorm on $W$, which is a
semi-ultranorm on $W$ when $N$ is a semi-ultranorm on $V$.  Thus a
collection $\mathcal{N}$ of seminorms on $V$ leads to a collection
$\mathcal{N}_W$ of seminorms on $W$ by restriction, for which the
corresponding collection $\mathcal{M}_W$ of semimetrics on $W$ is the
same as the collection of restrictions to $W$ of the semimetrics on
$V$ in the collection $\mathcal{M}$ associated to $\mathcal{N}$.  As
in Section \ref{seminorms}, the topology on $W$ associated to
$\mathcal{N}_W$ is the same as the topology on $W$ induced by the
topology on $V$ associated to $\mathcal{N}$ as before.

        Suppose that $N$ is a semi-ultranorm on $V$.  If $v, w \in V$
satisfy
\begin{equation}
\label{N(v - w) < N(v)}
        N(v - w) < N(v),
\end{equation}
then
\begin{equation}
\label{N(v) = N(w)}
        N(v) = N(w).
\end{equation}
This follows from (\ref{d(x, y) = d(x, z)}), where $d(\cdot, \cdot)$
is as in (\ref{d(v, w) = N(v - w), 2}), and with $x = 0$, $y = v$, and
$z = w$.  This is also analogous to (\ref{|y| = |z|}), and can be
verified more directly in this case, as before.

\section{$q$-Semimetrics}
\label{q-semimetrics}
\setcounter{equation}{0}

        Let $X$ be a set, and let $q$ be a positive real number.
A nonnegative real-valued function $d(x, y)$ defined for $x, y \in X$
is said to be a
\emph{$q$-semimetric}\index{qsemimetrics@$q$-semimetrics} on $X$ if
$d(x, y)^q$ is a semimetric on $X$.  Equivalently, this means that
$d(x, y)$ satisfies (\ref{d(x, x) = 0}), (\ref{d(x, y) = d(y, x)}),
and
\begin{equation}
\label{d(x, z)^q le d(x, y)^q + d(y, z)^q}
        d(x, z)^q \le d(x, y)^q + d(y, z)^q
\end{equation}
for every $x, y, z \in X$, instead of the usual triangle inequality
(\ref{d(x, z) le d(x, y) + d(y, z)}).  Similarly, if $d(x, y)^q$ is a
metric on $X$, then $d(x, y)$ is said to be a
\emph{$q$-metric}\index{qmetrics@$q$-metrics} on $X$.  Note that
(\ref{d(x, z)^q le d(x, y)^q + d(y, z)^q}) is the same as saying that
\begin{equation}
\label{d(x, z) le (d(x, y)^q + d(y, z)^q)^{1/q}}
        d(x, z) \le (d(x, y)^q + d(y, z)^q)^{1/q}
\end{equation}
for every $x, y, z \in X$.

        One can define open and closed balls in $X$ with respect to
a $q$-semimetric $d(x, y)$ on $X$ in exactly the same way as for an
ordinary semimetric, as in (\ref{B(x, r) = B_d(x, r) = {y in X : d(x,
    y) < r}}) and (\ref{overline{B}(x, r) = overline{B}_d(x, r) = {y
    in X : d(x, y) le r}}).  Observe that
\begin{equation}
\label{B_d(x, r) = B_{d^q}(x, r^q)}
        B_d(x, r) = B_{d^q}(x, r^q)
\end{equation}
for every $x \in X$ and $r > 0$, and that
\begin{equation}
\label{overline{B}_d(x, r) = overline{B}_{d^q}(x, r^q)}
        \overline{B}_d(x, r) = \overline{B}_{d^q}(x, r^q)
\end{equation}
for every $x \in X$ and $r \ge 0$.  Using these open balls in $X$ with
respect to $d(\cdot, \cdot)$, one can define a topology on $X$
associated to $d(\cdot, \cdot)$ in the same way as for ordinary
semimetrics.  The topology on $X$ associated to $d(\cdot, \cdot)$ is
the same as the topology on $X$ associated to $d(\cdot, \cdot)^q$,
because of (\ref{B_d(x, r) = B_{d^q}(x, r^q)}).  This implies that the
topology associated to a $q$-semimetric has the same properties as the
topology associated to an ordinary semimetric.  In particular, open
balls with respect to $d(\cdot, \cdot)$ are open sets with respect to
the topology on $X$ associated to $d(\cdot, \cdot)$, because of
(\ref{B_d(x, r) = B_{d^q}(x, r^q)}) and the analogous statement for
$d(\cdot, \cdot)^q$.  Similarly, closed balls with respect to
$d(\cdot, \cdot)$ are closed sets with respect to the topology on $X$
associated to $d(\cdot, \cdot)$, because of (\ref{overline{B}_d(x, r)
  = overline{B}_{d^q}(x, r^q)}) and the analogous statement for
$d(\cdot, \cdot)$.  If $Y \subseteq X$, then the restriction of
$d(x, y)$ to $x, y \in Y$ defines a $q$-semimetric on $Y$, with the
same types of properties as in the case of ordinary semimetrics.

        Let $\mathcal{M}$ be a  collection of $q$-semimetrics on $X$,
and let $\widetilde{\mathcal{M}}$ be the corresponding collection of
ordinary semimetrics on $X$ obtained by taking the $q$th powers of the
elements of $\mathcal{M}$.  More precisely, one can allow $q > 0$ to
depend on the element of $\mathcal{M}$.  One can define a topology on
$X$ associated to $\mathcal{M}$ in the same way as in Section
\ref{collections of semimetrics}.  This is the same as the topology on
$X$ associated to $\widetilde{\mathcal{M}}$ as before, because of
(\ref{B_d(x, r) = B_{d^q}(x, r^q)}), and so the topology on $X$
associated to $\mathcal{M}$ has the same types of properties as
before.  One can also define nondegeneracy of $\mathcal{M}$ on $X$ in
the same way as before, which is equivalent to nondegeneracy of
$\widetilde{\mathcal{M}}$
on $X$.\index{nondegenerate collections of semimetrics}

        Let $a$, $b$ be nonnegative real numbers, and observe that
\begin{equation}
\label{max(a, b) le (a^q + b^q)^{1/q} le 2^{1/q} max(a, b)}
        \max(a, b) \le (a^q + b^q)^{1/q} \le 2^{1/q} \, \max(a, b).
\end{equation}
If $0 < q \le r < \infty$, then it follows that
\begin{equation}
\label{a^r + b^r le ... = (a^q + b^q)^{r/q}}
 \quad  a^r + b^r \le \max(a, b)^{r - q} \, (a^q + b^q)
                  \le (a^q + b^q)^{(r - q)/q + 1} = (a^q + b^q)^{r/q},
\end{equation}
so that
\begin{equation}
\label{(a^r + b^r)^{1/r} le (a^q + b^q)^{1/q}}
        (a^r + b^r)^{1/r} \le (a^q + b^q)^{1/q}.
\end{equation}
This implies that every $r$-semimetric on $X$ is also a $q$-semimetric
on $X$ when $0 < q \le r < \infty$, using (\ref{d(x, z) le (d(x, y)^q
  + d(y, z)^q)^{1/q}}).  Similarly, every semi-ultrametric on $X$ is a
$q$-semimetric on $X$ for each $q > 0$, by the first inequality in
(\ref{max(a, b) le (a^q + b^q)^{1/q} le 2^{1/q} max(a, b)}).
Semi-ultrametrics on $X$ may be considered as $q$-semimetrics on $X$
with $q = \infty$, since
\begin{equation}
\label{lim_{q to infty} (a^q + b^q)^{1/q} = max(a, b)}
        \lim_{q \to \infty} (a^q + b^q)^{1/q} = \max(a, b)
\end{equation}
for every $a, b \ge 0$, by (\ref{max(a, b) le (a^q + b^q)^{1/q} le
  2^{1/q} max(a, b)}).

\section{$q$-Absolute value functions}
\label{q-absolute value functions}
\setcounter{equation}{0}

        Let $k$ be a field, and let $q$ be a positive real number.
A nonnegative real-valued function $|x|$ on $k$ is said to be a
\emph{$q$-absolute value function}\index{qabsolute value
  functions@$q$-absolute value functions} on $k$ if it satisfies
(\ref{|x| = 0 if and only if x = 0}), (\ref{|x y| = |x| |y|}), and
\begin{equation}
\label{|x + y|^q le |x|^q + |y|^q}
        |x + y|^q \le |x|^q + |y|^q
\end{equation}
for every $x, y \in k$, instead of (\ref{|x + y| le |x| + |y|}).
Equivalently, $|x|$ is a $q$-absolute value function on $k$ if and
only if $|x|^q$ is an ordinary absolute value function on $k$.  In
particular, the previous notion of an absolute value function is the
same as a $q$-absolute value function with $q = 1$.  If $|\cdot|$ is a
$q$-absolute value function on $k$, then
\begin{equation}
\label{d(x, y) = |x - y|, 2}
        d(x, y) = |x - y|
\end{equation}
defines a $q$-metric on $k$.

        As before, (\ref{|x + y|^q le |x|^q + |y|^q}) is the same as saying
that
\begin{equation}
\label{|x + y| le (|x|^q + |y|^q)^{1/q}}
        |x + y| \le (|x|^q + |y|^q)^{1/q}
\end{equation}
for every $x, y \in k$.  If $0 < q \le r < \infty$ and $|\cdot|$ is an
$r$-absolute value function on $k$, then it is easy to see that
$|\cdot|$ is a $q$-absolute value function on $k$ too, using
(\ref{(a^r + b^r)^{1/r} le (a^q + b^q)^{1/q}}) and the reformulation
(\ref{|x + y| le (|x|^q + |y|^q)^{1/q}}) of (\ref{|x + y|^q le |x|^q +
  |y|^q}).  Similarly, an ultrametric absolute value function on $k$
is a $q$-absolute value function on $k$ for every $q > 0$, because of
the first inequality in (\ref{max(a, b) le (a^q + b^q)^{1/q} le
  2^{1/q} max(a, b)}).  As in the previous section, ultrametric
absolute value functions may be considered as $q$-absolute value
functions with $q = \infty$, because of (\ref{lim_{q to infty} (a^q +
  b^q)^{1/q} = max(a, b)}).

        Of course, the archimedian\index{archimedian absolute value functions}
and non-archimedian\index{non-archimedian absolute value functions}
properties can be defined for $q$-absolute value functions in the same
way as for ordinary absolute value functions, as in Section
\ref{absolute value functions}.  Equivalently, a $q$-absolute value
function $|\cdot|$ on $k$ is archimedian or non-archimedian exactly
when $|\cdot|^q$ has the corresponding property as an ordinary
absolute value function on $k$.  If $|\cdot|$ is a non-archimedian
$q$-absolute value function on $k$ for some $q > 0$, so that
$|\cdot|^q$ is non-archimedian as an ordinary absolute value function
on $k$, then $|\cdot|^q$ is an ultrametric absolute value function on
$k$, as mentioned in Section \ref{absolute value functions}.  This
implies that $|\cdot|$ is an ultrametric absolute value function on
$k$ as well.

        Let $q_1$, $q_2$ be positive real numbers, and suppose that
$|\cdot|_1$, $|\cdot|_2$ are $q_1$, $q_2$-absolue value functions on $k$,
respectively.  Let us say that $|\cdot|_1$, $|\cdot|_2$ are 
\emph{equivalent}\index{equivalent absolute value functions} on $k$
if there is a positive real number $a$ such that
\begin{equation}
\label{|x|_2 = |x|_1^a}
        |x|_2 = |x|_1^a
\end{equation}
for every $x \in k$.  Of course, this implies that the corresponding
$q_1$, $q_2$-metrics on $k$ as in (\ref{d(x, y) = |x - y|, 2}) satisfy
an analogous relation, and hence that they determine the same topology
on $k$.  Conversely, if the corresponding $q_1$, $q_2$-metrics on $k$
determine the same topology on $k$, then $|\cdot|_1$, $|\cdot|_2$ are
equivalent in this sense.  This follows from Lemma 3.2 on p20 of
\cite{c} or Lemma 3.1.2 on p42 of \cite{fg} when $q_1 = q_2 = 1$, and
otherwise one can reduce to this case by considering
$|\cdot|_1^{q_1}$, $|\cdot|_2^{q_2}$ instead of $|\cdot|_1$,
$|\cdot|_2$.

        A famous theorem of Ostrowski\index{Ostrowski's theorem}
implies that every absolute value function on ${\bf Q}$ is either
trivial, equivalent to the standard Euclidean absolute value function
on ${\bf Q}$, or equivalent to the $p$-adic absolute value function on
${\bf Q}$ for some prime number $p$.  See Theorem 2.1 on p16 of
\cite{c}, or Theorem 3.1.3 on p44 of \cite{fg}.  If $|\cdot|$ is a
$q$-absolute value function on ${\bf Q}$ for some positive real number
$q$, then the same conclusion holds, since Ostrowski's theorem can be
applied to $|\cdot|^q$ as an ordinary absolute value function on ${\bf
  Q}$.

        If $k$ is a field of characteristic $0$, then it is well known
that there is a natural embedding of ${\bf Q}$ into $k$.  If $|\cdot|$
is a $q$-absolute value function on $k$ for some $q > 0$, then
$|\cdot|$ induces a $q$-absolute value function on ${\bf Q}$, using
this embedding.  Note that $|\cdot|$ is archimedian or non-archimedian
on $k$ exactly when the induced $q$-absolute value function on ${\bf Q}$
has the same property.  In particular, if $|\cdot|$ is archimedian on
$k$, then the induced $q$-absolute value function is archimedian on
${\bf Q}$, and hence is equivalent to the standard Euclidean absolute
value function on ${\bf Q}$, by Ostrowski's theorem.

        Let $k$ be any field again, and let $|\cdot|$ be a nonnegative
real-valued function on $k$.  If $|\cdot|$ satisfies (\ref{|x| = 0 if
  and only if x = 0}) and (\ref{|x y| = |x| |y|}), then $|\cdot|^a$
has the same properties for every $a > 0$.  If $|\cdot|$ is a
$q$-absolute value function on $k$ for some $q > 0$, then $|\cdot|^a$
is a $(q/a)$-absolute value function on $k$.  Similarly, if $|\cdot|$
is an ultrametric absolute value function on $k$, then $|\cdot|^a$ is
an ultrametric absolute value function on $k$ for every $a > 0$, which
was implicitly mentioned earlier.  If $|\cdot|$ is the standard
Euclidean absolute value function on ${\bf Q}$, then $|\cdot|^a$ is
not an absolute value function on ${\bf Q}$ for any $a > 1$, which is
the same as saying that $|\cdot|$ is not a $q$-absolute value function
on ${\bf Q}$ for any $q > 1$.

\section{$q$-Seminorms}
\label{q-seminorms}
\setcounter{equation}{0}

        Let $k$ be a field, and let $|\cdot|$ be a $q$-absolute value
function on $k$ for some $q > 0$.  Also let $V$ be a vector space over
$k$, and let $N$ be a nonnegative real-valued function on $V$.  Let us
say that $N$ is a \emph{$q$-seminorm}\index{qseminorms@$q$-seminorms}
on $V$ with respect to $|\cdot|$ on $k$ if $N$ satisfies the
homogeneity condition (\ref{N(t v) = |t| N(v), 2}), and if
\begin{equation}
\label{N(v + w)^q le N(v)^q + N(w)^q}
        N(v + w)^q \le N(v)^q + N(w)^q
\end{equation}
for every $v, w \in V$.  If $q = 1$, then $|\cdot|$ is an ordinary
absolute value function on $k$, and this is the same as saying that
$N$ is an ordinary seminorm on $V$ with respect to $|\cdot|$ on $k$,
as in Section \ref{seminorms, continued}.  If $q$ is any positive real
number, then $|\cdot|$ is a $q$-absolute value function on $k$ if and
only if $|\cdot|^q$ is an ordinary absolute value function on $k$, in
which case $N$ is a $q$-seminorm on $V$ with respect to $|\cdot|$ on
$k$ if and only if $N(v)^q$ is an ordinary seminorm on $V$ with
respect to $|\cdot|^q$ on $k$.  If $N$ is a $q$-seminorm on $V$ with
respect to $|\cdot|$ on $k$ for any $q > 0$, then
\begin{equation}
\label{d(v, w) = N(v - w), 3}
        d(v, w) = N(v - w)
\end{equation}
defines a $q$-semimetric on $V$.  If we also have that $N(v) > 0$ for
every $v \in V$ with $v \ne 0$, then $N$ is said to be a
\emph{$q$-norm}\index{qnorms@$q$-norms} on $V$.  In this case,
(\ref{d(v, w) = N(v - w), 3}) defines a $q$-metric on $V$.

        As usual, (\ref{N(v + w)^q le N(v)^q + N(w)^q}) is the same
as saying that
\begin{equation}
\label{N(v + w) le (N(v)^q + N(w)^q)^{1/q}}
        N(v + w) \le (N(v)^q + N(w)^q)^{1/q}
\end{equation}
for every $v, w \in V$.  Suppose that $|\cdot|$ is an $r$-absolute
value function on $k$ for some $r > 0$, and that $N$ is an
$r$-seminorm with respect to $|\cdot|$ on $k$.  If $0 < q \le r$, then
$|\cdot|$ is a $q$-absolute value function on $k$ too, as in the
previous section.  It is easy to see that $N$ is a $q$-seminorm on $V$
with respect to $|\cdot|$ on $k$ under these conditions, using
(\ref{(a^r + b^r)^{1/r} le (a^q + b^q)^{1/q}}), and the reformulation
(\ref{N(v + w) le (N(v)^q + N(w)^q)^{1/q}}) of (\ref{N(v + w)^q le
  N(v)^q + N(w)^q}).  Similarly, if $|\cdot|$ is an ultrametric
absolute value function on $k$, then $|\cdot|$ is a $q$-absolute value
function on $k$ for every $q > 0$.  If $N$ is a semi-ultranorm on $V$
with respect to $|\cdot|$ on $k$, then it is easy to see that $N$ is a
$q$-seminorm on $V$ with respect to $|\cdot|$ on $k$ for every $q >
0$, using the first inequality in (\ref{max(a, b) le (a^q + b^q)^{1/q}
  le 2^{1/q} max(a, b)}).  As before, semi-ultranorms may be
considered as $q$-seminorms with $q = \infty$, because of (\ref{lim_{q
    to infty} (a^q + b^q)^{1/q} = max(a, b)}).

        Let $|\cdot|$ be a $q$-absolute value function on $k$ for
some $q > 0$ again, and hence for some range of $q$'s, which may
include $q = \infty$.  Also let $\mathcal{N}$ be a collection of
$q$-seminorms on $V$ with respect to $|\cdot|$ on $k$, where $q$ is
allowed to depend on the element of $\mathcal{N}$, as long as
$|\cdot|$ is a $q$-absolute value function on $k$.  This leads to a
collection $\mathcal{M}$ of semimetrics on $V$, associated to the
elements of $\mathcal{N}$ as in (\ref{d(v, w) = N(v - w), 3}), and
where $q$ is allowed to depend on the element of $\mathcal{M}$.  Using
$\mathcal{M}$, we can get a topology on $V$, as in Sections
\ref{collections of semimetrics} and \ref{q-semimetrics}.  As before,
$\mathcal{N}$ is said to be
\emph{nondegenerate}\index{nondegenerate collections of seminorms} on
$V$ if for each $v \in V$ with $v \ne 0$ there is an $N \in
\mathcal{N}$ such that $N(v) > 0$, in which case $\mathcal{M}$ is
nondegenerate on $V$ too.

        Let $|\cdot|$ be a nonnegative real-valued function on $k$,
let $N$ be a nonnegative real-valued function on $V$, and suppose that
$N(v) > 0$ for some $v \in V$.  If $N$ satisfies the homogeneity
condition (\ref{N(t v) = |t| N(v), 2}), then $|\cdot|$ has to satisfy
the multiplicative property (\ref{|x y| = |x| |y|}).  One can also use
this to get that $|x| > 0$ for every $x \in k$ with $x \ne 0$.  If $N$
satisfies (\ref{N(v + w)^q le N(v)^q + N(w)^q}) for some $q > 0$ too,
then the homogeneity condition (\ref{N(t v) = |t| N(v), 2}) implies
that $|\cdot|$ satisfies (\ref{|x + y|^q le |x|^q + |y|^q}) in the
previous section.  Similarly, if $N$ satisfies the ultrametric version
of the triangle inequality (\ref{N(v + w) le max(N(v), N(w))}), then
the homogeneity condition (\ref{N(t v) = |t| N(v), 2}) implies that
$|\cdot|$ has the analogous property (\ref{|x + y| le max(|x|, |y|)}).

\section{Balanced sets}
\label{balanced sets}
\setcounter{equation}{0}

        Let $k$ be a field with a $q$-absolute value function $|\cdot|$
for some $q > 0$, and let $V$ be a vector space over $k$.  If $E
\subseteq V$ and $t \in k$, then put
\begin{equation}
\label{t E = {t v : v in E}}
        t \, E = \{t \, v : v \in E\}
\end{equation}
and
\begin{equation}
\label{-E = (-1) E = {-v : v in E}}
        -E = (-1) \, E = \{-v : v \in E\}.
\end{equation}
A set $E \subseteq V$ is said to be \emph{balanced}\index{balanced sets}
with respect to $|\cdot|$ on $k$ if
\begin{equation}
\label{t E subseteq E}
        t \, E \subseteq E
\end{equation}
for every $t \in k$ with $|t| \le 1$.  In particular, this implies
that
\begin{equation}
\label{t E = E}
        t \, E = E
\end{equation}
for every $t \in k$ with $|t| = 1$, since (\ref{t E subseteq E}) holds
with $t$ replaced by $1/t$ when $|t| = 1$.  Sometimes balanced sets
are said to be \emph{circled},\index{circled sets} although this could
also be used to refer to (\ref{t E = E}), especially when $k = {\bf
  C}$ with the standard absolute value function.  Similarly, $E
\subseteq V$ is said to be \emph{symmetric}\index{symmetric sets}
(about the origin) if (\ref{t E subseteq E}) holds with $t = -1$,
which is the same as saying that (\ref{t E = E}) holds with $t = -1$.
Note that nonempty balanced sets automatically contain $0$, by taking
$t = 0$ in (\ref{t E subseteq E}).  If $|\cdot|$ is the trivial
absolute value function on $k$, then $E \subseteq V$ is balanced if
and only if either $E = \emptyset$ or $0 \in E$ and $E$ satisfies
(\ref{t E = E}) for every $t \in k$ with $t \ne 0$.  If $|\cdot|$ is
any $q$-absolute value function on $k$, then the union and
intersection of any collection of balanced subsets of $V$ are balanced
in $V$ too.  If $E$ is any subset of $V$, then
\begin{equation}
\label{bigcup {t E : t in k, |t| le 1}}
        \bigcup \{t \, E : t \in k, \ |t| \le 1\}
\end{equation}
is a balanced subset of $V$ that contains $E$.  This is the smallest
balanced subset of $V$ that contains $E$, which may be described as
the \emph{balanced hull}\index{balanced hull} of $E$ in $V$.

        Suppose for the moment that $k = {\bf R}$ or ${\bf C}$, with
the standard absolute value function.  A set $E \subseteq V$ is said
to be \emph{starlike about $0$}\index{starlike sets} if (\ref{t E
  subseteq E}) holds for every $t \in {\bf R}$ with $0 \le t \le 1$.
As before, this implies that $0 \in E$ when $E \ne \emptyset$, by
taking $t = 0$.  If $k = {\bf R}$, then $E \subseteq V$ is balanced if
and only if $E$ is both symmetric and starlike about $0$.  If $k =
{\bf C}$, then $E$ is balanced if and only if $E$ is starlike about
$0$ and $E$ satisfies (\ref{t E = E}) for every $t \in {\bf C}$ with
$|t| = 1$.  In both cases, the union and intersection of any
collection of starlike subsets of $V$ about $0$ are starlike in $V$
about $0$ as well.  In particular, for any $E \subseteq V$,
\begin{equation}
\label{bigcup {t E : t in {bf R}, 0 le t le 1}}
        \bigcup \{t \, E : t \in {\bf R}, \ 0 \le t \le 1\}
\end{equation}
is starlike about $0$ in $V$ and contains $E$.  As before, this is the
smallest subset of $V$ that is starlike about $0$ and contains $E$,
and which may be described as the \emph{starlike hull}\index{starlike
  hull} of $E$ in $V$ about $0$.  If $k = {\bf R}$ and $E$ is
symmetric about $0$, then the balanced hull of $E$ in $V$ is the same
as the starlike hull of $E$ in $V$ about $0$.  Similarly, if $k = {\bf
  C}$ and $E$ satisfies (\ref{t E = E}) for every $t \in {\bf C}$ with
$|t| = 1$, then the balanced hull of $E$ in $V$ is the same as the
starlike hull of $E$ in $V$ about $0$.

        Let $|\cdot|$ be a $q$-absolute value function on any field $k$
again, and let $N$ be a nonnegative real-valued function on $V$ that
satisfies the usual homogeneity condition (\ref{N(t v) = |t| N(v), 2})
with respect to $|\cdot|$ on $k$.  Observe that
\begin{equation}
\label{B_N(0, r) = {v in V : N(v) < r}}
        B_N(0, r) = \{v \in V : N(v) < r\}
\end{equation}
is a balanced subset of $V$ for each $r > 0$, and that
\begin{equation}
\label{overline{B}_N(0, r) = {v in V : N(v) le r}}
        \overline{B}_N(0, r) = \{v \in V : N(v) \le r\}
\end{equation}
is balanced for every $r \ge 0$. More precisely,
\begin{equation}
\label{t B_N(0, r) = B_N(0, |t| r)}
        t \, B_N(0, r) = B_N(0, |t| \, r)
\end{equation}
for every $r > 0$ and $t \in k$ with $t \ne 0$, and
\begin{equation}
\label{t overline{B}_N(0, r) = overline{B}_N(0, |t| r)}
        t \, \overline{B}_N(0, r) = \overline{B}_N(0, |t| \, r)
\end{equation}
for every $r \ge 0$ and $t \in k$ with $t \ne 0$.  If $t = 0$, then
(\ref{t overline{B}_N(0, r) = overline{B}_N(0, |t| r)}) may not hold,
since there may be $v \in V$ with $v \ne 0$ and $N(v) = 0$.  

        Let us return now to the case where $k = {\bf R}$ or ${\bf C}$
with the standard absolute value function, and let $N$ be a
nonnegative real-valued function on $V$.  Suppose that $N$ is
homogeneous of degree $1$ with respect to multiplication by
nonnegative real numbers, which is to say that
\begin{equation}
\label{N(t v) = t N(v)}
        N(t \, v) = t \, N(v)
\end{equation}
for every $v \in V$ and $t \ge 0$.  As before, (\ref{B_N(0, r) = {v in
    V : N(v) < r}}) is starlike about $0$ in $V$ for every $r > 0$,
and (\ref{overline{B}_N(0, r) = {v in V : N(v) le r}}) is starlike
about $0$ in $V$ for every $r \ge 0$.  Similarly, (\ref{t B_N(0, r) =
  B_N(0, |t| r)}) holds for every $r, t > 0$, and (\ref{t
  overline{B}_N(0, r) = overline{B}_N(0, |t| r)}) holds for every $r
\ge 0$ and $t > 0$, where $|t|$ reduces to $t$ on the right sides of
these equations.  If $k = {\bf R}$, then $N$ satisfies (\ref{N(t v) =
  |t| N(v), 2}) exactly when $N$ satisfies (\ref{N(t v) = t N(v)}) and
$N$ is invariant under multiplication by $-1$.  If $k = {\bf C}$, then
$N$ satisfies (\ref{N(t v) = |t| N(v), 2}) exactly when $N$ satisfies
(\ref{N(t v) = t N(v)}) and $N$ is invariant under multiplication by
complex numbers with absolute value equal to $1$.  In both cases, we
are back to the situation discussed in the preceding paragraph.

\section{Absorbing sets}
\label{absorbing sets}
\setcounter{equation}{0}

        Let $k$ be a field with a $q$-absolute value function $|\cdot|$
for some $q > 0$ again, and let $V$ be a vector space over $k$.  A set
$A \subseteq V$ is said to be \emph{absorbing}\index{absorbing sets}
with respect to $|\cdot|$ on $k$ if for each $v \in V$ there is a
$t_0(v) \in k$ such that $t_0(v) \ne 0$ and
\begin{equation}
\label{t v in A}
        t \, v \in A
\end{equation}
for every $t \in k$ with $|t| \le |t_0(v)|$.  Equivalently, $A
\subseteq V$ is absorbing if for each $v \in V$ there is a $t_1(v) \in
k$ such that
\begin{equation}
\label{v in t A}
        v \in t \, A
\end{equation}
for every $t \in k$ with $|t| \ge |t_1(v)|$.  More precisely, both
versions of the absorbing condition imply that $0 \in A$.  If $A$
satisfies the first version of the absorbing condition, then $A$
satisfies the second version of the absorbing condition with $t_1(v) =
1/t_0(v)$.  Conversely, if $A$ satisfies the second version of the
absorbing condition, then $A$ satisfies the first version of the
absorbing condition with $t_0(v) = 1/t_1(v)$ when $t_1(v) \ne 0$.
This uses the fact that $0 \in A$ to get that (\ref{t v in A}) holds
when $t = 0$.  If $t_1(v) = 0$, then (\ref{v in t A}) holds for every
$t \in k$, including $t = 0$.  This implies that $v = 0$, and that
(\ref{t v in A}) holds for every $t \in k$, so that one can take
$t_0(v)$ to be any nonzero element of $k$, such as the multiplicative
identity element $1$.

        An absorbing set $A \subseteq V$ is also said to be
\emph{radial}\index{radial sets} (at $0$).  This is especially natural
when $k = {\bf R}$ with the standard absolute value function.  In this
case, it suffices to consider only nonnegative real numbers $t$ in
(\ref{t v in A}) and (\ref{v in t A}), since analogous conditions for
negative real numbers may be obtained by considering $-v$ instead of
$v$.  If $k = {\bf C}$ with the standard absolute value function, then
one might consider the absorbing or radial property of a subset $A$ of
$V$ as a real vector space, instead of the stronger condition for $V$
as a complex vector space.  The two conditions are equivalent in some
situations, such as when $A$ is invariant under multiplication by
complex numbers with absolute value equal to $1$, and when $A$ is
convex.

        As a weaker version of the absorbing property, one can ask
that for each $v \in V$ there be a $t \in k$ such that $t \ne 0$ and
(\ref{t v in A}) holds.  This is equivalent to asking that $0 \in A$
and for each $v \in V$ there be a $t \in k$ that satisfies (\ref{v in
  t A}), for the same types of reasons as before.  If $A$ is balanced
with respect to $|\cdot|$ on $k$, then this weaker condition implies
that $A$ is absorbing.  If $k = {\bf R}$, then one might refine this
weaker condition a bit by restricting one's attention to $t \ge 0$.
If $A$ is starlike about $0$ and satisfies this refined version of the
weaker condition, then $A$ is absorbing with respect to the standard
absolute value function on ${\bf R}$.

        If $|\cdot|$ is the trivial absolute value function on $k$,
then $A \subseteq V$ is absorbing if and only if $A = V$.  Suppose now
that $|\cdot|$ is a nontrivial $q$-absolute value function on $k$.  In
this case, the first version of the absorbing property may be
reformulated as saying that for each $v \in V$ there is a positive
real number $r_0(v)$ such that (\ref{t v in A}) holds for every $t \in
k$ with $|t| \le r_0(v)$.  Similarly, the second version of the
absorbing condition can be reformulated as saying that for each $v \in
V$ there is a nonnegative real number $r_1(v)$ such that (\ref{v in t
  A}) holds for every $t \in k$ with $|t| \ge r_1(v)$.  This uses the
fact that $|t|$ can take arbitrarily large and small positive values
with $t \in k$ when $|\cdot|$ is nontrivial on $k$.  If $N$ is a
nonnegative real-valued function on $N$ that satisfies the homogeneity
condition (\ref{N(t v) = |t| N(v), 2}), then the open and closed balls
(\ref{B_N(0, r) = {v in V : N(v) < r}}) and (\ref{overline{B}_N(0, r)
  = {v in V : N(v) le r}}) centered at $0$ in $V$ with respect to $N$
are absorbing in $V$ for every $r > 0$.  If $k = {\bf R}$ with the
standard absolute value function, and if $N$ is a nonnegative
real-valued function on $V$ that is homogeneous of degree $1$ with
respect to multiplication by nonnegative real numbers, as in (\ref{N(t
  v) = t N(v)}), then the corresponding open and closed balls
(\ref{B_N(0, r) = {v in V : N(v) < r}}) and (\ref{overline{B}_N(0, r)
  = {v in V : N(v) le r}}) are absorbing in $V$ for every $r > 0$ too.
If $k = {\bf C}$ with the standard absolute value function, then one
can apply the previous statement to $V$ as a vector space over ${\bf
  R}$.

\section{Discrete absolute value functions}
\label{discrete absolute value functions}
\setcounter{equation}{0}

        Let $k$ be a field, and let $|\cdot|$ be a $q$-absolute
value function on $k$ for some $q > 0$.  Let us say that $|\cdot|$ is
\emph{discrete}\index{discrete absolute value functions} on $k$ if
there is a nonnegative real number $\rho < 1$ such that
\begin{equation}
\label{|x| le rho}
        |x| \le \rho
\end{equation}
for every $x \in k$ with $|x| < 1$.  Equivalently, if $y, z \in k$
satisfy $|y| < |z|$, then
\begin{equation}
\label{|y| le rho |z|}
         |y| \le \rho \, |z|,
\end{equation}
by applying (\ref{|x| le rho}) to $x = y/z$.  In particular, if
$|z| > 1$, then (\ref{|y| le rho |z|}) implies that
\begin{equation}
\label{|z| ge 1/rho > 1}
        |z| \ge 1/\rho > 1,
\end{equation}
by taking $y = 1$.  More precisely, if $\rho = 0$, then $1/\rho$ is
interpreted as being $+\infty$, and $|\cdot|$ is the trivial absolute
value function on $k$.

        Observe that
\begin{equation}
\label{{|x| : x in k, x ne 0}}
        \{|x| : x \in k, \, x \ne 0\}
\end{equation}
is a subgroup of the multiplicative group ${\bf R}_+$\index{R_+@${\bf
    R}_+$} of positive real numbers.  If $|\cdot|$ is discrete on $k$,
then it is easy to see that (\ref{{|x| : x in k, x ne 0}}) has no
limit points in ${\bf R}_+$ with respect to the standard Euclidean
metric on ${\bf R}$, because of (\ref{|y| le rho |z|}).  However, $0$
is a limit point of (\ref{{|x| : x in k, x ne 0}}) with respect to the
standard Euclidean metric on ${\bf R}$ when $|\cdot|$ is nontrivial on
$k$.  Conversely, if $1$ is not a limit point of (\ref{{|x| : x in k,
    x ne 0}}) with respect to the standard Euclidean metric on ${\bf
  R}_+$, then $|\cdot|$ satisfies (\ref{|x| le rho}) for some $\rho <
1$, and so $|\cdot|$ is discrete on $k$.  If $|\cdot|$ is not discrete
on $k$, then (\ref{{|x| : x in k, x ne 0}}) is dense in ${\bf R}_+$
with respect to the standard Euclidean metric.

        If $|\cdot|$ is archimedian on $k$, then $k$ has characteristic
$0$, so that there is a natural embedding of ${\bf Q}$ in $k$.  The
induced $q$-absolute value function on ${\bf Q}$ is also archimedian
under these conditions, which implies that the induced $q$-absolute
value function on ${\bf Q}$ is equivalent to the standard Euclidean
absolute value function on ${\bf Q}$, by Ostrowski's theorem.  Of
course, the set of positive values of the standard Euclidean absolute
value function on ${\bf Q}$ is the same as the set ${\bf
  Q}_+$\index{Q_+@${\bf Q}_+$} of positive rational numbers, which is
dense in ${\bf R}_+$ with respect to the standard Euclidean metric.
It follows that (\ref{{|x| : x in k, x ne 0}}) is dense in ${\bf R}_+$
too, since (\ref{{|x| : x in k, x ne 0}}) includes the positive values
of a $q$-absolute value function on ${\bf Q}$ that is equivalent to
the standard Euclidean absolute value function on ${\bf Q}$.  Thus a
discrete $q$-absolute value function on a field $k$ has to be
non-archimedian, and hence an ultrametric absolute value function.

        Let $|\cdot|$ be a nontrivial discrete $q$-absolute value
function on a field $k$, and put
\begin{equation}
\label{rho_1 = sup {|x| : x in k, |x| < 1}}
        \rho_1 = \sup \{|x| : x \in k, \ |x| < 1\},
\end{equation}
so that $0 < \rho_1 < 1$.  It is easy to see that there is an $x_1 \in
k$ such that
\begin{equation}
\label{|x_1| = rho_1}
        |x_1| = \rho_1,
\end{equation}
because $\rho_1$ is an element of the closure of (\ref{{|x| : x in k,
    x ne 0}}) in ${\bf R}_+$, and (\ref{{|x| : x in k, x ne 0}}) has
no limit points in ${\bf R}_+$, by hypothesis.  Of course,
\begin{equation}
\label{|x_1^j| = |x_1|^j = rho_1^j}
        |x_1^j| = |x_1|^j = \rho_1^j
\end{equation}
for each integer $j$, which means that (\ref{{|x| : x in k, x ne 0}})
contains all integer powers of $\rho_1$.  Note that (\ref{|x| le rho})
holds with $\rho = \rho_1$, by construction, so that (\ref{|y| le rho
  |z|}) holds with $\rho = \rho_1$ too.  Using this, one can check
that every element of (\ref{{|x| : x in k, x ne 0}}) is an integer
power of $\rho_1$ in this situation.

\section{Balanced $q$-convexity}
\label{balanced q-convexity}
\setcounter{equation}{0}

        Let $k$ be a field with a $q$-absolute value function $|\cdot|$
for some $q > 0$, and let $V$ be a vector space over $k$.  Let us say
that a balanced set $E \subseteq V$ is \emph{$q$-convex}\index{qconvex
  sets@$q$-convex sets}\index{balanced sets!  qconvex@$q$-convex} with
respect to $|\cdot|$ on $k$ if for every $v_1, v_2 \in E$ and $t_1,
t_2 \in k$ with
\begin{equation}
\label{|t_1|^q + |t_2|^q le 1}
        |t_1|^q + |t_2|^q \le 1,
\end{equation}
we have that
\begin{equation}
\label{t_1 v_1 + t_2 v_2 in E}
        t_1 \, v_1 + t_2 \, v_2 \in E.
\end{equation}
Suppose that $N$ is a nonnegative real-valued function on $V$ that
satisfies the homogeneity condition (\ref{N(t v) = |t| N(v), 2}) with
respect to $|\cdot|$ on $k$.  If $B_N(0, r)$ and $\overline{B}_N(0,
r)$ are as in (\ref{B_N(0, r) = {v in V : N(v) < r}}) and
(\ref{overline{B}_N(0, r) = {v in V : N(v) le r}}), respectively, then
$B_N(0, r)$ is balanced in $V$ for every $r > 0$, and
$\overline{B}_N(0, r)$ is balanced in $V$ for every $r \ge 0$, as
before.  If $N$ is a $q$-seminorm on $V$ with respect to $|\cdot|$ on
$k$, then it is easy to see that $B_N(0, r)$ is $q$-convex for every
$r > 0$, and $\overline{B}_N(0, r)$ is $q$-convex for every $r \ge 0$.

        Suppose for the moment that $k = {\bf R}$ or ${\bf C}$ with the
standard absolute value function, and that $0 < q \le 1$.  Let us say
that a starlike set $E \subseteq V$ about $0$ is \emph{real
  $q$-convex}\index{real qconvex sets@real $q$-convex
  sets}\index{starlike sets! real qconvex@real $q$-convex} if for
every $v_1, v_2 \in E$ and nonnegative real numbers $t_1$, $t_2$ with
\begin{equation}
\label{t_1^q + t_2^q le 1}
        t_1^q + t_2^q \le 1,
\end{equation}
we have that (\ref{t_1 v_1 + t_2 v_2 in E}) holds.  If $E \subseteq V$
is balanced, then $E$ is starlike about $0$ in particular.  In this
case, real $q$-convexity is equivalent to $q$-convexity as defined in
the preceding paragraph.  Remember that a set $E \subseteq V$ is
convex\index{convex sets} in the ordinary sense if (\ref{t_1 v_1 + t_2
  v_2 in E}) holds for every $v_1, v_2 \in E$ and nonnegative real
numbers $t_1$, $t_2$ such that
\begin{equation}
\label{t_1 + t_2 = 1}
        t_1 + t_2 = 1.
\end{equation}
If $E \subseteq V$ is convex and $0 \in E$, then $E$ is starlike about
$0$, and real $1$-convex.  Of course, real $1$-convexity implies
ordinary convexity.

        Let $N$ be a nonnegative real-valued function on $V$ that is
homogeneous of degree $1$ with respect to multiplication by
nonnegative real numbers, as in (\ref{N(t v) = t N(v)}).  Thus $B_N(0,
r)$ and $\overline{B}_N(0, r)$ can be defined as in (\ref{B_N(0, r) =
  {v in V : N(v) < r}}) and (\ref{overline{B}_N(0, r) = {v in V : N(v)
    le r}}), respectively, $B_N(0, r)$ is starlike about $0$ in $V$
for every $r > 0$, and $\overline{B}_N(0, r)$ is starlike about $0$ in
$V$ for every $r \ge 0$.  Let us say that $N$ is
\emph{$q$-subadditive}\index{qsubadditive functions@$q$-subadditive
  functions} on $V$ if
\begin{equation}
\label{N(v + w)^q le N(v)^q + N(w)^q, 2}
        N(v + w)^q \le N(v)^q + N(w)^q
\end{equation}
for every $v, w \in V$.  In this case, $B_N(0, r)$ is real $q$-convex
in $V$ for every $r > 0$, and $\overline{B}_N(0, r)$ is real
$q$-convex for every $r \ge 0$.  Note that $N$ is a $q$-seminorm on
$V$ when $N$ is also invariant under multiplication by $-1$ in the
real case, and when $N$ is invariant under multiplication by complex
numbers with absolute value equal to $1$ in the complex case.

        Let $|\cdot|$ be a $q$-absolute value function on a field $k$
for some $q > 0$ again, and observe that (\ref{|t_1|^q + |t_2|^q le 1})
is equivalent to
\begin{equation}
\label{(|t_1|^q + |t_2|^q)^{1/q} le 1}
        (|t_1|^q + |t_2|^q)^{1/q} \le 1.
\end{equation}
If $0 < \widetilde{q} \le q$, then $|\cdot|$ is also a
$\widetilde{q}$-absolute value function on $k$, as in Section
\ref{q-absolute value functions}.  Similarly, if $E \subseteq V$ is
balanced and $q$-convex with respect to $|\cdot|$ on $k$, and if $0 <
\widetilde{q} \le q$, then $E$ is $\widetilde{q}$-convex.  This uses
the fact that the left side of (\ref{(|t_1|^q + |t_2|^q)^{1/q} le 1})
is monotonically decreasing in $q$, as in Section \ref{q-semimetrics}.
If $k = {\bf R}$ or ${\bf C}$ with the standard absolute value
function, $E \subseteq V$ is starlike about $0$ and real $q$-convex,
and $0 < \widetilde{q} \le q$, then $E$ is real $\widetilde{q}$-convex
too, for essentially the same reasons.  As usual, (\ref{N(v + w)^q le
  N(v)^q + N(w)^q, 2}) can be reformulated as saying that
\begin{equation}
\label{N(v + w) le (N(v)^q + N(w)^q)^{1/q}, 2}
        N(v + w) \le (N(v)^q + N(w)^q)^{1/q}
\end{equation}
for every $v, w \in V$.  If $0 < \widetilde{q} \le q$, then
$q$-subadditivity implies $\widetilde{q}$-subadditivity, just as for
$q$-seminorms.

        Suppose now that $|\cdot|$ is an ultrametric absolute value
function on a field $k$, and let $E$ be a balanced subset of $V$.  The
analogue of $q$-convexity with $q = \infty$ with respect to $|\cdot|$
on $k$ asks that (\ref{t_1 v_1 + t_2 v_2 in E}) hold for every $v_1,
v_2 \in E$ and $t_1, t_2 \in k$ with
\begin{equation}
\label{max(|t_1|, |t_2|) le 1}
        \max(|t_1|, |t_2|) \le 1.
\end{equation}
This is the same as saying that
\begin{equation}
\label{v_1 + v_2 in E}
        v_1 + v_2 \in E
\end{equation}
for every $v_1, v_2 \in E$, because $E$ is balanced.  If $E \ne
\emptyset$, so that $0 \in E$, then this means that $E$ is a subgroup
of $V$ with respect to addition, because $E$ is symmetric about the
origin.  Note that this property implies that $E$ is $q$-convex for
every $q > 0$.

        Let $N$ be a nonnegative real-valued function on $V$ that
satisfies the homogeneity condition (\ref{N(t v) = |t| N(v), 2}) with
respect to $|\cdot|$ on $k$ again.  If $N$ is a semi-ultranorm on $V$,
then $B_N(0, r)$ has the property described in the preceding paragraph
for every $r > 0$, and $\overline{B}_N(0, r)$ has this property for
every $r \ge 0$.  Conversely, if $\overline{B}_N(0, r)$ has this
property for every $r \ge 0$, then $N$ satisfies the ultrametric
version of the triangle inequality, as in (\ref{N(v + w) le max(N(v),
  N(w))}).  Similarly, $N$ satisfies the ultrametric version of the
triangle inequality when $B_N(0, r)$ has the property just mentioned
for every $r > 0$.

        Suppose for the moment that $|\cdot|$ is trivial on $k$.
If $0 < q < \infty$, then (\ref{|t_1|^q + |t_2|^q le 1}) holds if and
only if at least one of $t_1$ and $t_2$ is equal to $0$.  This implies
that every balanced set $E \subseteq V$ with respect to $|\cdot|$ on
$k$ is $q$-convex when $0 < q < \infty$.  In this case, $E \subseteq
V$ is balanced with respect to $|\cdot|$ on $k$ if and only if $t \, E
\subseteq E$ for every $t \in k$.  It follows that a balanced set $E
\subseteq V$ satisfies the $q = \infty$ version of $q$-convexity with
respect to $|\cdot|$ on $k$ if and only if $E$ is either empty or a
linear subspace of $V$.

        Let $|\cdot|$ be an arbitrary $q$-absolute value function
on $k$ again, for some $q > 0$.  Note that the properties of being
balanced and $q$-convex with respect to $|\cdot|$ on $k$ are preserved
by linear mappings, and by scalar multiplication in particular.
Remember that the union and intersection of any collection of balanced
subsets of $V$ are also balanced in $V$, as in Section \ref{balanced
  sets}.  The intersection of any collection of balanced $q$-convex
subsets of $V$ is $q$-convex in $V$ too.  If a collection of balanced
$q$-convex subsets of $V$ is linearly ordered by inclusion, then the
union of these sets is $q$-convex in $V$ as well.  If $k = {\bf R}$ or
${\bf C}$ with the standard absolute value function and $0 < q \le 1$,
then there are analogous statements for starlike subsets of $V$ about
$0$ and starlike sets that are real $q$-convex.  Of course, this is
also analogous to the situation for subsets of $V$ that are convex in
the ordinary sense.

\section{Minkowski functionals}
\label{minkowski functionals}
\setcounter{equation}{0}

        Let $k$ be a field with a $q$-absolute value function $|\cdot|$
for some $q > 0$ again, and let $V$ be a vector space over $k$.  Also
let $A$ be a balanced absorbing subset of $V$, and put
\begin{equation}
\label{N_A(v) = inf {|t| : t in k, v in t A}}
        N_A(v) = \inf \{|t| : t \in k, \ v \in t \, A\}
\end{equation}
for each $v \in V$.  If $|\cdot|$ is nontrivial on $k$, then $N_A(v)$
can be defined equivalently by
\begin{eqnarray}
\label{N_A(v) = inf {|t| : t in k setminus {0}, v in t A} = ...}
  N_A(v) & = & \inf \{|t| : t \in k \setminus \{0\}, \ v \in t \, A\} \\
 & = & \inf \{|t| : t \in k \setminus \{0\}, \ t^{-1} \, v \in A\} \nonumber
\end{eqnarray}
for every $v \in V$.  Otherwise, (\ref{N_A(v) = inf {|t| : t in k
    setminus {0}, v in t A} = ...}) can only differ from (\ref{N_A(v)
  = inf {|t| : t in k, v in t A}}) when $v = 0$, which is the only
case where $t = 0$ can be used in (\ref{N_A(v) = inf {|t| : t in k, v
    in t A}}).  If $|\cdot|$ is the trivial absolute value function on
$k$, then we have seen that $V$ is the only absorbing subset of
itself.  If $|\cdot|$ is the trivial absolute value function on $k$
and $A = V$, then (\ref{N_A(v) = inf {|t| : t in k, v in t A}}) is the
trivial ultranorm on $V$, while (\ref{N_A(v) = inf {|t| : t in k
    setminus {0}, v in t A} = ...}) is equal to $1$ for every $v \in
V$.  However, if $|\cdot|$ is nontrivial on $k$ and $A = V$, then both
(\ref{N_A(v) = inf {|t| : t in k, v in t A}}) and (\ref{N_A(v) = inf
  {|t| : t in k setminus {0}, v in t A} = ...}) are equal to $0$ for
every $v \in V$.

        More precisely, in order to define $N_A(v)$ as in
(\ref{N_A(v) = inf {|t| : t in k, v in t A}}), it suffices to ask that
for each $v \in V$ there be a $t \in k$ such that $v \in t \, A$, so
that the infimum in (\ref{N_A(v) = inf {|t| : t in k, v in t A}}) is
taken over a nonempty set.  Similarly, in order to define $N_A(v)$ as
in (\ref{N_A(v) = inf {|t| : t in k setminus {0}, v in t A} = ...}),
it suffices to ask that for each $v \in V$ there be a $t \in k
\setminus \{0\}$ such that $v \in t \, A$.  Of course, the second
condition holds if and only if the first condition holds and $0 \in
A$.  If $A$ has either of these properties, then the balanced hull of
$A$ in $V$ is both balanced and absorbing, as in Section
\ref{absorbing sets}.  In both cases, one can check that (\ref{N_A(v)
  = inf {|t| : t in k, v in t A}}) and (\ref{N_A(v) = inf {|t| : t in
    k setminus {0}, v in t A} = ...}) are the same for $A$ as for the
balanced hull of $A$ in $V$.  Thus we may as well ask that $A$ be
balanced and absorbing, as in the preceding paragraph.  Note that $N_A$
automatically satisfies the homogeneity condition (\ref{N(t v) = |t|
  N(v), 2}).

        Suppose for the moment that $k = {\bf R}$ with the standard
absolute value function.  Instead of (\ref{N_A(v) = inf {|t| : t in k,
    v in t A}}), it is customary to consider
\begin{equation}
\label{widetilde{N}_A(v) = inf {t in {bf R}_+ cup {0} : v in t A}}
 \widetilde{N}_A(v) = \inf \{t \in {\bf R}_+ \cup \{0\} : v \in t \, A\}
\end{equation}
for each $v \in V$.  As before, one should ask that for each $v \in V$
there be a $t \ge 0$ such that $v \in t \, A$, so that the infimum in
(\ref{widetilde{N}_A(v) = inf {t in {bf R}_+ cup {0} : v in t A}}) is
taken over a nonempty set of nonnegative real numbers.  Similarly, the
analogue of (\ref{N_A(v) = inf {|t| : t in k setminus {0}, v in t A} =
  ...}) in this case is
\begin{equation}
\label{widetilde{N}_A(v) = inf {t in {bf R}_+ : v in t A} = ...}
 \widetilde{N}_A(v) = \inf \{t \in {\bf R}_+ : v \in t \, A\}
                    = \inf \{t \in {\bf R}_+ : t^{-1} \, v \in A\}
\end{equation}
for each $v \in V$.  In this case, one should ask that for each $v \in
V$ there be a $t \in {\bf R}_+$ such that $v \in t \, A$, so that the
infimum is taken over a nonempty set of positive real numbers.  As in
the previous situation, this condition is equivalent to the
corresponding condition for (\ref{widetilde{N}_A(v) = inf {t in {bf
      R}_+ cup {0} : v in t A}}) together with the additional
requirement that $0 \in A$.  If this condition holds, then
(\ref{widetilde{N}_A(v) = inf {t in {bf R}_+ cup {0} : v in t A}}) is
equal to (\ref{widetilde{N}_A(v) = inf {t in {bf R}_+ : v in t A} =
  ...}) for every $v \in V$, for essentially the same reasons as before.

        If $A$ is absorbing in $V$, then $A$ satisfies the conditions
needed to define (\ref{widetilde{N}_A(v) = inf {t in {bf R}_+ cup {0}
    : v in t A}}) and (\ref{widetilde{N}_A(v) = inf {t in {bf R}_+ : v
    in t A} = ...}).  In the other direction, if $A$ satisfies the
condition needed to define (\ref{widetilde{N}_A(v) = inf {t in {bf
      R}_+ cup {0} : v in t A}}) for every $v \in V$, and if $A$ is
starlike about $0$, then $A$ is absorbing in $V$.  More precisely,
this uses the fact that one can define the absorbing property for
subsets of real vector spaces in terms of nonnegative real numbers, as
in Section \ref{absorbing sets}.  If $A$ satisfies either of
the conditions needed to define (\ref{widetilde{N}_A(v) = inf {t in
    {bf R}_+ cup {0} : v in t A}}) or (\ref{widetilde{N}_A(v) = inf {t
    in {bf R}_+ : v in t A} = ...}), then the starlike hull of $A$ in
$V$ about $0$ is absorbing in $V$ and starlike about $0$.  As in the
previous situation, (\ref{widetilde{N}_A(v) = inf {t in {bf R}_+ cup
    {0} : v in t A}}) and (\ref{widetilde{N}_A(v) = inf {t in {bf R}_+
    : v in t A} = ...}) are the same for $A$ as for the starlike hull
of $A$ in $V$ about $0$, and so one may as well suppose that $A$ is
absorbing in $V$ and starlike about $0$.  It is easy to see that
$\widetilde{N}_A(v)$ is homogeneous of degree $1$ with respect to
multiplication by nonnegative real numbers on $V$, as in (\ref{N(t v)
  = t N(v)}).  If $A$ is symmetric about $0$, then (\ref{N_A(v) = inf
  {|t| : t in k, v in t A}}) is the same as (\ref{widetilde{N}_A(v) =
  inf {t in {bf R}_+ cup {0} : v in t A}}), and (\ref{N_A(v) = inf
  {|t| : t in k setminus {0}, v in t A} = ...}) is the same as
(\ref{widetilde{N}_A(v) = inf {t in {bf R}_+ : v in t A} = ...}).  In
particular, $\widetilde{N}_A(v)$ is invariant under multiplication by
$-1$ on $V$ in this case.

        If $k = {\bf C}$ with the standard absolute value function,
then one can also treat $V$ as a vector space over ${\bf R}$, so that
the previous remarks can be applied.  In this case, if $A \subseteq V$
is invariant under multiplication by complex numbers with absolute
value equal to $1$, then (\ref{N_A(v) = inf {|t| : t in k, v in t A}})
is the same as (\ref{widetilde{N}_A(v) = inf {t in {bf R}_+ cup {0} :
    v in t A}}), (\ref{N_A(v) = inf {|t| : t in k setminus {0}, v in t
    A} = ...}) is the same as (\ref{widetilde{N}_A(v) = inf {t in {bf
      R}_+ : v in t A} = ...}), and the conditions under which they
are defined are equivalent.  It follows that $\widetilde{N}_A(v)$ is
invariant under multiplication by complex numbers with absolute value
equal to $1$ too, which could also be verified directly from the
definitions.  As in Section \ref{absorbing sets}, $A$ is absorbing in
$V$ as a complex vector space if and only if $A$ is absorbing in $V$
as a real vector space under these conditions.  We also have that $A$
is balanced in $V$ as a complex vector space when $A$ is starlike
about $0$ in this situation, and that the balanced hull of $A$ in $V$
as a complex vector space is the same as the starlike hull of $A$ in
$V$ about $0$.

        Let $|\cdot|$ be a $q$-absolute value function on a field $k$
again, and let $A$ be a balanced absorbing subset of $V$.  It is easy
to see that
\begin{equation}
\label{B_{N_A}(0, 1) subseteq A subseteq overline{B}_{N_A}(0, 1)}
        B_{N_A}(0, 1) \subseteq A \subseteq \overline{B}_{N_A}(0, 1),
\end{equation}
where $B_{N_A}(0, 1)$ and $\overline{B}_{N_A}(0, 1)$ are as in
(\ref{B_N(0, r) = {v in V : N(v) < r}}) and (\ref{overline{B}_N(0, r)
  = {v in V : N(v) le r}}).  More precisely, the second inclusion in
(\ref{B_{N_A}(0, 1) subseteq A subseteq overline{B}_{N_A}(0, 1)})
simply says that $N_A(v) \le 1$ when $v \in A$, which holds by the
definition (\ref{N_A(v) = inf {|t| : t in k, v in t A}}) of $N_A(v)$.
To get the first inclusion in (\ref{B_{N_A}(0, 1) subseteq A subseteq
  overline{B}_{N_A}(0, 1)}), observe that if $v \in V$ and $N_A(v) <
1$, then there is a $t \in k$ such that $v \in t \, A$ and $|t| < 1$,
by the definition (\ref{N_A(v) = inf {|t| : t in k, v in t A}}) of
$N_A(v)$.  If $A$ is balanced, then this implies that $v \in A$, as
desired.  Similarly, if $k = {\bf R}$ with the standard absolute value
function, and $A \subseteq V$ is starlike about $0$ and absorbing,
then
\begin{equation}
\label{B_{widetilde{N}_A}(0, 1) subseteq A subseteq ...}
        B_{\widetilde{N}_A}(0, 1) \subseteq A
                      \subseteq \overline{B}_{\widetilde{N}_A}(0, 1).
\end{equation}
If $k = {\bf C}$ with the standard absolute value function, then one
can also treat $V$ as a real vector space, and apply the previous
statement.

        Suppose now that $|\cdot|$ is a nontrivial discrete $q$-absolute
value function on a field $k$.  This implies that
\begin{equation}
\label{{|x| : x in k}}
        \{|x| : x \in k\}
\end{equation}
is a closed subset of the real line with respect to the standard
topology, and that $0$ is the only limit point of (\ref{{|x| : x in
    k}}) in ${\bf R}$.  If $A \subseteq V$ is balanced and absorbing,
as before, then it follows that $N_A(v)$ is an element of (\ref{{|x| :
    x in k}}) for every $v \in V$.  More precisely, if $N_A(v) > 0$,
then the infimum in (\ref{N_A(v) = inf {|t| : t in k, v in t A}}) is
attained, and so there is a $t \in k$ such that $v \in t \, A$ and
$N_A(v) = |t|$.  In this case, we have that
\begin{equation}
\label{A = overline{B}_{N_A}(0, 1)}
        A = \overline{B}_{N_A}(0, 1),
\end{equation}
instead of (\ref{B_{N_A}(0, 1) subseteq A subseteq
  overline{B}_{N_A}(0, 1)}).

        Let $|\cdot|$ be an ultrametric absolute value function on
a field $k$, and suppose that $A \subseteq V$ is balanced, absorbing,
and $\infty$-convex, as in the previous section.  We would like to
check that $N_A$ is a semi-ultranorm on $V$ under these conditions.
We already know that $N_A$ satisfies the homogeneity condition
(\ref{N(t v) = |t| N(v), 2}), and so it suffices to show that $N_A$
satisfies the ultrametric version (\ref{N(v + w) le max(N(v), N(w))})
of the triangle inequality.  Let $r > 0$ be given, and suppose that
$v_1, v_2 \in V$ satisfy
\begin{equation}
\label{N_A(v_1), N_A(v_2) < r}
        N_A(v_1), \ N_A(v_2) < r.
\end{equation}
This implies that there are $t_1, t_2 \in k$ such that $|t_1|, |t_2| <
r$ and
\begin{equation}
\label{v_1 in t_1 A, v_2 in t_2 A}
        v_1 \in t_1 \, A, \quad v_2 \in t_2 \, A,
\end{equation}
by the definition (\ref{N_A(v) = inf {|t| : t in k, v in t A}}) of
$N_A$.  Let $t$ be $t_1$ or $t_2$, in such a way that
\begin{equation}
\label{|t| = max(|t_1|, |t_2|)}
        |t| = \max(|t_1|, |t_2|).
\end{equation}
Using (\ref{v_1 in t_1 A, v_2 in t_2 A}), we get that $v_1, v_2 \in t
\, A$, since $A$ is supposed to be balanced.  It follows that
\begin{equation}
\label{v_1 + v_2 in t A}
        v_1 + v_2 \in t \, A,
\end{equation}
because $A$ is $\infty$-convex in $V$, by hypothesis.  Thus
\begin{equation}
\label{N_A(v_1 + v_2) le |t| = max(|t_1|, |t_2|) < r}
        N_A(v_1 + v_2) \le |t| = \max(|t_1|, |t_2|) < r,
\end{equation}
which implies that $N_A$ satisfies the ultrametric version of the
triangle inequality, as desired.

\section{Convexity and subadditivity}
\label{convexity, subadditivity}
\setcounter{equation}{0}

        Let $k$ be a field with a $q$-absolute value function $|\cdot|$
for some $q > 0$, and suppose that $|\cdot|$ is not discrete on $k$.
As in Section \ref{discrete absolute value functions}, this implies
that the set (\ref{{|x| : x in k, x ne 0}}) of positive values of
$|\cdot|$ on $k$ is dense in ${\bf R}_+$ with respect to the standard
topology on ${\bf R}$.  Also let $V$ be a vector space over $k$, and
let $N$ be a nonnegative real-valued function on $V$ that satisfies
the homogeneity condition (\ref{N(t v) = |t| N(v), 2}) with respect to
$|\cdot|$ on $k$.  Remember that the corresponding open and closed
balls $B_N(0, r)$ and $\overline{B}_N(0, r)$ centered at $0$ in $V$
are defined in (\ref{B_N(0, r) = {v in V : N(v) < r}}) and
(\ref{overline{B}_N(0, r) = {v in V : N(v) le r}}), and are balanced
subsets of $V$ with respect to $|\cdot|$ on $k$.  If either $B_N(0,
r)$ or $\overline{B}_N(0, r)$ is $q$-convex in $V$ with respect to
$|\cdot|$ on $k$ for any $r > 0$, then one can check that $B_N(0, r)$
and $\overline{B}_N(0, r)$ are both $q$-convex for every $r > 0$.
This uses (\ref{t B_N(0, r) = B_N(0, |t| r)}) and (\ref{t
  overline{B}_N(0, r) = overline{B}_N(0, |t| r)}) to first go from a
single positive radius to a dense set of positive radii in this
situation.  One can then use the remarks about unions and
intersections of collections of $q$-convex sets at the end of Section
\ref{balanced q-convexity} to get all positive radii, and to switch
between open and closed balls.  If every positive real number occurs
as a value of $|\cdot|$ on $k$, then one can go directly from a single
positive radius to all positive radii in the first step.

        Suppose that for every $v_1, v_2 \in V$ with
\begin{equation}
\label{N(v_1), N(v_2) < 1}
        N(v_1), \ N(v_2) < 1
\end{equation}
and $t_1, t_2 \in k$ with
\begin{equation}
\label{|t_1|^q + |t_2|^q le 1, 2}
        |t_1|^q + |t_2|^q \le 1
\end{equation}
we have that
\begin{equation}
\label{N(t_1 v_1 + t_2 v_2) le 1}
        N(t_1 \, v_1 + t_2 \, v_2) \le 1.
\end{equation}
In particular, this condition holds when either $B_N(0, 1)$ or
$\overline{B}_N(0, 1)$ is $q$-convex in $V$ with respect to $|\cdot|$
on $k$.  Let $w_1, w_2 \in V$ be given, and let us check that
\begin{equation}
\label{N(w_1 + w_2)^q le N(w_1)^q + N(w_2)^q}
        N(w_1 + w_2)^q \le N(w_1)^q + N(w_2)^q,
\end{equation}
so that $N$ is a $q$-seminorm on $V$.  If $\tau_1, \tau_2 \in k$ satisfy
\begin{equation}
\label{N(w_1) < |tau_1|, N(w_2) < |tau_2|}
        N(w_1) < |\tau_1|, \quad N(w_2) < |\tau_2|,
\end{equation}
then
\begin{equation}
\label{v_1 = tau_1^{-1} w_1, v_2 = tau_2^{-1} w_2}
        v_1 = \tau_1^{-1} \, w_1, \quad v_2 = \tau_2^{-1} \, w_2
\end{equation}
satisfy (\ref{N(v_1), N(v_2) < 1}).  If $\tau_3 \in k$ satisfies
\begin{equation}
\label{|tau_1|^q + |tau_2|^q le |tau_3|^q}
        |\tau_1|^q + |\tau_2|^q \le |\tau_3|^q,
\end{equation}
then
\begin{equation}
\label{t_1 = tau_1/tau_3, t_2 = tau_2/tau_3}
        t_1 = \tau_1/\tau_3, \quad t_2 = \tau_2/\tau_3
\end{equation}
satisfy (\ref{|t_1|^q + |t_2|^q le 1, 2}).  Thus (\ref{N(t_1 v_1 + t_2
  v_2) le 1}) holds under these conditions, by hypothesis.  We also
have that
\begin{equation}
\label{t_1 v_1 + t_2 v_2 = tau_3^{-1} (w_1 + w_2)}
        t_1 \, v_1 + t_2 \, v_2 = \tau_3^{-1} \, (w_1 + w_2),
\end{equation}
by construction, and so (\ref{N(t_1 v_1 + t_2 v_2) le 1}) implies that
\begin{equation}
\label{|tau_3|^{-1} N(w_1 + w_2) = N(tau_3^{-1} (w_1 + w_2)) le 1}
        |\tau_3|^{-1} \, N(w_1 + w_2) = N(\tau_3^{-1} \, (w_1 + w_2)) \le 1.
\end{equation}
Equivalently,
\begin{equation}
\label{N(w_1 + w_2)^q le |tau_3|^q}
        N(w_1 + w_2)^q \le |\tau_3|^q.
\end{equation}
If the positive values of $|\cdot|$ on $k$ are dense in ${\bf R}_+$,
then we can choose $\tau_1, \tau_2 \in k$ so that $|\tau_1|$ and
$|\tau_2|$ are as close as we want to $N(w_1)$, $N(w_2)$,
respectively.  Similarly, we can choose $\tau_3 \in k$ so that
$|\tau_3|^q$ is as close as we want to $|\tau_1|^q + |\tau_2|^q$,
which is as close as we want to $N(w_1)^q + N(w_2)^q$.  Thus
(\ref{N(w_1 + w_2)^q le |tau_3|^q}) implies (\ref{N(w_1 + w_2)^q le
  N(w_1)^q + N(w_2)^q}) in this situation, as desired.

        Let $A$ be a balanced absorbing subset of $V$ with respect
to $|\cdot|$ on $k$, and let $N_A$ be the corresponding Minkowski
functional on $V$, as in (\ref{N_A(v) = inf {|t| : t in k, v in t
    A}}).  Thus $N_A$ satisfies the usual homogeneity condition
(\ref{N(t v) = |t| N(v), 2}), and the open and closed unit balls in
$V$ with respect to $N_A$ are related to $A$ as in (\ref{B_{N_A}(0, 1)
  subseteq A subseteq overline{B}_{N_A}(0, 1)}).  If $A$ is also
$q$-convex in $V$ with respect to $|\cdot|$ on $k$, then
(\ref{B_{N_A}(0, 1) subseteq A subseteq overline{B}_{N_A}(0, 1)})
implies that $N_A$ satisfies the condition mentioned at the beginning
of the preceding paragraph.  More precisely, this means that
(\ref{N(v_1), N(v_2) < 1}) and (\ref{|t_1|^q + |t_2|^q le 1, 2}) imply
(\ref{N(t_1 v_1 + t_2 v_2) le 1}) when $N = N_A$.  It follows that
$N_A$ is a $q$-seminorm on $V$ in this situation.

        Suppose now that $k = {\bf R}$ or ${\bf C}$ with the standard
absolute value function, in which case the previous arguments are
quite classical and can be simplified a bit.  As a variant of the
earlier discussion, let $N$ be a nonnegative real-valued function on
$V$ which is homogeneous of degree $1$ with respect to multiplication
by nonnegative real numbers, as in (\ref{N(t v) = t N(v)}).  As
before, the open and closed balls $B_N(0, r)$ and $\overline{B}_N(0,
r)$ centered at $0$ in $V$ associated to $N$ can be defined as in
(\ref{B_N(0, r) = {v in V : N(v) < r}}) and (\ref{overline{B}_N(0, r)
  = {v in V : N(v) le r}}), and are starlike about $0$ in $V$.  If
either $B_N(0, r)$ or $\overline{B}_N(0, r)$ is real $q$-convex for
some $q \in (0, 1]$ and $r > 0$, then one can check that $B_N(0, r)$
  and $\overline{B}_N(0, r)$ are both real $q$-convex for every $r >
  0$, for the essentially same reasons as before.

        Let $0 < q \le 1$ be given, and suppose that
(\ref{N(t_1 v_1 + t_2 v_2) le 1}) holds for every $v_1, v_2 \in V$
that satisfy (\ref{N(v_1), N(v_2) < 1}) and nonnegative real numbers
$t_1$, $t_2$ such that
\begin{equation}
\label{t_1^q + t_2^q le 1, 2}
        t_1^q + t_2^q \le 1.
\end{equation}
Under these conditions, one can check that (\ref{N(w_1 + w_2)^q le
  N(w_1)^q + N(w_2)^q}) holds for every $w_1, w_2 \in V$, so that $N$
is $q$-subadditive on $V$.  The argument is essentially the same as
before, except that one can take $\tau_1$, $\tau_2$, and $\tau_3$ to
be positive real numbers.  One can also choose $\tau_3$ so that
equality holds in (\ref{|tau_1|^q + |tau_2|^q le |tau_3|^q}),
which means that it suffices to consider $t_1, t_2 \ge 0$ such that
\begin{equation}
\label{t_1^q + t_2^q = 1}
        t_1^q + t_2^q = 1,
\end{equation}
instead of (\ref{t_1^q + t_2^q le 1, 2}).  Note that this condition on
$N$ holds when either $B_N(0, 1)$ or $\overline{B}_N(0, 1)$ is real
$q$-convex in $V$, as before.

        Let $A \subseteq V$ be starlike about $0$ and absorbing in $V$
as a real vector space, and let $\widetilde{N}_A$ be the associated
Minkowski functional on $V$, as in (\ref{widetilde{N}_A(v) = inf {t in
    {bf R}_+ cup {0} : v in t A}}).  Remember that $\widetilde{N}_A$
is homogeneous of degree $1$ with respect to multiplication by
nonnegative real numbers, and that the open and closed unit balls in
$V$ with respect to $\widetilde{N}_A$ are related to $A$ as in
(\ref{B_{widetilde{N}_A}(0, 1) subseteq A subseteq ...}).  If $A$ is
also real $q$-convex for some $0 < q \le 1$, then it is easy to see
that $\widetilde{N}_A$ satisfies the condition described at the
beginning of the previous paragraph, because of
(\ref{B_{widetilde{N}_A}(0, 1) subseteq A subseteq ...}).  This
implies that $\widetilde{N}_A$ is $q$-subadditive on $V$.

        Of course, statements like these are often formulated in terms
of ordinary convexity when $q = 1$.  As in Section \ref{balanced
  q-convexity}, convex subsets of $V$ that contain $0$ are starlike
about $0$, and real $1$-convex.  If $A \subseteq V$ is convex and
absorbing, then $0 \in A$, and $\widetilde{N}_A$ is $1$-subadditive on
$V$, as in the preceding paragraph.  Similarly, if $N$ is a
nonnegative real-valued function on $V$ which is homogeneous of degree
$1$ with respect to multiplication by nonnegative real numbers, then
real $1$-convexity of open or closed balls in $V$ centered at $0$ with
respect to $N$ is equivalent to ordinary convexity.

\section{Minkowski functionals, continued}
\label{minkowski functionals, continued}
\setcounter{equation}{0}

        Let $k$ be a field with a $q$-absolute value function $|\cdot|$
for some $q > 0$, and $V$ be a vector space over $k$.  Also let $B$,
$C$ be balanced absorbing subsets of $V$, and let $N_B$, $N_C$ be the
corresponding Minkowski functionals on $V$, respectively, as in
(\ref{N_A(v) = inf {|t| : t in k, v in t A}}).  If
\begin{equation}
\label{B subseteq C}
        B \subseteq C,
\end{equation}
then it is easy to see that
\begin{equation}
\label{N_C(v) le N_B(v)}
        N_C(v) \le N_B(v)
\end{equation}
for every $v \in V$.  Similarly, suppose that $k = {\bf R}$ or ${\bf
  C}$ with the standard absolute value function, and that $B, C
\subseteq V$ are starlike about $0$, absorbing in $V$ as a real vector
space, and satisfy (\ref{B subseteq C}).  If $\widetilde{N}_B$,
$\widetilde{N}_C$ are the corresponding Minkowski functionals on $V$,
as in (\ref{widetilde{N}_A(v) = inf {t in {bf R}_+ cup {0} : v in t
    A}}), then
\begin{equation}
\label{widetilde{N}_C(v) le widetilde{N}_B(v)}
        \widetilde{N}_C(v) \le \widetilde{N}_B(v)
\end{equation}
for every $v \in V$.

        Let $|\cdot|$ be a nontrivial $q$-absolute value function on
a field $k$, and let $N$ be a nonnegative real-valued function on $V$
that satisfies the usual homogeneity condition (\ref{N(t v) = |t|
  N(v), 2}) with respect to $|\cdot|$ on $k$.  The nontriviality of
$|\cdot|$ on $k$ implies that the open and closed balls (\ref{B_N(0,
  r) = {v in V : N(v) < r}}) and (\ref{overline{B}_N(0, r) = {v in V :
    N(v) le r}}) in $V$ centered at $0$ with radius $r > 0$ with
respect to $N$ are absorbing in $V$, and we have also seen that they
are balanced in $V$.  Let us take
\begin{equation}
\label{C = overline{B}_N(0, 1)}
        C = \overline{B}_N(0, 1)
\end{equation}
for the moment, and consider the corresponding Minkowski functional
$N_C$ on $V$, as in (\ref{N_A(v) = inf {|t| : t in k, v in t A}}).
Suppose that $v \in V$ and $t \in k$ satisfy $v \in t \, C$.  If $t =
0$, then it follows that $v = 0$.  Otherwise, if $t \ne 0$, then
\begin{equation}
\label{v in t C = overline{B}_N(0, |t|)}
        v \in t \, C = \overline{B}_N(0, |t|),
\end{equation}
using (\ref{t overline{B}_N(0, r) = overline{B}_N(0, |t| r)}) in the
second step.  This implies that
\begin{equation}
\label{N(v) le |t|}
        N(v) \le |t|,
\end{equation}
which also works when $t = 0$.  Taking the infimum over $t$, we get
that
\begin{equation}
\label{N(v) le N_C(v)}
        N(v) \le N_C(v)
\end{equation}
for every $v \in V$ under these conditions.

        Suppose now that $|\cdot|$ is not discrete on $k$, so that the
set (\ref{{|x| : x in k, x ne 0}}) of positive values of $|\cdot|$ on $k$
is dense in ${\bf R}_+$ with respect to the standard topology on ${\bf R}$,
as in Section \ref{discrete absolute value functions}.  As before,
\begin{equation}
\label{B = B_N(0, 1)}
        B = B_N(0, 1)
\end{equation}
is balanced and absorbing in $V$, and we let $N_B$ be the
corresponding Minkowski functional on $V$.  Let $v \in V$ be given,
and suppose that $t \in k$ satisfies $N(v) < |t|$, so that
\begin{equation}
\label{v in B_N(0, |t|) = |t| B}
        v \in B_N(0, |t|) = |t| \, B,
\end{equation}
using (\ref{t B_N(0, r) = B_N(0, |t| r)}) in the second step.  This
implies that
\begin{equation}
\label{N_B(v) le |t|}
        N_B(v) \le |t|,
\end{equation}
and hence that
\begin{equation}
\label{N_B(v) le N(v)}
        N_B(v) \le N(v),
\end{equation}
by taking the infimum over $t$ in (\ref{N_B(v) le |t|}).  If $C$ is as
in (\ref{C = overline{B}_N(0, 1)}), then (\ref{B subseteq C}) holds,
and we get that
\begin{equation}
\label{N(v) = N_B(v) = N_C(v)}
        N(v) = N_B(v) = N_C(v)
\end{equation}
for every $v \in V$, by combining (\ref{N_C(v) le N_B(v)}), (\ref{N(v)
  le N_C(v)}), and (\ref{N_B(v) le N(v)}).

        Similarly, let $k = {\bf R}$ or ${\bf C}$ with the standard
absolute value function again, and let $N$ be a nonnegative
real-valued function on $V$ that is homogeneous of degree $1$ with
respect to multiplication by nonnegative real numbers, as in (\ref{N(t
  v) = t N(v)}).  Also let $B$, $C$ be as in (\ref{C =
  overline{B}_N(0, 1)}) and (\ref{B = B_N(0, 1)}), so that $B$, $C$
are starlike about $0$, absorbing in $V$ as a real vector space, and
satisfy (\ref{B subseteq C}).  If $\widetilde{N}_B$, $\widetilde{N}_C$
are the corresponding Minkowski functionals on $V$, then one can check
that
\begin{equation}
\label{N(v) = widetilde{N}_B(v) = widetilde{N}_C(v)}
        N(v) = \widetilde{N}_B(v) = \widetilde{N}_C(v)
\end{equation}
for every $v \in V$, as in the previous two paragraphs.  More
precisely, in this situation, one should restrict one's attention to
nonnegative real numbers $t$ in the earlier arguments.

\section{Continuity of semimetrics}
\label{continuity of semimetrics}
\setcounter{equation}{0}

        Let $X$ be a set, and let $d(x, y)$ be a $q$-semimetric on
$X$ for some $q > 0$.  Thus
\begin{equation}
\label{d(x, z)^q - d(y, z)^q le d(x, y)^q}
        d(x, z)^q - d(y, z)^q \le d(x, y)^q
\end{equation}
for every $x, y, z \in X$, by the $q$-semimetric version (\ref{d(x,
  z)^q le d(x, y)^q + d(y, z)^q}) of the triangle inequality.
Similarly,
\begin{equation}
\label{d(y, z)^q - d(x, z)^q le d(x, y)^q}
        d(y, z)^q - d(x, z)^q \le d(x, y)^q
\end{equation}
for every $x, y, z \in X$, and hence
\begin{equation}
\label{|d(x, z)^q - d(y, z)^q| le d(x, y)^q}
        |d(x, z)^q - d(y, z)^q| \le d(x, y)^q,
\end{equation}
using the standard absolute value function on ${\bf R}$ on the left
side of (\ref{|d(x, z)^q - d(y, z)^q| le d(x, y)^q}).  This implies
that $d(x, z)$ is continuous as a real-valued function of $x \in X$
for each $z \in X$, with respect to the topology determined on $X$ by
$d(\cdot, \cdot)$ as in Sections \ref{semimetrics} and
\ref{q-semimetrics}, and using the standard topology on ${\bf R}$ in
the range of this function.  Of course, this was also implicit in some
of the earlier discussions.  Using the analogous estimate in both
variables, we get that
\begin{eqnarray}
\label{|d(x, w)^q - d(y, z)^q| le ... le d(w, z)^q + d(x, y)^q}
 \qquad |d(x, w)^q - d(y, z)^q|
          & \le & |d(x, w)^q - d(x, z)^q| + |d(x, z)^q - d(y, z)^q| \\
          & \le & d(w, z)^q + d(x, y)^q \nonumber
\end{eqnarray}
for every $w, x, y, z \in X$.  In particular, this implies that $d(x,
w)$ is continuous on $X \times X$, with respect to the product
topology associated to the topology determined on $X$ by $d(\cdot,
\cdot)$.

        Let $\tau$ be a topology on $X$.  Suppose that for each
$u \in X$ and $\epsilon > 0$ there is an open set $U \subseteq X$
with respect to $\tau$ such that $u \in U$ and
\begin{equation}
\label{d(u, v) < epsilon}
        d(u, v) < \epsilon
\end{equation}
for every $v \in U$.  Equivalently, (\ref{d(u, v) < epsilon}) says that
\begin{equation}
\label{U subseteq B_d(u, epsilon)}
        U \subseteq B_d(u, \epsilon),
\end{equation}
where $B_d(u, \epsilon)$ is the open ball in $X$ centered at $u$ with
radius $\epsilon$ with respect to $d(\cdot, \cdot)$, as in (\ref{B(x,
  r) = B_d(x, r) = {y in X : d(x, y) < r}}).  This is also the same as
saying that $u$ is an element of the interior of $B_d(u, \epsilon)$
with respect to $\tau$ for every $u \in X$ and $\epsilon > 0$.  This
condition implies that every open set in $X$ with respect to the
topology on $X$ determined by $d(\cdot, \cdot)$ is an open set with
respect to $\tau$ as well.  Conversely, suppose that every open set in
$X$ with respect to the topology determined by $d(\cdot, \cdot)$ is an
open set in $X$ with respect to $\tau$ too.  Remember that open balls
in $X$ are open sets in $X$ with respect to the topology determined by
$d(\cdot, \cdot)$, as in Sections \ref{semimetrics} and
\ref{q-semimetrics}.  In this case, open balls in $X$ with respect to
$d(\cdot, \cdot)$ are also open sets with respect to $\tau$, which
obviously implies the previous condition.

        The condition described at the beginning of the preceding
paragraph is the same as saying that for each $u \in X$, $d(u, v)$ is
continuous as a real-valued function of $v \in X$ with respect to
$\tau$ at $u$.  This implies that $d(x, y)$ is continuous as a
real-valued function of $x$ or $y$ with respect to $\tau$ on $X$, and
in fact that $d(x, y)$ is continuous with respect to the product
topology on $X \times X$ associated to $\tau$ on each copy of $X$.
This can be derived from (\ref{|d(x, z)^q - d(y, z)^q| le d(x, y)^q})
and (\ref{|d(x, w)^q - d(y, z)^q| le ... le d(w, z)^q + d(x, y)^q}) as
before.  Alternatively, this can be obtained from the analogous
continuity properties with respect to the topology determined on $X$
by $d(\cdot, \cdot)$, and the fact that every open set in $X$ with
respect to the topology determined by $d(\cdot, \cdot)$ is also an
open set with respect to $\tau$.

        Now let $k$ be a field with a $q$-absolute value function
$|\cdot|$ for some $q > 0$, and let $V$ be a vector space over $k$.
If $N$ is a $q$-seminorm on $V$ with respect to $|\cdot|$ on $k$,
then we have that
\begin{equation}
\label{|N(v)^q - N(w)^q| le N(v - w)^q}
        |N(v)^q - N(w)^q| \le N(v - w)^q
\end{equation}
for every $v, w \in V$.  This can be verified in the same way as for
(\ref{|d(x, z)^q - d(y, z)^q| le d(x, y)^q}).  This can also be
considered as a special case of (\ref{|d(x, z)^q - d(y, z)^q| le d(x,
  y)^q}), using the $q$-semimetric associated to $N$ on $V$, as in
(\ref{d(v, w) = N(v - w), 3}).  As before, this implies in particular
that $N(v)$ is continuous as a real-valued function on $V$, with
respect to the topology determined on $V$ by the $q$-semimetric
associated to $N$ as in (\ref{d(v, w) = N(v - w), 3}).

\section{Commutative topological groups}
\label{commutative topological groups}
\setcounter{equation}{0}

        Let $A$ be a commutative group, with the group operations
expressed additively.  Suppose that $A$ is also equipped with a
topology.  We say that $A$ is a \emph{topological group}\index{commutative
topological groups} if the group operations on $A$ are continuous.
More precisely, this means that addition on $A$ should be continuous
as a mapping from $A \times A$ into $A$, where $A \times A$ is
equipped with the product topology associated to the given topology
on $A$.  Similarly,
\begin{equation}
\label{x mapsto -x}
        x \mapsto -x
\end{equation}
should be continuous as a mapping from $A$ into itself, which implies
that this mapping is a homeomorphism, since it is its own inverse.

        Continuity of addition on $A$ implies that
\begin{equation}
\label{x mapsto x + a}
        x \mapsto x + a
\end{equation}
is a continuous mapping from $A$ into itself for each $a \in A$.  This
corresponds to continuity of addition on $A$ in each variable
separately, instead of joint continuity of addition as a mapping from
$A \times A$ into $A$.  Of course, the inverse of the translation
mapping (\ref{x mapsto x + a}) is given by translation by $-a$, so
that continuity of translations implies that the translation mappings
(\ref{x mapsto x + a}) are homeomorphisms on $A$.  If $A$ is equipped
with a topology for which the translation mappings (\ref{x mapsto x +
  a}) are continuous, and if addition on $A$ is continuous as a
mapping from $A \times A$ into $A$ with respect to the product
topology on $A \times A$ at the point $(0, 0)$ in $A \times A$, then
one can check that addition on $A$ is continuous as a mapping from $A
\times A$ into $A$ at every point in $A \times A$.  Similarly, if the
translation mappings (\ref{x mapsto x + a}) are continuous on $A$, and
if (\ref{x mapsto -x}) is continuous at $0$, then (\ref{x mapsto -x})
is continuous everywhere on $A$.

        Let us suppose for the rest of the section that $A$ is a
commutative topoloigcal group.  If $a \in A$ and $B \subseteq A$, then
put
\begin{equation}
\label{a + B = B + a = {a + b : b in B}}
        a + B = B + a = \{a + b : b \in B\},
\end{equation}
which is the same as the image of $B$ under the translation mapping
(\ref{x mapsto x + a}).  If $B, C \subseteq A$, then we put
\begin{equation}
\label{B + C = {b + c : b in B, c in C} = ...}
        B + C = \{b + c : b \in B, \ c \in C\}
              = \bigcup_{b \in B} (b + C) = \bigcup_{c \in C} (B + c).
\end{equation}
In particular, if either $B$ or $C$ is an open set in $A$, then $B +
C$ is an open set in $A$ as well, because it is a union of open sets
in $A$.  Let us also put
\begin{eqnarray}
\label{-C = {-c : c in C}}
        -C & = & \{-c : c \in C\},                \\
         b - C & = & b + (-C),
\end{eqnarray}
and
\begin{equation}
\label{B - C = B + (-C)}
        B - C = B + (-C)
\end{equation}
for every $b \in A$ and $B, C \subseteq A$.

        Let $E$ be any subset of $A$, and let $V \subseteq A$ be an open
set that contains $0$.  We would like to check that
\begin{equation}
\label{overline{E} subseteq E + V}
        \overline{E} \subseteq E + V,
\end{equation}
where $\overline{E}$ is the closure of $E$ in $V$, as usual.  If $x
\in \overline{E}$, then
\begin{equation}
\label{(x - V) cap E ne emptyset}
        (x - V) \cap E \ne \emptyset,
\end{equation}
because $x - V$ is an open set in $V$ that contains $x$.  This is the
same as saying that $x \in E + V$, which implies (\ref{overline{E}
  subseteq E + V}).  Similarly, if $x \in E + V$ for every open set $V
\subseteq A$ that contains $0$, then $x \in \overline{E}$, which
implies that
\begin{equation}
\label{overline{E} = bigcap {E + V : V subseteq A is an open set with 0 in V}}
        \overline{E} = \bigcap \{E + V : V \subseteq A
                                 \hbox{ is an open set with } 0 \in V\}.
\end{equation}

        Let $W$ be an open subset of $A$ that contains $0$.  Continuity
of addition on $A$ as a mapping from $A \times A$ into $A$ at $(0, 0)$
says exactly that there are open sets $U, V \subseteq A$ that contain
$0$ and satisfy
\begin{equation}
\label{U + V subseteq W}
        U + V \subseteq W.
\end{equation}
In particular, this implies that
\begin{equation}
\label{overline{U} subseteq W}
        \overline{U} \subseteq W,
\end{equation}
as in (\ref{overline{E} subseteq E + V}).  If $x$ is any element of
$A$, then it follows that every open set in $A$ that contains $x$ also
contains the closure of another open subset of $A$ that contains $x$,
by continuity of translations on $A$.  This means that $A$ is regular
as a topological space in the strict sense, without including the
first or $0$th separation condition.

        Sometimes the condition that $\{0\}$ be a closed set in $A$ is
included in the definition of a commutative topological
group.\index{commutative topological groups} This implies that every
subset of $A$ with exactly one element is a closed set, by continuity
of translations, so that $A$ satisfies the first separation condition.
Combining this with the remarks in the preceding paragraph, we get
that $A$ is regular in the stronger sense that includes the first
separation condition, which is sometimes referred to as the third
separation condition.  In particular, this implies that $A$ is Hausdorff
as a topological space.

        Here is another way to look at the Hausdorff property of $A$
when $\{0\}$ is a closed set in $A$.  Let $x, y \in A$ be given, with
$x \ne y$, and put
\begin{equation}
\label{W = A setminus {x - y}}
        W = A \setminus \{x - y\},
\end{equation}
so that $0 \in W$.  Note that $W$ is an open set in $A$, since $A$
satisfies the first separation condition.  Thus there are open sets
$U, V \subseteq A$ that contain $0$ and satisfy (\ref{U + V subseteq
  W}), which means that
\begin{equation}
\label{x - y not in U + V}
        x - y \not\in U + V.
\end{equation}
Equivalently, this implies that
\begin{equation}
\label{(x - U) cap (y + V) = emptyset}
        (x - U) \cap (y + V) = \emptyset,
\end{equation}
and hence that $A$ is Hausdorff as a topological space, because $x -
U$ and $y + V$ are disjoint open subsets of $A$ that contain $x$ and
$y$, respectively.

        Similarly, let $x \in A$ and $E \subseteq A$ be given, with
$x \not\in E$, and put
\begin{equation}
\label{W = A setminus (x - E)}
        W = A \setminus (x - E),
\end{equation}
so that $0 \in W$.  If $E$ is a closed subset of $A$, then $x - E$ is
a closed set in $A$ too, and hence $W$ is an open set in $A$.  This
implies that there are open sets $U, V \subseteq A$ that contain $0$
and satisfy (\ref{U + V subseteq W}), as before.  In this case,
(\ref{U + V subseteq W}) says that
\begin{equation}
\label{(U + V)  cap (x - E) =  emptyset}
        (U + V) \cap (x - E) = \emptyset,
\end{equation}
which implies that
\begin{equation}
\label{(x - U)  cap (E + V) = emptyset}
        (x - U) \cap (E + V) = \emptyset.
\end{equation}
This is another way to look at the regularity of $A$ as a topological
space, since $x - U$ and $E + V$ are disjoint open subsets of $A$ that
contain $x$ and $E$, respectively.

\section{Translation-invariant semimetrics}
\label{translation-invariant semimetrics}
\setcounter{equation}{0}

        Let $A$ be a commutative group, and let $d(x, y)$ be a semimetric
on $A$.  If
\begin{equation}
\label{d(x + a, y + a) = d(x, y)}
        d(x + a, y + a) = d(x, y)
\end{equation}
for every $a, x, y \in A$, then we say that $d(\cdot, \cdot)$ is
\emph{invariant under translations}\index{translation-invariant
  semimetrics} on $A$.  In this case, the translation mappings (\ref{x
  mapsto x + a}) automaticvally define homeomorphisms on $A$ with
respect to the topology on $A$ determined by $d(\cdot, \cdot)$.  One
can also check that addition on $A$ is continuous as a mapping from $A
\times A$ into $A$ under these conditions, using the triangle
inequality.  Similarly, it is easy to see that (\ref{x mapsto -x}) is
continuous on $A$, because $d(x, y)$ is symmetric in $x$ and $y$.
This means that $A$ satisfies the requirements of a commutative
topological group with respect to the topology determined by $d(\cdot,
\cdot)$.  If $d(\cdot, \cdot)$ is a translation-invariant metric on
$A$, then $A$ is also Hausdorff with respect to this topology.  Of
course, invariance under translations can be defined for
$q$-semimetrics on $A$ in the same way, for any $q > 0$.
Equivalently, a $q$-semimetric $d(x, y)$ on $A$ is invariant under
translations if and only if $d(x, y)^q$ is invariant under
translations as a semimetric on $A$, in which case there are analogues
of the previous statements for the corresponding topology.

        Let $\mathcal{M}$ be a collection of $q$-semimetrics on $A$,
where $q > 0$ is allowed to depend on the element of $\mathcal{M}$.
This leads to a topology on $A$, as in Sections \ref{collections of
  semimetrics} and \ref{q-semimetrics}.  If the elements of
$\mathcal{M}$ are invariant under translations on $A$, then the group
operations on $A$ are continuous with respect to this topology, as in
the preceding paragraph, so that $A$ becomes a commutative topological
group.  If $\mathcal{M}$ is also nondegenerate on $A$, then we have
seen that $A$ is Hausdorff with respect to this topology.  If $A$ is
any commutative topological group, then it is well known that the
topology on $A$ corresponds to a collection of translation-invariant
semimetrics on $A$ in this way.  Note that any commutative group is a
topological group with respect to the discrete topology.  The discrete
metric on any group is invariant under translations, and determines
the discrete topology on the group.

        Let $A$ be a commutative topological group, and let $d(x, y)$
be a translation-invariant $q$-semimetric on $A$ for some $q > 0$.
Suppose that $d(0, v)$ is continuous as a real-valued function of
$v \in A$ at $0$, with respect to the standard topology on ${\bf R}$
in the range of this function.  This implies that for each $u \in A$,
$d(u, v)$ is continuous as a real-valued function of $v \in A$ at $u$,
because of translation-invariance of $d(\cdot, \cdot)$ and the
topology on $A$.  As in Section \ref{continuity of semimetrics},
it follows that $d(x, y)$ is continuous as a real-valued function
of $x$ or $y$ on $A$, and in fact that $d(x, y)$ is continuous as
a real-valued function on $A \times A$ with respect to the corresponding
product topology.  This condition also implies that every open set
in $A$ with respect to the topology determined by $d(\cdot, \cdot)$
is an open set in $A$ with respect to its given topology, as before.

        If the topology on $A$ is determined by a collection $\mathcal{M}$
of translation-invariant $q$-semimetrics on $A$, then each element of
$\mathcal{M}$ automatically has the property described in the preceding
paragraph.  If $d(x, y)$ is any other translation-invariant $q$-semimetric
on $A$ with this property, then one can add $d(x, y)$ to $\mathcal{M}$,
and get the same topology on $A$.  If $A$ is any commutative topological
group, then one might like to find a collection $\mathcal{M}$ of 
translation-invariant $q$-semimetrics on $A$ that determines the same
topology on $A$, as mentioned earlier.  Each element of $\mathcal{M}$
should be compatible with the given topology on $A$, as in the previous
paragraph.  The point of the theorem mentioned earlier is that one can
always find enough translation-invariant semimetrics on $A$ that are
compatible with the given topology on $A$ to generate this topology on $A$.

\section{Topological vector spaces}
\label{topological vector spaces}
\setcounter{equation}{0}

        Let $k$ be a field with a $q$-absolute value function
$|\cdot|$ for some $q > 0$, and let $V$ be a vector space over $k$.
Remember that $|x - y|$ defines a $q$-metric on $k$, which determines
a topology on $k$ in the usual way.  We say that $V$ is a
\emph{topological vector space}\index{topological vector spaces} with
respect to $|\cdot|$ on $k$ if $V$ is equipped with a topology for
which the vector space operations are continuous.  More precisely, the
continuity of addition on $V$ means that addition defines a continuous
mapping from $V \times V$ into $V$, with respect to the product
topology on $V \times V$ associated to the given topology on $V$.
Similarly, continuity of scalar multiplication on $V$ means that
scalar multiplication defines a continuous mapping from $k \times V$
into $V$, where $k \times V$ is equipped with the product topology
associated to the topology on $k$ mentioned earlier and the given
topology on $V$.

        In particular, continuity of scalar multiplication implies
that for each $t \in k$,
\begin{equation}
\label{v mapsto t v}
        v \mapsto t \, v
\end{equation}
is a continuous mapping on $V$.  This includes the continuity of
\begin{equation}
\label{v mapsto -v}
        v \mapsto -v
\end{equation}
on $V$, by taking $t = -1$ in $k$ in (\ref{v mapsto t v}).  Of course,
vector spaces are commutative groups with respect to addition, and
this implies that topological vector spaces are topological groups
with respect to addition.  If $\{0\}$ is a closed set $V$, then it
follows that $V$ is Hausdorff as a topological space, as in Section
\ref{commutative topological groups}.  This additional condition is
sometimes included in the definition of a topological vector
space,\index{topological vector spaces} as before.

        If $k = {\bf R}$ or ${\bf C}$ with the standard absolute
value function, then this corresponds to the usual notion of a real or
complex topological vector space.  If $|\cdot|$ is a nontrivial
absolute value function on any field $k$, then this definition
corresponds to the one in \cite{sch-w}.  The restriction to $q = 1$ is
not too serious, since $|x|^q$ defines an absolute value function on
$k$ when $|x|$ is a $q$-absolute value function on $k$, and which
leads to the same topology on $k$.  If $|\cdot|$ is the trivial
absolute value function on $k$, then the corresponding topology on $k$
is discrete, and the continuity of scalar multiplication on $V$ as a
mapping from $k \times V$ into $V$ with respect to the product toplogy
on $k \times V$ is the same as the continuity of (\ref{v mapsto t v})
on $V$ for each $t \in k$.  One may wish to exclude this case in some
situations, as in \cite{sch-w}.\index{topological vector spaces}

        Let $|\cdot|$ be a $q$-absolute value function on a field $k$
for some $q > 0$ again, and let $V$ be a topological vector space over
$k$.  Note that (\ref{v mapsto t v}) is a homeomorphism on $V$ for
each $t \in k$ with $t \ne 0$, since the inverse mapping is given by
multiplication by $1/t$.  Let $W$ be an open subset of $V$ that
contains $0$.  Continuity of scalar multiplication at $(0, 0)$ in $k
\times V$ implies that there is an open set $U \subseteq V$ with $0
\in U$ and a $\delta > 0$ such that
\begin{equation}
\label{t U subseteq W}
        t \, U \subseteq W
\end{equation}
for every $t \in k$ with $|t| < \delta$, where $t \, U$ is as in
(\ref{t E = {t v : v in E}}).  This statement is vacuous when
$|\cdot|$ is the trivial absolute value function on $k$, and so let us
suppose for the moment that $|\cdot|$ is nontrivial on $k$.  This
implies that there are $t \in k$ with $t \ne 0$ and $|t| < \delta$,
and we put
\begin{equation}
\label{U_1 = bigcup {t U : t in k, 0 < |t| < delta}}
        U_1 = \bigcup \{t \, U : t \in k, \ 0 < |t| < \delta\},
\end{equation}
which contains $0$ because $0 \in U$.  We also have that
\begin{equation}
\label{U_1 subseteq W}
        U_1 \subseteq W,
\end{equation}
by (\ref{t U subseteq W}), and that $U_1$ is an open set in $V$, since
it is a union of open sets.  By construction, $U_1$ is a balanced
subset of $V$ too, as in Section \ref{balanced sets}.

        Let us rephrase this a bit, as follows.  If $W$ is an open set
in $V$ that contains $0$, and if $|\cdot|$ is nontrivial on $k$, then
there is an open set $\widetilde{U} \subseteq V$ such that $0 \in
\widetilde{U}$ and
\begin{equation}
\label{t widetilde{U} subseteq W}
        t \, \widetilde{U} \subseteq W
\end{equation}
for every $t \in k$ with $|t| \le 1$.  More precisely, if $U$
satisfies the conditions described in the preceding paragraph, then
one can take $\widetilde{U}$ to be $U_1$, as in (\ref{U_1 = bigcup {t
    U : t in k, 0 < |t| < delta}}).  Conversely, if $\widetilde{U}$
satisfies these conditions, then one can simply take $U$ to be
$\widetilde{U}$ in the previous paragraph.  In this case,
\begin{equation}
\label{widetilde{U}_1 = bigcup {t widetilde{U} : t in k, |t| le 1}}
 \widetilde{U}_1 = \bigcup \{t \, \widetilde{U} : t \in k, \ |t| \le 1\}
\end{equation}
is a nonempty balanced open subset of $V$ that is contained in $W$, as
before.

        If $|\cdot|$ is the trivial absolute value function on $k$,
then the existence of $\widetilde{U}$ as in (\ref{t widetilde{U}
  subseteq W}) still makes sense, even if it is not necessarily
included in the definition of a topological vector space.  It would
also be stronger than before, in the sense that it applies to every $t
\in k$, and similarly balanced subsets of $V$ would be invariant under
multiplication by any $t \in k$ with $t \ne 0$, as in Section
\ref{balanced sets}.  If $|\cdot|$ is any $q$-absolute value function
on $k$ and $V$ is a topological vector space with respect to $|\cdot|$
on $k$, then $V$ is also a topological vector space with respect to
the trivial absolute value function on $k$, because (\ref{v mapsto t
  v}) is still continuous on $V$ for every $t \in k$.  In this case,
(\ref{t widetilde{U} subseteq W}) would not normally work with respect
to the trivial absolute value function on $k$.  In particular, one can
take $V = k$ with the topology associated to a nontrivial absolute
value function on $k$.

        Let $|\cdot|$ be any $q$-absolute value function on a field
$k$ again.  If $V$ is a topological vector space over $k$, then continuity
of scalar multiplication implies that
\begin{equation}
\label{t mapsto t v}
        t \mapsto t \, v
\end{equation}
is a continuous mapping from $k$ into $V$ for every $v \in V$.  Let
$W$ be an open subset of $V$ that contains $0$, and let $v \in V$ be
given.  The continuity of (\ref{t mapsto t v}) at $0$ implies that
there is a $\delta(v) > 0$ such that
\begin{equation}
\label{t v in W}
        t \, v \in W
\end{equation}
for every $t \in k$ with $|t| < \delta(v)$.  This is the same as
saying that $W$ is absorbing in $V$ when $|\cdot|$ is nontrivial on
$k$, as in Section \ref{absorbing sets}.  If $|\cdot|$ is the trivial
absolute value function on $k$, then this condition on $W$ is vacuous,
and we have seen that $V$ is the only absorbing subset of itself in
this case.

\section{Compatible seminorms}
\label{compatible seminorms}
\setcounter{equation}{0}

        Let $k$ be a field equipped with a $q$-absolute value function
$|\cdot|$ for some $q > 0$, and let $V$ be a vector space over $k$.
Also let $\mathcal{N}$ be a collection of $q$-seminorms on $V$ with
respect to $|\cdot|$ on $k$.  More precisely, one can let $q > 0$
depend on the element of $\mathcal{N}$, as long as $|\cdot|$ is a
$q$-absolute value function on $k$ for each such $q$.  Every element
of $\mathcal{N}$ determines a $q$-semimetric on $V$ in the usual way,
as in (\ref{d(v, w) = N(v - w), 3}).  Note that these $q$-semimetrics
are automatically invariant under translations on $V$.  Consider the
topology on $V$ determined by the collection $\mathcal{M}$ of
$q$-semimetrics associated to the elements of $\mathcal{N}$.  One can
check that scalar multiplication on $V$ defines a continuous mapping
from $k \times V$ into $V$ with respect to this topology on $V$, using
the product topology on $k \times V$ corresponding to this topology on
$V$ and the topology on $k$ determined by the $q$-metric associated to
$|\cdot|$.  This implies that $V$ is a topological vector space with
respect to $|\cdot|$ on $k$ under these conditions.  If $\mathcal{N}$
is nondegenerate on $V$, then $V$ is Hausdorff with respect to this
topology, as before.  If $|\cdot|$ is the trivial absolute value
function on $k$, then this topology on $V$ also satisfies the
condition (\ref{t widetilde{U} subseteq W}) discussed in the previous
section, because of the corresponding homogeneity property of
$q$-seminorms on $V$.

        Suppose now that $V$ is equipped with a topology that makes
it a topological vector space with respect to $|\cdot|$ on $k$, and
let $N$ be a $q$-seminorm with respect to $|\cdot|$ on $k$.  Suppose
that $N$ is continuous at $0$ as a real-valued function on $V$, where
${\bf R}$ is equipped with the standard topology in the range of this
function.  It is easy to see that this implies that $N$ is continuous
on all of $V$, using (\ref{|N(v)^q - N(w)^q| le N(v - w)^q}).  In
particular, it follows that open balls in $V$ with respect to $N$ are
open sets in $V$ with respect to the given topology on $V$.  One can
look at this in the context of Section \ref{translation-invariant
  semimetrics} as well, since $V$ is a commutative topological group
with respect to addition.  Let $d(v, w)$ be the $q$-semimetric on $V$
associated to $N$ as in (\ref{d(v, w) = N(v - w), 3}), which is
automatically invariant under translations on $V$, as mentioned
earlier.  The hypothesis that $N$ be continuous at $0$ is the same as
saying that $d(0, w)$ is continuous as a real-valued function of $w
\in V$ at $0$.  As in Section \ref{translation-invariant semimetrics},
this implies that $d(v, w)$ is continuous in $v$ and $w$, which
corresponds to the continuity of $N$ on $V$.  This also implies that
open subsets of $V$ with respect to the topology determined by
$d(\cdot, \cdot)$ are open sets with respect to the given topology on
$V$, as before.

        As in Section \ref{continuity of semimetrics}, the condition
that $N$ be continuous at $0$ is the same as saying that for each
$\epsilon > 0$, $0$ is an element of the interior of $B_N(0,
\epsilon)$ in $V$, where $B_N(0, \epsilon)$ is as in (\ref{B_N(0, r) =
  {v in V : N(v) < r}}).  This implies that $N$ is continuous on all
of $V$, as in the preceding paragraph, and hence that $B_N(0, r)$ is
an open set in $V$ for every $r > 0$.  Remember that
\begin{equation}
\label{t B_N(0, epsilon) = B_N(0, |t| epsilon)}
        t \, B_N(0, \epsilon) = B_N(0, |t| \, \epsilon),
\end{equation}
for every $\epsilon > 0$ and $t \in k$ with $t \ne 0$, as in (\ref{t
  B_N(0, r) = B_N(0, |t| r)}).  If $0$ is an element of the interior
of $B_N(0, \epsilon)$ in $V$ for some $\epsilon > 0$, then it follows
that $0$ is also an element of the interior of (\ref{t B_N(0, epsilon)
  = B_N(0, |t| epsilon)}) for every $t \in k$ with $t \ne 0$.  This
implies that $N$ is continuous at $0$ when $|\cdot|$ is not the
trivial absolute value function on $k$.

        Let $A$ be a balanced absorbing subset of $V$, and let $N_A$
be the corresponding Minkowski functional, as in (\ref{N_A(v) = inf
  {|t| : t in k, v in t A}}).  Thus
\begin{equation}
\label{N_A(v) le |t|}
        N_A(v) \le |t|
\end{equation}
for every $t \in k$ and $v \in t \, A$, by definition of $N_A(v)$.  If
$|\cdot|$ is nontrivial on $k$, and if $0$ is an element of the
interior of $A$ in $V$, then it follows from (\ref{N_A(v) le |t|})
that $N_A$ is continuous at $0$ as a real-valued function on $V$.  Of
course, if $|\cdot|$ is trivial on $k$, then $V$ is the only absorbing
subset of itself, and $N_V$ is equal to $0$ on all of $V$, as in
Section \ref{minkowski functionals}.  Conversely, if $N_A$ is
continuous at $0$ on $V$, then $0$ is an element of the interior of
$B_{N_A}(0, 1)$ in $V$, where $B_{N_A}(0, 1)$ is defined as in
(\ref{B_N(0, r) = {v in V : N(v) < r}}).  We have also seen that
$B_{N_A}(0, 1)$ is contained in $A$ when $A$ is balanced, as in
(\ref{B_{N_A}(0, 1) subseteq A subseteq overline{B}_{N_A}(0, 1)}).
This implies that $0$ is an element of the interior of $A$ in $V$
when $N_A$ is continuous at $0$ on $V$.

\section{Subadditive functions}
\label{subadditive functions}
\setcounter{equation}{0}

        Let us take $k = {\bf R}$ with the standard absolute value
function in this section, and let $V$ be a vector space over ${\bf
  R}$.  Also let $N$ be a nonnegative real-valued function on $V$
which is homogeneous of degree $1$ with respect to multiplication by
nonnegative real numbers, as in (\ref{N(t v) = t N(v)}).  Suppose that
$N$ is $q$-subadditive on $V$ for some $q > 0$, as in (\ref{N(v + w)^q
  le N(v)^q + N(w)^q, 2}).  More precisely, we may as well take $0 < q
\le 1$ here, so that the standard absolute value function on ${\bf R}$
is a $q$-absolute value function, as in Section \ref{q-absolute value
  functions}.  Observe that $N(-v)$ satisfies these same properties
on $V$, which is to say that $N(-v)$ is also homogeneous with respect
to multiplication by nonnegative real numbers and $q$-subadditive.
This implies that
\begin{equation}
\label{N^{sym}(v) = max(N(v), N(-v))}
        N^{sym}(v) = \max(N(v), N(-v))
\end{equation}
has the same properties too, because the maximum of two $q$-suabdditive
functions is $q$-subadditive.  By construction,
\begin{equation}
\label{N^{sym}(-v) = N^{sym}(v)}
        N^{sym}(-v) = N^{sym}(v)
\end{equation}
for every $v \in V$, so that $\widetilde{N}$ is a $q$-seminorm on $V$.
As in Section \ref{continuity of semimetrics}, we have that
\begin{equation}
\label{N(v)^q - N(w)^q le N(v - w)^q}
        N(v)^q - N(w)^q \le N(v - w)^q
\end{equation}
and
\begin{equation}
\label{N(w)^q - N(v)^q le N(w - v)^q}
        N(w)^q - N(v)^q \le N(w - v)^q
\end{equation}
for every $v, w \in V$, by $q$-subadditivity.  It follows that
\begin{equation}
\label{|N(v)^q - N(w)^q| le N^{sym}(v - w)^q}
        |N(v)^q - N(w)^q| \le N^{sym}(v - w)^q
\end{equation}
for every $v, w \in V$, so that $N$ is continuous as a real-valued
function on $V$ with respect to the topology determined on $V$ by the
$q$-semimetric associated to $N^{sym}$.

        Suppose now that $V$ is equipped with a topology which makes
it into a topological vector space over ${\bf R}$.  If $N(v)$ is
continuous at $0$ as a real-valued function on $V$ with respect to
this topology, then $N(-v)$ is also continuous at $0$, and hence
$N^{sym}(v)$ is continuous at $0$.  This implies that $N$ is continuous
on all of $V$, by (\ref{|N(v)^q - N(w)^q| le N^{sym}(v - w)^q}).
In particular, this means that open balls in $V$ with respect to
$N$ are open sets in $V$.  The continuity of $N^{sym}$ at $0$ on $V$
also implies that $N^{sym}$ is continuous on all of $V$, as in the
previous section.

        Let $A$ be a subset of $V$ that is starlike about $0$ and
absorbing, and let $\widetilde{N}_A$ be the corresponding Minkowski
functional, as in (\ref{widetilde{N}_A(v) = inf {t in {bf R}_+ cup {0}
    : v in t A}}).  It is easy to see that $-A$ is starlike about $0$
and absorbing in $V$ too, and we let $\widetilde{N}_{-A}$ be the
Minkowski functional associated to $-A$ on $V$ as in
(\ref{widetilde{N}_A(v) = inf {t in {bf R}_+ cup {0} : v in t A}}).
Similarly,
\begin{equation}
\label{A^{sym} = A cap (-A)}
        A^{sym} = A \cap (-A)
\end{equation}
is starlike about $0$ and absorbing in $V$, and in fact $A^{sym}$ is
balanced in $V$, since it is automatically symmetric about $0$.  Thus
the Minkowski functional $N_{A^{sym}}$ associated to $A^{sym}$ as in
(\ref{N_A(v) = inf {|t| : t in k, v in t A}}) is the same as its
analogue $\widetilde{N}_{A^{sym}}$ as in (\ref{widetilde{N}_A(v) = inf
  {t in {bf R}_+ cup {0} : v in t A}}).  One can check that
\begin{equation}
\label{widetilde{N}_A(-v) = widetilde{N}_{-A}(v)}
        \widetilde{N}_A(-v) = \widetilde{N}_{-A}(v)
\end{equation}
and
\begin{equation}
\label{N_{A^{sym}}(v) = max(widetilde{N}_A(v), widetilde{N}_A(-v))}
        N_{A^{sym}}(v) = \max(\widetilde{N}_A(v), \widetilde{N}_A(-v))
\end{equation}
for every $v \in V$.  

        As in the previous section,
\begin{equation}
\label{widetilde{N}_A(v) le t}
        \widetilde{N}_A(v) \le t
\end{equation}
for every $v \in t \, A$ when $t \ge 0$, by the definition of
$\widetilde{N}_A(v)$.  If $0$ is an element of the interior of $A$ in
$V$, then it follows that $\widetilde{N}_A(v)$ is continuous at $0$ on
$V$, as before.  Conversely, if $\widetilde{N}_A(v)$ is continuous at
$0$ on $V$, then $0$ is an element of the interior of
$B_{\widetilde{N}_A}(0, 1)$ in $V$, where $B_{\widetilde{N}_A}(0, 1)$
is as defined in (\ref{B_N(0, r) = {v in V : N(v) < r}}).  This
implies that $0$ is an element of the interior of $A$ in $V$, because
$B_{\widetilde{N}_A}(0, 1)$ is contained in $A$ when $A$ is starlike
about $0$, as in (\ref{B_{widetilde{N}_A}(0, 1) subseteq A subseteq
  ...}).  Of course, if $0$ is an element of the interior of $A$ in
$V$, then $0$ is also an element of the interior of $-A$, which
implies that $0$ is an element of the interior of $A^{sym}$ as well.

\section{Cartesian products}
\label{cartesian products}
\setcounter{equation}{0}

        Let $I$ be a nonempty set, let $X_j$ be a set for each $j \in I$,
and let
\begin{equation}
\label{X = prod_{j in I} X_j}
        X = \prod_{j \in I} X_j
\end{equation}
be the Cartesian product of the $X_j$'s.  If $x \in X$ and $j \in I$,
then it will be convenient to let $x_j$ be the $j$th coordinate of $x$
in $X_j$.  Let $\mathcal{M}_j$ be a collection of $q$-semimetrics on
$X_j$ for each $j \in I$, where $q > 0$ is allowed to depend on the
element of $\mathcal{M}_j$, as before.  If $d_j(\cdot, \cdot)$ is an
element of $\mathcal{M}_j$ for some $j \in I$, then it is easy to see
that
\begin{equation}
\label{widehat{d}_j(x, y) = d_j(x_j, y_j)}
        \widehat{d}_j(x, y) = d_j(x_j, y_j)
\end{equation}
defines a $q$-semimetric on $X$, with the same $q$ as for $d_j(\cdot,
\cdot)$.  Put
\begin{equation}
\label{widehat{mathcal{M}}_j = {widehat{d}_j : d_j in mathcal{M}_j}}
 \widehat{\mathcal{M}}_j = \{\widehat{d}_j(\cdot, \cdot) :
                              d_j(\cdot, \cdot) \in \mathcal{M}_j\}
\end{equation}
for each $j \in I$, and
\begin{equation}
\label{mathcal{M} = bigcup_{j in I} widehat{mathcal{M}}_j}
        \mathcal{M} = \bigcup_{j \in I} \widehat{\mathcal{M}}_j.
\end{equation}
Let $X_j$ be equipped with the topology associated to $\mathcal{M}_j$
as in Sections \ref{collections of semimetrics} and
\ref{q-semimetrics}, for each $j \in I$.  One can check that the
topology on $X$ associated to $\mathcal{M}$ in this way is the same as
the corresponding product topology on $X$.  Note that $\mathcal{M}$ is
nondegenerate on $X$ when $\mathcal{M}_j$ is nondegenerate on $X_j$
for every $j \in I$.

        Now let $A_j$ be a commutative group for each $j \in I$, and
let
\begin{equation}
\label{A = prod_{j in I} A_j}
        A = \prod_{j \in I} A_j
\end{equation}
be the Cartesian product of the $A_j$'s.  This is also a commutative
group, where the group operations are defined coordinatewise.  This
group is known as the \emph{direct product}\index{direct products} of
the $A_j$'s.  If $A_j$ is a commutative topological group for each $j
\in I$, then $A$ is a commutative topological group too, with respect
to the corresponding product topology.  Let $\mathcal{M}_j$ be a
collection of translation-invariant $q$-semimetrics on $A$, and let
$\widehat{\mathcal{M}}_j$ be the corresponding collection of
$q$-semimetrics on $A$ for each $j \in I$, as in
(\ref{widehat{mathcal{M}}_j = {widehat{d}_j : d_j in mathcal{M}_j}}).
Note that the elements of $\widehat{\mathcal{M}}_j$ are invariant
under translations on $A$ for each $j$, so that the union
$\mathcal{M}$ of the $\widehat{\mathcal{M}}_j$'s consists of
translation-invariant $q$-semimetrics on $A$ as well.  If $A_j$ is
equipped with the topology determined by $\mathcal{M}_j$ for each $j
\in I$, then the topology on $A$ associated to $\mathcal{M}$ is the
same as the corresponding product topology, as in the previous
paragraph.

        Let $k$ be a field, and let $V_j$ be a vector space over $k$
for each $j \in I$.  The Cartesian product
\begin{equation}
\label{V = prod_{j in I} V_j}
        V = \prod_{j \in I} V_j
\end{equation}
is a vector space over $k$ too, with respect to coordinatewise
addition and scalar multiplication.  This vector space is known as the
\emph{direct product}\index{direct products} of the $V_j$'s.  Let
$|\cdot|$ be a $q$-absolute value function on $k$ for some $q > 0$,
and suppose that $V_j$ is a topological vector space with respect to
$|\cdot|$ on $k$ for each $j \in I$.  Under these conditions, one can
check that $V$ is a topological vector space too, with respect to the
corresponding product topology.

        Let $\mathcal{N}_j$ be a collection of $q$-seminorms on $V_j$
for each $j \in I$, where $q > 0$ is allowed to depend on the element
of $\mathcal{N}_j$, as long as $|\cdot|$ is a $q$-absolute value function
on $k$.  If $N_j$ is an element of $\mathcal{N}_j$ for some $j \in I$,
then
\begin{equation}
\label{widehat{N}_j(v) = N_j(v_j)}
        \widehat{N}_j(v) = N_j(v_j)
\end{equation}
defines a $q$-seminorm on $V$ with the same $q$ as for $N_j$, where
$v_j \in V_j$ is the $j$th coordinate of $v \in V$, as before.  In
analogy with (\ref{widehat{mathcal{M}}_j = {widehat{d}_j : d_j in
    mathcal{M}_j}}) and (\ref{mathcal{M} = bigcup_{j in I}
  widehat{mathcal{M}}_j}), put
\begin{equation}
\label{widehat{mathcal{N}}_j = {widehat{N}_j : N_j in mathcal{N}_j}}
        \widehat{\mathcal{N}}_j = \{\widehat{N}_j : N_j \in \mathcal{N}_j\}
\end{equation}
for each $j \in I$, and
\begin{equation}
\label{mathcal{N} = bigcup_{j in I} widehat{mathcal{N}}_j}
        \mathcal{N} = \bigcup_{j \in I} \widehat{\mathcal{N}}_j.
\end{equation}
Let $\mathcal{M}_j$ be the collection of $q$-semimetrics on $V_j$
corresponding to elements of $\mathcal{N}_j$ as in (\ref{d(v, w) = N(v
  - w), 3}) for each $j \in I$.  If $\widehat{\mathcal{M}}_j$ is
associated to $\mathcal{M}_j$ as in (\ref{widehat{mathcal{M}}_j =
  {widehat{d}_j : d_j in mathcal{M}_j}}) for each $j \in I$, then
$\widehat{\mathcal{M}}_j$ is the same as the collection of
$q$-semimetrics on $V$ that correspond to elements of
$\widehat{\mathcal{N}}_j$ as in (\ref{d(v, w) = N(v - w), 3}) for each
$j \in I$.  Similarly, if $\mathcal{M}$ is as in (\ref{mathcal{M} =
  bigcup_{j in I} widehat{mathcal{M}}_j}), then $\mathcal{M}$ is the
same as the collection of $q$-semimetrics on $V$ that correspond to
elements of $\mathcal{N}$ as in (\ref{d(v, w) = N(v - w), 3}) for each
$j \in I$.  If $V_j$ is equipped with the topology determined by
$\mathcal{N}_j$ for each $j \in I$, then it follows that the topology
on $V$ associated to $\mathcal{N}$ is the same as the corresponding
product topology.

\section{Bounded semimetrics}
\label{bounded semimetrics}
\setcounter{equation}{0}

        Let $X$ be a set, and let $d(x, y)$ be a $q$-semimetric on $X$
for some $q > 0$.  If $r_0$ is any positive real number, then
it is easy to see that
\begin{equation}
\label{d'(x, y) = min(d(x, y), r_0)}
        d'(x, y) = \min(d(x, y), r_0)
\end{equation}
is also a $q$-semimetric on $X$.  We also have that
\begin{equation}
\label{B_{d'}(x, r) = B_d(x, r)}
        B_{d'}(x, r) = B_d(x, r)
\end{equation}
for every $x \in X$ and $0 < r \le r_0$, and that
\begin{equation}
\label{B_{d'}(x, r) = X}
        B_{d'}(x, r) = X
\end{equation}
for every $x \in X$ when $r > r_0$, where the open balls are defined
as in (\ref{B(x, r) = B_d(x, r) = {y in X : d(x, y) < r}}).  This
implies that $d(x, y)$ and $d'(x, y)$ determine the same topologies on
$X$, since only balls of small radius are important for the topology.
Similarly, if $X$ is already equipped with a topology $\tau$, then
continuity properties of $d(x, y)$ with respect to $\tau$ as in
Section \ref{continuity of semimetrics} are equivalent to the
analogous continuity properties of (\ref{d'(x, y) = min(d(x, y),
  r_0)}) with respect to $\tau$.

        Let $A$ be a commutative group, and suppose that $d(x, y)$ is
a $q$-semimetric on $A$ for some $q > 0$ that is invariant under
translations on $A$.  In this case, (\ref{d'(x, y) = min(d(x, y),
  r_0)}) is also invariant under translations on $A$ for every $r_0 >
0$.  If $A$ is a commutative topological group, then continuity
properties of $d(x, y)$ on $A$ as in Section \ref{continuity of
  semimetrics} can be reformulated as continuity conditions at $0$, as
in Section \ref{translation-invariant semimetrics}.  These continuity
conditions are equivalent to their analogues for (\ref{d'(x, y) =
  min(d(x, y), r_0)}), as in the preceding paragraph.  Note that
\begin{equation}
\label{d(x, y) = d(x - (x + y), y - (x + y)) = d(-y, -x) = d(-x, -y)}
        d(x, y) = d(x - (x + y), y - (x + y)) = d(-y, -x) = d(-x, -y)
\end{equation}
for every $x, y \in A$, using invariance under translations in the
first step.

        Let $k$ be a field with a $q$-absolute value function $|\cdot|$
for some $q > 0$, and let $V$ be a vector space over $k$.  Also let
$N$ be a $q$-seminorm on $V$, and let $d(v, w)$ be the associated
$q$-semimetric, as in (\ref{d(v, w) = N(v - w), 3}).  Thus
\begin{equation}
\label{d(t v, t w) = N(t v - t w) = |t| N(v - w) = |t| d(v, w)}
 d(t \, v, t \, w) = N(t \, v - t \, w) = |t| \, N(v - w) = |t| \, d(v, w)
\end{equation}
for every $v, w \in V$ and $t \in k$, using the homogeneity of $N$ in
the second step.  In particular, this implies that
\begin{equation}
\label{d(t v, t w) le d(v, w)}
        d(t \, v, t \, w) \le d(v, w)
\end{equation}
for every $v, w \in V$ when $|t| \le 1$, with equality when $|t| = 1$.
If $d'(v, w)$ is defined as in (\ref{d'(x, y) = min(d(x, y), r_0)})
for some $r_0 > 0$, then it follows that
\begin{equation}
\label{d'(t v, t w) le d'(v, w)}
        d'(t \, v, t \, w) \le d'(v, w)
\end{equation}
for every $v, w \in V$ when $t \in k$ satisfies $|t| \le 1$, with
equality when $|t| = 1$.

        Now let $d(v, w)$ be any $q$-semimetric on $V$ that is
invariant under translations.  Observe that $d(v, w)$ satisfies
(\ref{d(t v, t w) le d(v, w)}) when $|t| \le 1$ exactly when open and
closed balls in $V$ with respect to $d$ centered at $0$ are balanced
subsets of $V$.  If $d(v, w)$ does not already have this property,
then one may wish to replace $d(v, w)$ with
\begin{equation}
\label{widetilde{d}(v, w) = sup {d(t v, t w) : t in k, |t| le 1}}
 \widetilde{d}(v, w) = \sup \{d(t \, v, t \, w) : t \in k, \ |t| \le 1\},
\end{equation}
at least when the supremum is finite.  In this case, it is easy to see
that (\ref{widetilde{d}(v, w) = sup {d(t v, t w) : t in k, |t| le 1}})
is also a $q$-semimetric on $V$ that is invariant under translations.
By construction,
\begin{equation}
\label{widetilde{d}(a v, a w) le widetilde{d}(v, w)}
        \widetilde{d}(a \, v, a \, w) \le \widetilde{d}(v, w)
\end{equation}
for every $v, w \in V$ and $a \in k$ with $|a| \le 1$, and
\begin{equation}
\label{d(v, w) le widetilde{d}(v, w)}
        d(v, w) \le \widetilde{d}(v, w)
\end{equation}
for every $v, w \in V$.  Of course, if $d(v, w)$ is bounded on $V$,
then the supremum in (\ref{widetilde{d}(v, w) = sup {d(t v, t w) : t
    in k, |t| le 1}}) is finite, with the same upper bounds as for
$d(v, w)$.  Otherwise, one can first replace $d(v, w)$ with
(\ref{d'(x, y) = min(d(x, y), r_0)}) for some $r_0 > 0$, to ensure
that it is bounded.

        Suppose that $V$ is a topological vector space with respect
to $|\cdot|$ on $k$, and that $d(v, w)$ is compatible with the
topology on $V$, in the sense that $d(0, w)$ is continuous as a
real-valued function of $w \in V$ at $0$.  Equivalently, this means
that for each $\epsilon > 0$, $0$ is an element of the interior of
$B_d(0, \epsilon)$ in $V$, where $B_d(0, \epsilon)$ is the open ball
in $V$ centered at $0$ with radius $\epsilon$ with respect to $d(v,
w)$, as in (\ref{B(x, r) = B_d(x, r) = {y in X : d(x, y) < r}}).  This
is the same as saying that for each $\epsilon > 0$ there is an open
set $W(\epsilon) \subseteq V$ such that $0 \in W(\epsilon)$ and
\begin{equation}
\label{W(epsilon) subseteq B_d(0, epsilon)}
        W(\epsilon) \subseteq B_d(0, \epsilon).
\end{equation}
If $|\cdot|$ is nontrivial on $k$, then for each $\epsilon > 0$ there
is an open set $U(\epsilon) \subseteq V$ such that $0 \in U(\epsilon)$
and
\begin{equation}
\label{t U(epsilon) subseteq W(epsilon)}
        t \, U(\epsilon) \subseteq W(\epsilon)
\end{equation}
for every $t \in k$ with $|t| \le 1$, as in (\ref{t widetilde{U}
  subseteq W}).  It follows that
\begin{equation}
\label{d(0, t u) < epsilon}
        d(0, t \, u) < \epsilon
\end{equation}
for every $\epsilon > 0$, $u \in U(\epsilon)$, and $t \in k$ with $|t|
\le 1$, by combining (\ref{W(epsilon) subseteq B_d(0, epsilon)}) and
(\ref{t U(epsilon) subseteq W(epsilon)}).

        Let $\widetilde{d}(v, w)$ be as in
(\ref{widetilde{d}(v, w) = sup {d(t v, t w) : t in k, |t| le 1}}),
and suppose that the supremum is finite for every $v, w \in V$.  As
before, this can always be arranged by replacing $d(v, w)$ with
(\ref{d'(x, y) = min(d(x, y), r_0)}) for some $r_0 > 0$, which would
not affect the continuity properties of $d(v, w)$.  Using (\ref{d(0, t
  u) < epsilon}), we get that
\begin{equation}
\label{widetilde{d}(0, u) le epsilon}
        \widetilde{d}(0, u) \le \epsilon
\end{equation}
for every $\epsilon > 0$ and $u \in U(\epsilon)$.  This shows that
$\widetilde{d}(0, u)$ is continuous as a real-valued function of $u
\in V$ at $0$ under these conditions.  Remember that this implies
additional continuity properties of $\widetilde{d}(v, w)$ on $V$, as
in Sections \ref{continuity of semimetrics} and
\ref{translation-invariant semimetrics}.  If $|\cdot|$ is the trivial
absolute value function on $k$, then (\ref{t widetilde{U} subseteq W})
may be considered as an additional hypothesis on the topology of $V$.
With this additional hypothesis, the rest of the previous argument
goes through as before.

        As mentioned in Section \ref{translation-invariant semimetrics},
there are well known methods for constructing translation-invariant
semimetrics on a commutative topological group that are compatible
with the given topology.  In particular, if $V$ is a topological
vector space with respect to $|\cdot|$ on $k$, as before, then this
can be applied to $V$ as a commutative topological group with respect
to addition.  Under suitable conditions, the arguments in the
preceding paragraphs can be used to get semimetrics on $V$ that also
satisfy (\ref{widetilde{d}(a v, a w) le widetilde{d}(v, w)}).
However, translation-invariant semimetrics on $V$ with this property
are typically obtained directly from the original construction under
the same conditions.  The main point is to use balanced neighborhoods
of $0$ in $V$ in the construction of the semimetrics, when $|\cdot|$
is nontrivial on $k$, or when $|\cdot|$ is trivial on $k$ and the
topology on $V$ satisfies (\ref{t widetilde{U} subseteq W}).

\section{Metrization}
\label{metrization}
\setcounter{equation}{0}

        Let $\mathcal{M}$ be a nonempty collection of $q$-seminorms
on a set $X$ for some $q > 0$.  Although we have often allowed $q$ to
depend on the element of $\mathcal{M}$, it is better to ask that each
element of $\mathcal{M}$ correspond to the same $q$ here.  As in
Section \ref{q-semimetrics}, this can always be arranged by taking
suitable powers anyway.  If there is a positive lower bound $q_0$ for
the $q$'s associated to elements of $\mathcal{M}$, then we can simply
treat each element of $\mathcal{M}$ as $q_0$-semimetric on $X$.  This
uses the fact that a $q$-semimetric on $X$ is also a $q_0$-semimetric
when $0 < q_0 \le q$, as in Section \ref{q-semimetrics}.

        If there are only finitely many elements of $\mathcal{M}$, then
their maximum is a $q$-semimetric on $X$ that determines the same
topology on $X$, as in Section \ref{collections of semimetrics}.  If
$\mathcal{M}$ is also nondegenerate on $X$, then the maximum of the
elements of $\mathcal{M}$ is a $q$-metric on $X$.

        Suppose now that $\mathcal{M}$ is countably infinite, and let
$d_j(x, y)$ with $j \in {\bf Z}_+$ be an enumeration of the elements of
$\mathcal{M}$.  Put
\begin{equation}
\label{d_j'(x, y) = min(d_j(x, y), 1/j)}
        d_j'(x, y) = \min(d_j(x, y), 1/j)
\end{equation}
for every $x, y \in X$ and $j \in {\bf Z}_+$, and
\begin{equation}
\label{d(x, y) = max_{j ge 1} d_j'(x, y)}
        d(x, y) = \max_{j \ge 1} d_j'(x, y)
\end{equation}
for every $x, y \in X$.  More precisely, $d(x, y) = 0$ when $d_j'(x,
y) = 0$ for every $j$.  Otherwise, the maximum in (\ref{d(x, y) =
  max_{j ge 1} d_j'(x, y)}) can be reduced to the maximum of finitely
many positive real numbers, because $d_j'(x, y) \le 1/j$ for every $x,
y \in X$ and $j \ge 1$, by construction.  As in the previous section,
$d_j'(x, y)$ is a $q$-semimetric on $X$ that determines the same
topology on $X$ as $d_j(x, y)$ for each $j \in {\bf Z}_+$.  Similarly,
one can check that $d(x, y)$ is a $q$-semimetric on $X$ under these
conditions, which is a $q$-metric on $X$ when $\mathcal{M}$ is
nondegenerate on $X$.  We would like to show that the topology
determined on $X$ by $d(x, y)$ is the same as the topology associated
to $\mathcal{M}$.

        Observe that
\begin{equation}
\label{B_d(x, r) = bigcap_{j = 1}^infty B_{d_j'}(x, r)}
        B_d(x, r) = \bigcap_{j = 1}^\infty B_{d_j'}(x, r)
\end{equation}
for every $x \in X$ and $r > 0$, by the definition (\ref{d(x, y) =
  max_{j ge 1} d_j'(x, y)}) of $d(x, y)$.  As usual, the open balls
are defined as in (\ref{B(x, r) = B_d(x, r) = {y in X : d(x, y) <
    r}}).  We also have that
\begin{equation}
\label{B_{d_j'}(x, r) = B_{d_j}(x, r)}
        B_{d_j'}(x, r) = B_{d_j}(x, r)
\end{equation}
when $r \le 1/j$, as in (\ref{B_{d'}(x, r) = B_d(x, r)}), and that
\begin{equation}
\label{B_{d_j'}(x, r) = X}
        B_{d_j'}(x, r) = X
\end{equation}
when $r > 1/j$, as in (\ref{B_{d'}(x, r) = X}).  If $r > 1$, then it
follows that (\ref{B_d(x, r) = bigcap_{j = 1}^infty B_{d_j'}(x, r)})
is equal to $X$ too.  Otherwise, if $r \le 1$, then let $l(r)$ be the
largest positive integer such that
\begin{equation}
\label{r le 1/l(r)}
        r \le 1/l(r).
\end{equation}
Thus (\ref{B_d(x, r) = bigcap_{j = 1}^infty B_{d_j'}(x, r)}) reduces
to
\begin{equation}
\label{B_d(x, r) = bigcap_{j = 1}^{l(r)} B_{d_j}(x, r)}
        B_d(x, r) = \bigcap_{j = 1}^{l(r)} B_{d_j}(x, r),
\end{equation}
by (\ref{B_{d_j'}(x, r) = B_{d_j}(x, r)}) and (\ref{B_{d_j'}(x, r) =
  X}).  Using this, one can verify that the topology determined on $X$
by $d(x, y)$ is the same as the one associated to $\mathcal{M}$, as
desired.

        If $A$ is a commutative group, then one can apply the
previous remarks to translation-invariant $q$-semimetrics on $A$, and
get translation-invariant $q$-semimetrics on $A$ as a result.  In
particular, if $A$ is a commutative topological group with a local
base for the topology at $0$ with only finitely or countably many
elements, then it is well known that there is a translation-invariant
semimetric on $A$ that determines the same topology on $A$.  More
precisely, if $A$ is any commutative topological group, then there is
a collection of translation-invariant semimetrics on $A$ that
determines the same topology on $A$, as mentioned in Section
\ref{translation-invariant semimetrics}.  If there is a local base for
the topology of $A$ at $0$ with only finitely or countably many
elements, then only finitely or countably many such semimetrics are
needed, which can be reduced to a single semimetric as before.
However, a single semimetric is often constructed directly in this
case, using the same types of arguments.  If $\{0\}$ is a closed set
in $A$, then this semimetric on $A$ is a metric.  Of course, if the
topology on any set $X$ is determined by a semimetric, then there is a
local base for the topology on $X$ at any point $x \in X$ with only
finitely or countably many elements, consisting of open balls centered
at $x$ with radius $1/j$ for $j \in {\bf Z}_+$, for instance.

        Let $k$ be a field with a $q$-absolute value function
$|\cdot|$ for some $q > 0$, and let $V$ be a topological vector
space over $k$ with respect to $|\cdot|$.  In particular, $V$ is a
commutative topological group with respect to addition, and so if
there is a local base for the topology of $V$ at $0$ with only
finitely or countably many elements, then there is a
translation-invariant semimetric $d(v, w)$ on $V$ that determines the
same topology on $V$, as in the preceding paragraph.  If $|\cdot|$ is
nontrivial on $k$, or if $|\cdot|$ is trivial on $V$ and the topology
on $V$ satisfies (\ref{t widetilde{U} subseteq W}), then one can
also choose $d(v, w)$ on $V$ so that
\begin{equation}
\label{d(t v, t w) le d(v, w), 2}
        d(t \, v, t \, w) \le d(v, w)
\end{equation}
for every $v, w \in V$ and $t \in k$ with $|t| \le 1$.  This is
typically obtained from the construction of the semimetric, following
the analogous construction for commutative topological groups, and
using the fact that one can choose a local base for the topology of
$V$ at $0$ consisting of balanced open sets under these conditions, as
in Section \ref{topological vector spaces}.  One can also get
(\ref{d(t v, t w) le d(v, w), 2}) as in the previous section.

        If the topology on $V$ is determined by a single $q$-seminorm
$N$, then one can simply use the semimetric associated to $N$ as in
(\ref{d(v, w) = N(v - w), 3}).  If the topology on $V$ is determined
by finitely many $q$-seminorms, then one can reduce to the case of a
single $q$-seminorm by taking their maximum.  Suppose now that for
each $j \in {\bf Z}_+$, $N_j$ is a $q_j$-seminorm on $V$ for some $q_j
> 0$, and that the topology on $V$ is determined by the $N_j$'s.  Thus
\begin{equation}
\label{N_j(v - w)}
        N_j(v - w)
\end{equation}
is a $q_j$-semimetric on $V$ for each $j \in {\bf Z}_+$, as in
(\ref{d(v, w) = N(v - w), 3}).  If there is a $q_0 > 0$ such that $q_j
\ge q_0$ for each $j \ge 1$, then (\ref{N_j(v - w)}) may be considered
as a $q_0$-semimetric on $V$ for each $j \ge 1$.  Otherwise, one can
replace (\ref{N_j(v - w)}) with
\begin{equation}
\label{N_j(v - w)^{q_j}}
        N_j(v - w)^{q_j}
\end{equation}
for each $j$, to get a sequence of semimetrics on $V$ that determines
the same topology on $V$.  In both cases, one can get a
translation-invariant $q_0$-semimetric $d(v, w)$ on $V$ for some $q_0
> 0$ that determines the same topology on $V$, as in (\ref{d(x, y) =
  max_{j ge 1} d_j'(x, y)}).  More precisely, one would first apply
(\ref{d_j'(x, y) = min(d_j(x, y), 1/j)}) to (\ref{N_j(v - w)}) or
(\ref{N_j(v - w)^{q_j}}) for each $j$, as appropriate, and then define
$d(v, w)$ as in (\ref{d(x, y) = max_{j ge 1} d_j'(x, y)}).  By
construction, $d(v, w)$ also satisfies (\ref{d(t v, t w) le d(v, w),
  2}), because (\ref{N_j(v - w)}) and (\ref{N_j(v - w)^{q_j}})
automatically have this property for each $j$, as in the previous
section.

\section{Totally bounded sets}
\label{totally bounded sets}
\setcounter{equation}{0}

        Let $X$ be a nonempty set, and let $d(x, y)$ be a $q$-semimetric
on $X$ for some $q > 0$.  A set $E \subseteq X$ is said to be
\emph{bounded}\index{bounded sets} with respect to $d$ if there
is an $x \in X$ and $r > 0$ such that
\begin{equation}
\label{E subseteq B_d(x, r)}
        E \subseteq B_d(x, r),
\end{equation}
where $B_d(x, r)$ is as in (\ref{B(x, r) = B_d(x, r) = {y in X : d(x,
    y) < r}}).  This implies that for each $x \in X$ there is an $r >
0$ such that (\ref{E subseteq B_d(x, r)}) holds, as one can check
using the $q$-semimetric version of the triangle inequality.  It
follows from this that if $E_1, \ldots, E_n$ are finitely many subsets
of $X$ that are bounded with respect to $d$, then their union
$\bigcup_{j = 1}^n E_j$ is bounded with respect to $E$ as well.  If $E
\subseteq X$ is compact with respect to the topology determined on $X$
by $d$, then it is easy to see that $E$ is bounded with respect to
$d$, by covering $E$ by balls of the form $B_d(x, r)$ for a fixed
$x \in X$ and with arbitrarily large $r$.

        A set $E \subseteq X$ is said to be
\emph{totally bounded}\index{totally bounded sets} with respect to $d$
if for every $r > 0$ there are finitely many elements $x_1, \ldots, x_n$
of $X$ such that
\begin{equation}
\label{E subseteq bigcup_{j = 1}^n B_d(x_j, r)}
        E \subseteq \bigcup_{j = 1}^n B_d(x_j, r).
\end{equation}
If this condition holds for any $r > 0$, then $E$ has to be bounded in
$X$ with respect to $d$, because the union of finitely many bounded
sets is still bounded, as in the preceding paragraph.  Similarly, the
union of finitely many subsets of $X$ that are totally bounded with
respect to $d$ is totally bounded with respect to $d$ too.  If $E
\subseteq X$ is compact with respect to the topology determined by
$d$, then $E$ is totally bounded with respect to $d$, as one can
verify by covering $E$ balls of radius $r$ for any $r > 0$.  If $E
\subseteq X$ is totally bounded with respect to $d$, then the closure
of $E$ with respect to the topology determined on $X$ by $d$ is
totally bounded with respect to $d$ as well.  If $d'$ is the
$q$-semimetric obtained from $d$ as in (\ref{d'(x, y) = min(d(x, y),
  r_0)}) for some $r_0 > 0$, then $d'$ is automatically bounded on
$X$, so that every subset of $X$ is bounded with respect to $d'$.
However, one can check that $E \subseteq X$ is totally bounded with
respect to $d'$ if and only if $E$ is totally bounded with respect to
$d$.

        If $E \subseteq X$ is nonempty and bounded with respect to $d$,
then the \emph{diameter}\index{diameters of sets} of $E$ with respect
to $d$ is defined by
\begin{equation}
\label{diam E = {diam}_d E = sup {d(x, y) : x, y in E}}
        \diam E = {\diam}_d E = \sup \{d(x, y) : x, y \in E\}.
\end{equation}
It is sometimes convenient to take this to be $+\infty$ when $E$ is
not bounded with respect to $d$, and to be $0$ when $E = \emptyset$.
One can check that $E \subseteq X$ is totally bounded with respect to
$d$ if and only if for every $t > 0$, $E$ is contained in the union of
finitely many subsets of $X$ with $d$-diameter less than or equal to
$t$.  More precisely, the ``if'' part of this statement uses the fact
that if $A \subseteq X$ has $d$-diameter less than $t$, then
\begin{equation}
\label{A subseteq B_d(x, t)}
        A \subseteq B_d(x, t)
\end{equation}
for every $x \in A$.  The converse uses the fact that any ball in $X$
of radius $r$ with respect to $d$ has $d$-diameter less than or equal
to $2^{1/q} \, r$, by the $q$-semimetric version of the triangle
inequality.

        Let $q_, \ldots, q_n$ be finitely many positive real numbers,
and let $d_j$ be a $q_j$-semimetric on $X$ for each $j = 1, \ldots,
n$.  Suppose that $E \subseteq X$ is totally bounded with respect to
$d_j$ for each $j = 1, \ldots, n$, and let $t_1, \ldots, t_n$ be
positive real numbers.  Thus for each $j = 1, \ldots, n$, $E$ can be
covered by finitely many subsets of $X$, each of which has
$d_j$-diameter less than or equal to $t_j$.  Using a common refinement
of these coverings, one can cover $E$ by finitely many subsets of $X$,
each of which has $d_j$-diameter less than or equal to $t_j$ for each
$j = 1, \ldots, n$ simultaneously.

\section{Totally bounded sets, continued}
\label{totally bounded sets, continued}
\setcounter{equation}{0}

        Let $A$ be a commutative topological group.  A set $E \subseteq A$
is said to be \emph{totally bounded}\index{totally bounded sets} in
$A$ if for each open set $U \subseteq A$ with $0 \in U$ there are
finitely many elements $x_1, \ldots, x_n$ of $A$ such that
\begin{equation}
\label{E subseteq bigcup_{j = 1}^n (x_j + U)}
        E \subseteq \bigcup_{j = 1}^n (x_j + U).
\end{equation}
As before, it is easy to see that the union of finitely many totally
bounded subsets of $A$ is also totally bounded, and that compact
subsets of $A$ are totally bounded.  If $E \subseteq A$ is totally
bounded, then every translate of $E$ in $A$ is totally bounded too, as
is $-E$.  If $\mathcal{B}_0$ is a local base for the topology of $A$
at $0$, then it suffices to show that $E \subseteq A$ can be covered
by finitely many translates of each $U \in \mathcal{B}_0$, in order to
verify that $E$ is totally bounded.

        Let $W$ be an open subset of $A$ that contains $0$, and let us
say that a set $C \subseteq A$ is \emph{$W$-small}\index{small sets}
if
\begin{equation}
\label{C - C subseteq W}
        C - C \subseteq W.
\end{equation}
Equivalently, this means that $x - y \in W$ for every $x, y \in C$,
which is the same as saying that
\begin{equation}
\label{C subseteq y + W}
        C \subseteq y + W
\end{equation}
for every $y \in C$.  Of course, if $C$ is $W$-small, then every
subset of $C$ is $W$-small as well.  Put
\begin{equation}
\label{widetilde{W} = W - W}
        \widetilde{W} = W - W,
\end{equation}
which is also an open subset of $A$ that contains $0$.  If (\ref{C
  subseteq y + W}) holds for any $y \in A$, then
\begin{equation}
\label{C - C subseteq W - W = widetilde{W}}
        C - C \subseteq W - W = \widetilde{W},
\end{equation}
so that $C$ is $\widetilde{W}$-small in $A$.

        Let $U \subseteq A$ be an open set that contains $0$ and
satisfies
\begin{equation}
\label{U - U subseteq W}
        U - U \subseteq W,
\end{equation}
which exists by the continuity of the group operations on $A$ at $0$.
If $E \subseteq A$ is totally bounded, then $E$ can be covered by
finitely many translates of $U$ in $A$, as in (\ref{E subseteq
  bigcup_{j = 1}^n (x_j + U)}).  Each translate of $U$ in $A$ is
$W$-small in $A$, so that $E$ can be covered by finitely many sets
which are $W$-small in $A$.  Conversely, if $E$ can be covered by
finitely many subsets of $A$ that are $W$-small, then $E$ can be
covered by finitely many translates of $W$.  Thus $E$ is totally
bounded in $A$ if and only if for every open set $W \subseteq A$ that
contains $0$, $E$ can be covered by finitely many subsets of $A$ that
are $W$-small.

        Let $B$ be a subgroup of $A$, which is also a commutative
topological group with respect to the topology induced by the one on
$A$.  If $W \subseteq A$ is an open set that contains $0$, then
\begin{equation}
\label{W_B = B cap W}
        W_B = B \cap W
\end{equation}
is a relatively open set in $B$ that contains $0$.  It is easy to see
that $C \subseteq B$ is $W$-small in $A$ if and only if $C$ is
$W_B$-small in $B$.  If $C \subseteq A$ is $W$-small, then $B \cap C$
is $W$-small in $A$ too, which implies that $B \cap C$ is $W_B$-small
in $B$.  This implies that $E \subseteq B$ can be covered by finitely
many subsets of $B$ that are $W_B$-small if and only if $E$ can be
covered by finitely many subsets of $A$ that are $W$-small.  It
follows that $E \subseteq B$ is totally bounded in $B$ if and only if
$E$ is totally bounded in $A$.  This also uses the fact that every
relatively open subset of $B$ that contains $0$ is of the form $W_B$
for some open set $W \subseteq A$ that contains $0$.

        Let $W_1, \ldots, W_n$ be finitely many open subsets of $A$
that contain $0$, so that
\begin{equation}
\label{W = bigcap_{j = 1}^n W_j}
        W = \bigcap_{j = 1}^n W_j
\end{equation}
is an open set that contains $0$ too.  Let $E \subseteq A$ be given,
and suppose that for each $j = 1, \ldots, n$, $E$ can be covered by
finitely many $W_j$-small subsets of $A$.  Using a common refinement
of these coverings, one can cover $E$ by finitely many $W$-small
sbsets of $A$.  Let $\mathcal{B}_0$ be a local sub-base for the
topology of $A$ at $0$, so that every open set in $A$ that contains
$0$ also contains the intersection of finitely many elements of
$\mathcal{B}_0$.  If for each $U \in \mathcal{B}_0$, $E$ can be
covered by finitely many $U$-small subsets of $A$, then it follows
that $E$ is totally bounded in $A$.

        Suppose for the moment that the topology on $A$ is determined by
a collection $\mathcal{M}$ of translation-invariant $q$-semimetrics on
$A$, where $q > 0$ is allowed to depend on the element of $\mathcal{M}$.
If $E \subseteq A$ is totally bounded as a subset of $A$ as a commutative
topological group, then it is easy to see that $E$ is totally bounded
with respect to every element of $\mathcal{M}$.  This uses the fact that
open balls in $A$ with respect to elements of $\mathcal{M}$ are open
sets in $A$, and that
\begin{equation}
\label{B_d(x, r) = x + B_d(0, r)}
        B_d(x, r) = x + B_d(0, r)
\end{equation}
for every $x \in A$, $r > 0$, and $d \in \mathcal{M}$, because $d$ is
supposed to be invariant under translations on $A$.  Conversely, if $E
\subseteq A$ is totally bounded with respect to every element of
$\mathcal{M}$, then one can check that $E$ is totally bounded in $A$
as a commutative topological group.  More precisely, let $d_1, \ldots,
d_n$ be finitely many elements of $\mathcal{M}$, and let $r_1, \ldots,
r_n$ be finitely many positive real numbers.  Thus
\begin{equation}
\label{U = bigcap_{j = 1}^n B_{d_j}(0, r_j)}
        U = \bigcap_{j = 1}^n B_{d_j}(0, r_j)
\end{equation}
is an open set in $A$ that contains $0$, and subsets of $A$ of this
form determine a local base for the topology of $A$ at $0$, as in
Section \ref{collections of semimetrics}.  If $E \subseteq A$ is
totally bounded with respect to $d_j$ for each $j = 1, \ldots, n$,
then $E$ can be covered by finitely many subsets of $X$, each of which
has $d_j$-diameter less than $r_j$ for every $j = 1, \ldots, n$, as in
the preceding section.  This implies that $E$ can be covered by
finitely many translates of (\ref{U = bigcap_{j = 1}^n B_{d_j}(0,
  r_j)}) in $A$, as desired.  Alternatively, this can be derived from
the remarks in paragraph, since the collection of open balls in $A$
centered at $0$ with respect to elements of $\mathcal{M}$ form a local
sub-base for the topology of $A$ at $0$.

        Suppose that $E \subseteq A$ is totally bounded in $A$ as a
commutative topological group, and let $W \subseteq A$ be an open set
that contains $0$.  Thus there are open sets $U_1, U_2 \subseteq A$
that contain $0$ and satisfy
\begin{equation}
\label{U_1 + U_2 subseteq W}
        U_1 + U_2 \subseteq W,
\end{equation}
by continuity of addition on $A$ at $0$.  Because $E$ is totally
bounded in $A$, there are finitely many elements $x_1, \ldots, x_n$ of
$A$ such that
\begin{equation}
\label{E subseteq bigcup_{j = 1}^n (x_j + U_1)}
        E \subseteq \bigcup_{j = 1}^n (x_j + U_1).
\end{equation}
This implies that
\begin{equation}
\label{overline{E} subseteq ... subseteq bigcup_{j = 1}^n (x_j + W)}
 \overline{E} \subseteq E + U_2 \subseteq \bigcup_{j = 1}^n (x_j + U_1 + U_2)
                                 \subseteq \bigcup_{j = 1}^n (x_j + W),
\end{equation}
where $\overline{E}$ is the closure of $E$ in $A$, and using
(\ref{overline{E} subseteq E + V}) in the first step.  It follows that
$\overline{E}$ is totally bounded in $A$ too.  Alternatively, let $U
\subseteq A$ be an open set that contains $0$ and satisfies
\begin{equation}
\label{overline{U} subseteq W, 2}
        \overline{U} \subseteq W,
\end{equation}
as in Section \ref{commutative topological groups}.  If $E$ is totally
bounded, then $E$ can be covered by finitely many translates of $U$,
which implies that $\overline{E}$ can be covered by finitely many
translates of $\overline{U}$, and hence by finitely many translates of
$W$.

        Suppose that $E_1, E_2 \subseteq A$ are totally bounded, and
let $W \subseteq A$ be an open set that contains $0$ again.  Thus
there are open sets $U_1, U_2 \subseteq A$ that contain $0$ and satisfy
(\ref{U_1 + U_2 subseteq W}), as before.  Because $E_1$ and $E_2$
are totally bounded, there are finitely many elements $x_1, \ldots, x_m$
and $y_1, \ldots, y_n$ of $A$ such that
\begin{equation}
\label{E_1 subseteq bigcup_{j = 1}^m (x_j + U_1) and ...}
 E_1 \subseteq \bigcup_{j = 1}^m (x_j + U_1) \quad\hbox{and}\quad
  E_2 \subseteq \bigcup_{l = 1}^n (y_l + U_2).
\end{equation}
This implies that
\begin{equation}
\label{E_1 + E_2 subseteq ...}
 E_1 + E_2 \subseteq \bigcup_{j = 1}^m \bigcup_{l = 1}^n (x_j + y_l + U_1 + U_2)
            \subseteq \bigcup_{j = 1}^m \bigcup_{l = 1}^n (x_j + y_l + W),
\end{equation}
using (\ref{U_1 + U_2 subseteq W}) in the second step.  It follows
that $E_1 + E_2$ is totally bounded in $A$ as well.

        Let $I$ be a nonempty set, and suppose that $A_j$ is a
commutative topological group for each $j \in I$.  As in Section
\ref{cartesian products},
\begin{equation}
\label{A = prod_{j in I} A_j, 2}
        A = \prod_{j \in I} A_j
\end{equation}
is also a commutative topological group, where the group operations
are defined coordinatewise, and using the product topology on $A$
associated to the given topology on $A_j$ for each $j \in I$.  If $E_j
\subseteq A_j$ is totally bounded for each $j \in I$, then one can
check that
\begin{equation}
\label{E = prod_{j in I} E_j}
        E = \prod_{j \in I} E_j
\end{equation}
is totally bounded in $A$, with respect to the product topology on
$A$.  More precisely, if $U_j \subseteq A_j$ is an open set that
contains $0$ for each $j \in I$, and if $U_j = A_j$ for all but
finitely many $j \in I$, then
\begin{equation}
\label{U = prod_{j in I} U_j}
        U = \prod_{j \in I} U_j
\end{equation}
is an open set in $A$ with respect to the product topology that
contains $0$.  In addition, the collection of subsets of $A$ of this
form is a local base for the product topology on $A$ at $0$.  In order
to show that $E \subseteq A$ is totally bounded, it suffices to verify
that $E$ can be covered by finitely many translates of $U$ in $A$, for
each $U \subseteq A$ of this form.  If $E_j \subseteq A_j$ is totally
bounded, then $E_j$ can be covered by finitely many translates of
$U_j$ in $A_j$.  If $E$ is as in (\ref{E = prod_{j in I} E_j}), then
one can cover $E$ by finitely many translates of $U$ in $A$, using the
coverings of $E_j$ by finitely many translates of $U_j$ in $A_j$ for
each of the finitely many $j \in I$ such that $U_j \ne A_j$.

\section{Bounded sets}
\label{bounded sets}
\setcounter{equation}{0}

        Let $k$ be a field, and let $|\cdot|$ be a nontrivial $q$-absolute
value function on $k$ for some $q > 0$.  Also let $V$ be a topological
vector space over $k$, with respect to $|\cdot|$ on $k$.  A set
$E \subseteq V$ is said to be \emph{bounded}\index{bounded sets}
if for each open set $U \subseteq V$ that contains $0$ there is a
$t_0 \in k$ such that
\begin{equation}
\label{E subseteq t_0 U}
        E \subseteq t_0 \, U.
\end{equation}
If $U$ is balanced in $V$, then it follows that
\begin{equation}
\label{E subseteq t U}
        E \subseteq t \, U
\end{equation}
for every $t \in k$ with $|t| \ge |t_0|$.  If $U$ is any open set in
$V$ that contains $0$, then there is a balanced open set $U_1
\subseteq A$ such that $0 \in U_1$ and $U_1 \subseteq U$, as in
Section \ref{topological vector spaces}.  If $E \subseteq V$ is a
bounded set, then there is a $t_1 \in k$ such that
\begin{equation}
\label{E subseteq t_1 U_1}
        E \subseteq t_1 \, U_1,
\end{equation}
as in (\ref{E subseteq t_0 U}).  This implies that
\begin{equation}
\label{E subseteq t U_1}
        E \subseteq t \, U_1
\end{equation}
for every $t \in k$ with $|t| \ge |t_1|$, as in (\ref{E subseteq t
  U}), because $U_1$ is balanced in $V$.  It follows that (\ref{E
  subseteq t U}) holds for every $t \in k$ with $|t| \ge |t_1|$,
because $U_1 \subseteq U$.

        Let $\mathcal{B}_0$ be a local base for the topology of
$V$ at $0$.  If $E \subseteq V$ has the property that for each
$U \in \mathcal{B}_0$ there is a $t_0 \in k$ that satisfies
(\ref{E subseteq t_0 U}), then it is easy to see that $E$ is
bounded in $V$.  In particular, one can take $\mathcal{B}_0$
to be the collection of all balanced open sets in $V$, which
basically corresponds to the discussion in the previous paragraph.

        Remember that open subsets of $V$ that contain $0$ are
absorbing, as in Section \ref{topological vector spaces}.  This
implies that subsets of $V$ with only one element are bounded, and in
fact that finite subsets of $V$ are bounded.  Similarly, one can check
that the union of finitely many bounded subsets of $V$ is bounded,
using the characterization of boundedness in terms of (\ref{E subseteq
  t U}) holding when $|t|$ is sufficiently large.

        Let $U \subseteq V$ be an open set that contains $0$.  If
$\{t_j\}_{j = 1}^\infty$ is a sequence of elements of $k$ such that
$|t_j| \to \infty$ as $j \to \infty$, then
\begin{equation}
\label{bigcup_{j = 1}^infty t_j U = V}
        \bigcup_{j = 1}^\infty t_j \, U = V,
\end{equation}
because $U$ is absorbing in $V$.  If $U$ is balanced in $V$, and if
$|t_j| \le |t_{j + 1}|$ for every $j$, then we get that
\begin{equation}
\label{t_j U subseteq t_{j + 1} U}
        t_j \, U \subseteq t_{j + 1} \, U
\end{equation}
for each $j$.  Let us also ask that $t_j \ne 0$ for each $j$, so that
$t_j \, U$ is an open set in $V$ for every $j$.  In particular, if $t
\in k$ and $|t| > 1$, then these conditions hold with $t_j = t^j$.  If
$E \subseteq V$ is compact, then (\ref{bigcup_{j = 1}^infty t_j U =
  V}) and (\ref{t_j U subseteq t_{j + 1} U}) imply that
\begin{equation}
\label{E subseteq t_j U}
        E \subseteq t_j \, U
\end{equation}
for some $j$.  This implies that $E$ is bounded in $V$, since we can
restrict our attention to balanced open subsets $U$ of $V$ in the
definition of boundedness.

        Suppose now that $E \subseteq V$ is totally bounded, as a subset
of $V$ as a commutative topological group with respect to addition.
Also let $W \subseteq V$ be an open set that contains $0$, and let
$U_1, U_2 \subseteq V$ be open sets that contain $0$ and satisfy
\begin{equation}
\label{U_1 + U_2 subseteq W, 2}
        U_1 + U_2 \subseteq W,
\end{equation}
as in (\ref{U + V subseteq W}).  We may as well ask that $U_1$, $U_2$
be balanced in $V$ too, since otherwise they can be replaced by
balanced open subsets, as in Section \ref{topological vector spaces}.
Because $E$ is totally bounded in $V$, there are finitely many
elements $v_1, \ldots, v_n$ of $V$ such that
\begin{equation}
\label{E subseteq bigcup_{j = 1}^n (v_j + U_2)}
        E \subseteq \bigcup_{j = 1}^n (v_j + U_2),
\end{equation}
as in (\ref{E subseteq bigcup_{j = 1}^n (x_j + U)}).  This implies
that
\begin{equation}
\label{E subseteq t U_1 + U_2}
        E \subseteq t \, U_1 + U_2
\end{equation}
for every $t \in k$ such that $|t|$ is sufficiently large, because
$U_1$ is absorbing in $V$.  Note that $U_2 \subseteq t \, U_2$ when
$|t| \ge 1$, since $U_2$ is balanced in $V$.  It follows that
\begin{equation}
\label{E subseteq t U_1 + t U_2 subseteq t W}
        E \subseteq t \, U_1 + t \, U_2 \subseteq t \, W
\end{equation}
when $|t|$ is sufficiently large, using (\ref{U_1 + U_2 subseteq W, 2})
in the second step.  This shows that totally bounded subsets of $V$
are bounded in $V$.

        Let $E \subseteq V$ be a bounded set, let $W \subseteq V$ be an
open set that contains $0$, and let $U_1, U_2 \subseteq V$ be open sets
that contain $0$ and satisfy (\ref{U_1 + U_2 subseteq W, 2}).  Thus
\begin{equation}
\label{E subseteq t U_1, 2}
        E \subseteq t \, U_1
\end{equation}
for every $t \in k$ such that $|t|$ is sufficiently large, because $E$
is bounded in $V$.  We also have that
\begin{equation}
\label{overline{E} subseteq E + t U_2}
        \overline{E} \subseteq E + t \, U_2
\end{equation}
for every $t \in k$ with $t \ne 0$, where $\overline{E}$ is the
closure of $E$ in $V$.  This follows from (\ref{overline{E} subseteq E
  + V}) and the fact that $t \, U_2$ is an open set in $V$ when $t \ne
0$.  Combining (\ref{E subseteq t U_1, 2}) and (\ref{overline{E}
  subseteq E + t U_2}), we get that
\begin{equation}
\label{overline{E} subseteq t U_1 + t U_2 = t W}
        \overline{E} \subseteq t \, U_1 + t \, U_2 = t \, W
\end{equation}
for every $t \in k$ such that $|t|$ is sufficiently large, which
implies that $\overline{E}$ is bounded in $V$ too.

        Let $E_1, E_2 \subseteq V$ be bounded sets, and let us check
that $E_1 + E_2$ is also a bounded subset of $V$.  To do this, let $W
\subseteq V$ be an open set that contains $0$, and let $U_1, U_2
\subseteq V$ be open sets that contain $0$ and satisfy (\ref{U_1 + U_2
  subseteq W, 2}).  Thus
\begin{equation}
\label{E_1 subseteq t U_1 and E_2 subseteq t U_2}
        E_1 \subseteq t \, U_1 \quad\hbox{and}\quad E_2 \subseteq t \, U_2
\end{equation}
when $t \in k$ and $|t|$ is sufficiently large, which implies that
\begin{equation}
\label{E_1 + E_2 subseteq t U_1 + t U_2 = t W}
        E_1 + E_2 \subseteq t \, U_1 + t \, U_2 = t \, W
\end{equation}
when $|t|$ is suffiently large, as desired.

        Let $I$ be a nonempty set, and let $V_j$ be a topological
vector space over $k$ for each $j \in I$.  Remember that
\begin{equation}
\label{V = prod_{j in I} V_j, 2}
        V = \prod_{j \in I} V_j
\end{equation}
is a topological vector space over $k$ as well, as in Section
\ref{cartesian products}, where the vector space operations are
defined coordinatewise, and using the product topology on $V$
associated to the given topology on $V_j$ for each $j \in I$.  Suppose
that $E_j \subseteq V_j$ is a bounded set for each $j \in I$, and let
us check that
\begin{equation}
\label{E = prod_{j in I} E_j, 2}
        E = \prod_{j \in I} E_j
\end{equation}
is bounded in $V$.  To do this, let $U_j \subseteq V_j$ be an open set
that contains $0$ for each $j \in I$, and suppose that $U_j = V_j$ for
all but finitely many $j \in I$.  Thus
\begin{equation}
\label{U = prod_{j in I} U_j, 2}
        U = \prod_{j \in I} U_j
\end{equation}
is an open set in $V$ with respect to the product topology that
contains $0$, and the collection of subsets of $V$ of this form is a
local base for the topology of $V$ at $0$.  In order to show that
$E$ is bounded in $V$, it suffices to verify that
\begin{equation}
\label{E subseteq t U, 2}
        E \subseteq t \, U
\end{equation}
for every $t \in k$ such that $|t|$ is sufficiently large.  Of course,
for each $j \in I$, we know that
\begin{equation}
\label{E_j subseteq t U_j}
        E_j \subseteq t \, U_j
\end{equation}
for every $t \in k$ such that $|t|$ is sufficiently large, because
$E_j$ is bounded in $V_j$.  By hypothesis, $U_j = V_j$ for all but
finitely many $j \in I$, in which case (\ref{E_j subseteq t U_j})
holds for every $t \in k$ with $t \ne 0$.  If $|t|$ is sufficiently
large, then it follows that (\ref{E_j subseteq t U_j}) holds for all
$j \in I$ simultaneously, which implies that (\ref{E subseteq t U, 2})
holds, as desired.

        Let $V$ be any topological vector space over $k$ again, and
let $N$ be a nonnegative real-valued function on $V$ that satisfies
the usual homogeneity condition (\ref{N(t v) = |t| N(v), 2}).
If $N$ is continuous at $0$, then there is an open set $U \subseteq V$
that contains $0$ and satisfies
\begin{equation}
\label{U subseteq B_N(0, 1)}
        U \subseteq B_N(0, 1),
\end{equation}
where $B_N(0, r)$ is as in (\ref{B_N(0, r) = {v in V : N(v) < r}}).
Conversely, this condition implies that $N$ is continuous at $0$,
because of (\ref{t B_N(0, r) = B_N(0, |t| r)}), and our standing
hypothesis in this section that $|\cdot|$ be nontrivial on $k$.  If $E
\subseteq V$ is bounded, then it follows that $N$ is bounded on $E$ in
this situation.  If $N$ is a $q$-seminorm on $V$ for some $q > 0$,
then continuity of $N$ at $0$ on $V$ implies that $N$ is continuous on
$V$, as in Section \ref{compatible seminorms}.  In particular, this
implies that open balls in $V$ with respect to $N$ are open sets in
$V$, as before.  In this case, the boundedness of $N$ on $E$ is the
same as the boundedness of $E$ with respect to the corresponding
$q$-semimetric (\ref{d(v, w) = N(v - w), 3}), as in the previous
section.

        Now let $V$ be a vector space over $k$, and let $\mathcal{N}$
be a collection of $q$-seminorms on $V$.  As usual, we can let $q$
depend on the element of $\mathcal{N}$, as long as $|\cdot|$ is a
$q$-absolute value function on $k$.  As in Section \ref{compatible
seminorms}, $V$ is a topological vector space with respect to the
topology associated to $\mathcal{N}$.  If $E \subseteq V$ is bounded
with respect to the topology associated to $\mathcal{N}$, then each
element of $\mathcal{N}$ is bounded on $E$.  This follows from the
remarks in the preceding paragraph, since every element of $\mathcal{N}$
is continuous with respect to the corresponding topology on $V$.
Conversely, suppose that $E \subseteq V$ has the property that every
element of $\mathcal{N}$ is bounded on $E$.  Let $U \subseteq V$ be
an open set with respect to the topology associated to $\mathcal{M}$
that contains $0$.  By construction, this means that there are finitely
many elements $N_1, \ldots, N_l$ of $\mathcal{N}$ and finitely many
positive real numbers $r_1, \ldots, r_l$ such that
\begin{equation}
\label{bigcap_{j = 1}^l B_{N_j}(0, r_j) subseteq U}
        \bigcap_{j = 1}^l B_{N_j}(0, r_j) \subseteq U.
\end{equation}
Using this, it is easy to see that $E \subseteq t \, U$ for every $t
\in k$ such that $|t|$ is sufficiently large, because $N_j$ is bounded
on $E$ for each $j = 1, \ldots, l$.  This implies that $E$ is bounded
with respect to the topology on $V$ associated to $\mathcal{N}$.

        Let $V$ be any topological vector space over $k$, and suppose
that $E \subseteq V$ is bounded.  If $\{v_j\}_{j = 1}^\infty$ is a
sequence of elements of $E$, and if $\{t_j\}_{j = 1}^\infty$ is a
sequence of elements of $k$ that converges to $0$ with respect to
$|\cdot|$, then it is easy to see that $\{t_j \, v_j\}_{j = 1}^\infty$
converges to $0$ in $V$.  Conversely, suppose that $E$ is not bounded
in $V$, so that there is an open set $U \subseteq V$ such that $0 \in
U$ and $E \not\subseteq t \, U$ for every $t \in k$.  Let $t_0$ be an
element of $k$ such that $|t_0| > 1$, which exists because $|\cdot|$
is nontrivial on $k$.  Thus $E \not\subseteq t_0^j \, U$ for every
positive integer $j$, and so we can choose $v_j \in E \setminus (t_0^j
\, U)$ for each $j$.  This implies that $t_0^{-j} \, v_j \not\in U$
for each $j$, so that $\{t_0^{-j} \, v_j\}_{j = 1}^\infty$ does not
converge to $0$ in $V$.  Of course, $\{t_0^{-j}\}_{j = 1}^\infty$ does
converge to $0$ in $k$, since $|t_0| > 1$.

\section{Continuous functions}
\label{continuous functions}
\setcounter{equation}{0}

        If $X$ and $Y$ are topological spaces, then we let
$C(X, Y)$\index{C(X, Y)@$C(X, Y)$} be the space of continuous
mappings from $X$ into $Y$.  Let $X$ be a nonempty topological space,
and let $k$ be a field with a $q$-absolute value function $|\cdot|$
for some $q > 0$.  This leads to a topology on $k$ in the usual way,
corresponding to the $q$-metric (\ref{d(x, y) = |x - y|, 2})
associated to $|\cdot|$ on $k$.  Note that $C(X, k)$ is a vector space
over $k$ with respect to pointwise addition and scalar multiplication.
More precisely, $C(X, k)$ is a commutative algebra over $k$ with
respect to pointwise multiplication of functions.

        Let $E$ be a nonempty compact subset of $X$, and let $f$
be a continuous $k$-valued function on $X$.  Thus $f(E)$ is a nonempty
compact subset of $k$, which is bounded with respect to $|\cdot|$ in
particular.  This permits us to put
\begin{equation}
\label{||f||_E = sup_{x in E} |f(x)|}
        \|f\|_E = \sup_{x \in E} |f(x)|,
\end{equation}
since the right side of (\ref{||f||_E = sup_{x in E} |f(x)|}) is the
supremum of a nonempty set of nonnegative real numbers.  Remember that
$|\cdot|$ is continuous as a real-valued function on $k$ with respect
to the topology determined on $k$ by the corresponding $q$-metric, as
in Section \ref{continuity of semimetrics}.  This implies that
$|f(x)|$ is continuous as a real-valued function on $X$, so that the
supremum in the right side of (\ref{||f||_E = sup_{x in E} |f(x)|}) is
attained.

        It is easy to see that (\ref{||f||_E = sup_{x in E} |f(x)|})
defines a $q$-seminorm on $C(X, k)$ for every nonempty compact set $E
\subseteq X$.  This is the \emph{supremum $q$-seminorm}\index{supremum
  seminorms} associated to $E$.  If $X$ is compact, then we can take
$E = X$ in (\ref{||f||_E = sup_{x in E} |f(x)|}), which defines a
$q$-norm on $C(X, k)$.  We also have that
\begin{equation}
\label{||f g||_E le ||f||_E ||g||_E}
        \|f \, g\|_E \le \|f\|_E \, \|g\|_E
\end{equation}
for every $f, g \in C(X, k)$ and nonempty compact set $E \subseteq X$,
because of the multiplicative property of absolute value functions.

        Let $\mathcal{N}$ be the collection of supremum $q$-seminorms
on $C(X, k)$ that are associated to nonempty compact subsets $E$ of
$X$ as in (\ref{||f||_E = sup_{x in E} |f(x)|}).  This is a
nondegenerate collection of $q$-seminorms on $C(X, k)$, because finite
subsets of $X$ are compact.  As in Section \ref{compatible seminorms},
$C(X, k)$ is a Hausdorff topological vector space over $k$ with
respect to the topology determined by $\mathcal{N}$.  Similarly, one
can check that multiplication of functions defines a continuous
mapping from $C(X, k) \times C(X, k)$ into $C(X, k)$, where $C(X, k)
\times C(X, k)$ is equipped with the product topology associated to
the topology just defined on $C(X, k)$, using (\ref{||f g||_E le
  ||f||_E ||g||_E}).  If $X$ is compact, then the same topology on
$C(X, k)$ is determined by the supremum $q$-norm $\|f\|_X$.

        Observe that $C(X, k \setminus \{0\})$ is a commutative group
with respect to multiplication.  In fact, $C(X, k \setminus \{0\})$ is
a commutative topological group, with respect to the topology induced
by the one on $C(X, k)$ described in the preceding paragraph.
Continuity of multiplication on $C(X, k \setminus \{0\})$ follows from
the analogous statement for $C(X, k)$ just mentioned.  The remaining
point is that $f \mapsto 1/f$ defines a continuous mapping on $C(X, k
\setminus \{0\})$ with respect to the induced topology, which can be
verified using standard arguments.  In particular, if $f$ is a continuous
function on $X$ with values in $k \setminus \{0\}$, then $1/f$ is bounded
on compact subsets of $X$, since it is continuous on $X$ too.

        Put
\begin{equation}
\label{{bf T}_k = {x in k : |x| = 1}}\index{T_k@${\bf T}_k$}
        {\bf T}_k = \{x \in k : |x| = 1\},
\end{equation}
which is a subgroup of $k \setminus \{0\}$ with respect to
multiplication.  Thus $C(X, {\bf T}_k)$ is a subgroup of $C(X, k
\setminus \{0\})$ with respect to multiplication of functions.  It
follows that $C(X, {\bf T}_k)$ is also a commutative topological group
with respect to the topology induced by the one on $C(X, k)$ mentioned
earlier.  Note that ${\bf T}_k$ is a closed set in $k$ with respect to
the topology determined by the $q$-metric associated to $|\cdot|$.
Using this, one can check that $C(X, {\bf T}_k)$ is a closed set in
$C(X, k)$, with respect to the usual topology.

        Let $E$ be a nonempty compact subset of $X$, and observe that
\begin{equation}
\label{||a f||_E = ||f||_E}
        \|a \, f\|_E = \|f\|_E
\end{equation}
for every $a \in C(X, {\bf T}_k)$ and $f \in C(X, k)$.  If
\begin{equation}
\label{d_E(f, g) = ||f - g||_E}
        d_E(f, g) = \|f - g\|_E
\end{equation}
is the $q$-semimetric associated to $\|\cdot\|_E$, then we get that
\begin{equation}
\label{d_E(a f, a g) = d_E(f, g)}
        d_E(a \, f, a \, g) = d_E(f, g)
\end{equation}
for every $a \in C(X, {\bf T}_k)$ and $f, g \in C(X, k)$.  The
restriction of (\ref{d_E(f, g) = ||f - g||_E}) to $f$, $g$ in $C(X,
{\bf T}_k)$ defines a $q$-semimetric on $C(X, {\bf T}_k)$, and
(\ref{d_E(a f, a g) = d_E(f, g)}) implies that this semimetric is
invariant under translations on $C(X, {\bf T}_k)$, as a group with
respect to multiplication.  By construction, the usual topology on
$C(X, k)$ is the one determined by the collection of $q$-semimetrics
(\ref{d_E(f, g) = ||f - g||_E}) associated to nonempty compact sets $E
\subseteq X$.  The induced topology on $C(X, {\bf T}_k)$ is the same
as the one determined by the collection of restrictions of these
$q$-semimetrics to $C(X, {\bf T}_k)$, as in Section \ref{collections
  of semimetrics}.

        Suppose for the moment that $X$ is any nonempty set equipped
with the discrete topology.  This implies that every function on $X$
is continuous, and that every compact subset of $X$ has only finitely
many elements.  In this case, $C(X, k)$ can be identified with the
Cartesian product of a family of copies of $k$ indexed by $X$, as in
Section \ref{cartesian products}.  The topology on $C(X, k)$ described
earlier corresponds exactly to the product topology in this situation,
using the topology on each factor of $k$ determined by the $q$-metric
associated to $|\cdot|$.  Similarly, $C(X, k \setminus \{0\})$
and $C(X, {\bf T}_k)$ can be identified with Cartesian products of
copies of $k \setminus \{0\}$ and ${\bf T}_k$ indexed by $X$, with
their appropriate product topologies.

        Suppose now that $X$ is a nonempty topological space which is
locally compact,\index{locally compact topological spaces} in the
sense that every point in $X$ is contained in an open subset of $X$
that is contained in a compact subset of $X$.  This implies that every
compact subset of $X$ is contained in an open subset of $X$ that is
contained in a compact subset of $X$, by standard arguments.  Suppose
also that $X$ is $\sigma$-compact,\index{sigma-compact topological
  spaces@$\sigma$-compact topological spaces} which means that there
is a sequence $E_1, E_2, E_3, \ldots$ of compact subsets of $X$ such that
\begin{equation}
\label{bigcup_{j = 1}^infty E_j = X}
        \bigcup_{j = 1}^\infty E_j = X.
\end{equation}
We may as well suppose that $E_j \subseteq E_{j + 1}$ for each $j$,
since otherwise we can replace $E_j$ with $\bigcup_{l = 1}^j E_l$ for
each $j$.  Using local compactness, we can refine this a bit, to get
that for each $j \ge 1$ there is an open set $U_j \subseteq X$ such
that
\begin{equation}
\label{E_j subseteq U_j subseteq E_{j + 1}}
        E_j \subseteq U_j \subseteq E_{j + 1}.
\end{equation}
More precisely, if $E_j$ has already been chosen for some $j$, then we
take $U_j$ to be an open set in $X$ that contains $E_j$ and is
contained in a compact subset of $X$.  In the next step, we expand
$E_{j + 1}$ so that it contains $U_j$, and continue as before.

        It follows that $U_j \subseteq U_{j + 1}$ for each $j$, and that
\begin{equation}
\label{bigcup_{j = 1}^infty U_j = X}
        \bigcup_{j = 1}^\infty U_j = X.
\end{equation}
If $E$ is any compact subset of $X$, then $E$ is contained in the
union of finitely many $U_j$'s, and hence in a single $U_j$.  This
implies that
\begin{equation}
\label{E subseteq E_j}
        E \subseteq E_j
\end{equation}
for some $j$.  We may as well choose $E_1$ to be nonempty, so that
$E_j \ne \emptyset$ for every $j$, which means that the supremum
$q$-seminorm associated to $E_j$ may be defined on $C(X, k)$ for each
$j$ as before.  Under these conditions, the sequence of supremum
$q$-seminorms associated to the $E_j$'s are sufficient to define the
usual topology on $C(X, k)$.

\section{Continuous linear mappings}
\label{continuous linear mappings}
\setcounter{equation}{0}

        Let $X_1$ and $X_2$ be topological spaces, and let $\phi$ be
a mapping from $X_1$ into $X_2$.  As usual, $\phi$ is said to be
\emph{sequentially continuous}\index{sequential continuity} at a point
$x \in X_1$ if for each sequence $\{x_j\}_{j = 1}^\infty$ of elements
of $X_1$ that converges to $x$, $\{\phi(x_j)\}_{j = 1}^\infty$
converges to $\phi(x)$ in $X_2$.  Of course, if $\phi$ is continuous
at $x$, then $\phi$ is sequentially continuous at $x$.  If there is a
local base for the topology of $X_1$ at $x$ with only finitely or
countably many elements, and if $\phi$ is sequentially continuous at
$x$, then it is well known that $\phi$ is continuous at $x$.  To see
this, suppose for the sake of a contradiction that there is an open
set $W \subseteq X_2$ such that $\phi(x) \in W$, but
\begin{equation}
\label{phi(U) not subseteq W}
        \phi(U) \not\subseteq W
\end{equation}
for every open set $U \subseteq X_1$ that contains $x$.  The
hypothesis that there be a local base for the topology of $X_1$ at $x$
with only finitely or countably many elements means that there is a
sequence $U_1(x), U_2(x), U_3(x), \ldots$ of open subsets of $X_1$
that contain $x$ such that if $U$ is any other open subset of $X_1$
that contains $x$, then $U_l(x) \subseteq U$ for some $l \ge 1$.  We
may as well ask also that
\begin{equation}
\label{U_{l + 1}(x) subseteq U_l(x)}
        U_{l + 1}(x) \subseteq U_l(x)
\end{equation}
for every $l$, since otherwise we can replace $U_l(x)$ with
$\bigcap_{j = 1}^l U_j(x)$ for each $l$.  Our hypothesis (\ref{phi(U)
  not subseteq W}) implies that for each positive integer $l$ there is
an $x_l \in U_l(x)$ such that
\begin{equation}
\label{phi(x_l) not in W}
        \phi(x_l) \not\in W.
\end{equation}
Under these conditions, $\{x_l\}_{l = 1}^\infty$ converges to $x$ in
$X_1$, but $\{\phi(x_l)\}_{l = 1}^\infty$ does not converge to
$\phi(x)$ in $X_2$, as desired.

        Now let $A_1$, $A_2$ be commutative topological groups, and
let $\phi$ be a group homomorphism from $A_1$ into $A_2$.  If $\phi$
is continuous at $0$, then for each open set $W \subseteq A_2$ that
contains $0$ there is an open set $U \subseteq A_1$ that contains $0$
such that
\begin{equation}
\label{phi(U) subseteq W}
        \phi(U) \subseteq W.
\end{equation}
This implies that
\begin{equation}
\label{phi(x + U) subseteq phi(x) + W}
        \phi(x + U) \subseteq \phi(x) + W
\end{equation}
for every $x \in A_1$, and hence that $\phi$ is continuous at every
point in $A_1$.  More precisely, this implies that $\phi$ is uniformly
continuous as a mapping from $A_1$ into $A_2$ in a suitable sense.  In
particular, if $E \subseteq A_1$ is totally bounded, then one can use
this to check that $\phi(E)$ is totally bounded in $A_2$.  If $\phi$
is sequentially continuous at $0$, then it is easy to see that $\phi$
is sequentially continuous at every point in $A_1$.  If there is also
a local base for the topology of $A_1$ at $0$ with only finitely or
countably many elements, then $\phi$ is continuous at $0$, as in the
preceding paragraph.

        Let $k$ be a field, and let $|\cdot|$ be a $q$-absolute value
function on $k$ for some $q > 0$.  Remember that topological vector
spaces over $k$ are commutative topological groups with respect to
addition, and that linear mappings between vector space over $k$ are
group homomorphisms with respect to addition.  Thus the remarks in the
previous paragraph can be applied to linear mappings between
topological vector spaces over $k$.  Note that $k$ may be considered
as a one-dimensional topological vector space over itself, using the
topology determined on $k$ by the $q$-metric associated to $|\cdot|$.

        Let $V$ be a topological vector space over $k$, and let $v \in V$
be given.  Put
\begin{equation}
\label{phi_v(t) = t v}
        \phi_v(t) = t \, v
\end{equation}
for each $t \in k$, which defines a continuous linear mapping from $k$
into $V$, because of continuity of scalar multiplication on $V$, as in
Section \ref{topological vector spaces}.  Suppose that $v \ne 0$,
and put
\begin{equation}
\label{L_v = phi_v(k) = {t v : t in k}}
        L_v = \phi_v(k) = \{t \, v : t \in k\},
\end{equation}
which is the one-dimensional linear subspace of $V$ spanned by $v$.
Let us ask also that $\{0\}$ be a closed set in $V$, so that $\{v\}$
is a closed set in $V$ too, by continuity of translations.  Thus $V
\setminus \{v\}$ is an open set in $V$ that contains $0$.

        If $|\cdot|$ is nontrivial on $k$, then there is a balanced open
set $U \subseteq V$ contained in $V \setminus \{v\}$, as in Section
\ref{topological vector spaces}.  Because $U$ is balanced, we have
that
\begin{equation}
\label{t v not in U}
        t \, v \not\in U
\end{equation}
for every $t \in k$ with $|t| \ge 1$.  If $a \in k$ and $a \ne 0$,
then $a \, U$ is also a balanced open set in $V$ that contains $0$, and
\begin{equation}
\label{t v not in a U}
        t \, v \not\in a \, U
\end{equation}
for every $t \in k$ with $|t| \ge |a|$, by (\ref{t v not in U}).
Using this, one can check that $\phi_v$ is a homeomorphism from $k$
onto $L_v$ under these conditions, with respect to the topology
induced on $L_v$ by the one on $V$.  More precisely, (\ref{t v not in
  a U}) implies that the inverse of $\phi_v$ is continuous at $0$ as a
mapping from $L_v$ into $k$.

        Let us ask that $|\cdot|$ be nontrivial on $k$ for the rest
of the section, and let $V_1$ and $V_2$ be topological vector spaces
over $k$.  A linear mapping $\phi$ from $V_1$ into $V_2$ is said to
be \emph{bounded}\index{bounded linear mappings} if for each bounded
set $E \subseteq V_1$, we have that $\phi(E)$ is a bounded subset of $V_2$,
where boundedness is defined as in Section \ref{bounded sets}.  If $\phi$
is continuous, then it is easy to see that $\phi$ is bounded, directly
from the definitions.

        Suppose for the moment that $U$ is an open subset of $V_1$
that contains $0$.  If $\phi(U)$ is a bounded subset of $V_2$, then
one can check that $\phi$ is continuous.  In particular, this holds
when $U$ is also bounded in $V_1$, and $\phi$ is a bounded linear
mapping.  If there is a bounded open set $W \subseteq V_2$ that
contains $0$, and if $\phi$ is continuous, then there is an open set
$U \subseteq V_1$ that contains $0$ such that $\phi(U) \subseteq W$,
which implies that $\phi(U)$ is bounded in $V_2$.  If the topology on
a vector space $V$ over $k$ is determined by a single $q$-seminorm
$N$, then open balls with respect to $N$ are both bounded and open
in $V$.

        Suppose now that there is a local base for the topology of
$V_1$ at $0$ with only finitely or countably many elements.  Thus
there is a sequence $U_1, U_2, U_3, \ldots$ of open subsets of $V_1$
that contain $0$ such that any other open subset of $V_1$ that contains
$0$ also contains $U_l$ for some $l$.  We may as well take $U_l$ to be
balanced in $V_1$ for each $l$, because $|\cdot|$ is nontrivial on $k$,
as in Section \ref{topological vector spaces}.  We may also ask that
\begin{equation}
\label{U_{l + 1} subseteq U_l}
        U_{l + 1} \subseteq U_l
\end{equation}
for each $l$, since otherwise we can replace $U_l$ with $\bigcap_{j =
  1}^l U_j$ for each $l$, as before.  Note that $U_l$ is absorbing in
$V_1$ for each $l$, as in Section \ref{topological vector spaces}, and
using the nontriviality of $|\cdot|$ on $V$ again.  This permits us to
define
\begin{equation}
\label{N_l(v) = N_{U_l}(v)}
        N_l(v) = N_{U_l}(v)
\end{equation}
for each $v \in V_1$ and $l \ge 1$ as in (\ref{N_A(v) = inf {|t| : t
    in k, v in t A}}).  Using (\ref{U_{l + 1} subseteq U_l}), we get
that
\begin{equation}
\label{N_l(v) le N_{l + 1}(v)}
        N_l(v) \le N_{l + 1}(v)
\end{equation}
for every $v \in V_1$ and $l \ge 1$, as in (\ref{N_C(v) le N_B(v)}).

        Let $\{v_j\}_{j = 1}^\infty$ be a sequence of elements of $V_1$.
If $\{v_j\}_{j = 1}^\infty$ converges to $0$ in $V_1$, then 
\begin{equation}
\label{lim_{j to infty} N_l(v_j) = 0}
        \lim_{j \to \infty} N_l(v_j) = 0
\end{equation}
for every $l \ge 1$, because $t \, U_l$ is an open set in $V$ that
contains $0$ for each $l \ge 1$ and $t \in k$ with $t \ne 0$.
Conversely, if (\ref{lim_{j to infty} N_l(v_j) = 0}) holds for every
$l$, then $\{v_j\}_{j = 1}^\infty$ converges to $0$ in $V_1$, because
the $U_l$'s determine a local base for the topology of $V_1$ at $0$.
In fact, it suffices that the $U_l$'s and their dilates by nonzero
elements of $k$ determine a local base for the topology of $V_1$ at
$0$ for this to work.

        If $\{l_j\}_{j = 1}^\infty$ is a sequence of positive integers
such that
\begin{equation}
\label{l_j to infty as j to infty}
        l_j \to \infty \quad\hbox{as } j \to \infty
\end{equation}
and
\begin{equation}
\label{lim_{j to infty} N_{l_j}(v_j) = 0}
        \lim_{j \to \infty} N_{l_j}(v_j) = 0,
\end{equation}
then (\ref{lim_{j to infty} N_l(v_j) = 0}) holds for each $l$, because
of (\ref{N_l(v) le N_{l + 1}(v)}).  This implies that $\{v_j\}_{j =
  1}^\infty$ converges to $0$ in $V_1$, as before.  Conversely, if
$\{v_j\}_{j = 1}^\infty$ converges to $0$ in $V_1$, then there is a
sequence $\{l_j\}_{j = 1}^\infty$ of positive integers that satisfies
(\ref{l_j to infty as j to infty}) and (\ref{lim_{j to infty}
  N_{l_j}(v_j) = 0}).  Of course, this uses the fact that (\ref{lim_{j
    to infty} N_l(v_j) = 0}) holds for each $l$.

        If $\{v_j\}_{j = 1}^\infty$ converges to $0$ in $V_1$, then there
is a sequence $\{t_j\}_{j = 1}^\infty$ of nonzero elements of $k$ that
converges to $0$ in $k$ and has the property that
\begin{equation}
\label{lim_{j to infty} t_j^{-1} v_j = 0}
        \lim_{j \to \infty} t_j^{-1} \, v_j = 0
\end{equation}
in $V_1$.  More precisely, if $\{l_j\}_{j = 1}^\infty$ is a sequence
of positive integers that satisfies (\ref{l_j to infty as j to infty})
and (\ref{lim_{j to infty} N_{l_j}(v_j) = 0}), then it suffices to
show that there is a sequence $\{t_j\}_{j = 1}^\infty$ of nonzero
elements of $k$ that converges to $0$ and satisfies
\begin{equation}
\label{N_{l_j}(t_j^{-1} v_j) = |t_j|^{-1} N_{l_j}(v_j) to 0 as j to infty}
        N_{l_j}(t_j^{-1} \, v_j) = |t_j|^{-1} \, N_{l_j}(v_j) \to 0
                                             \quad\hbox{as } j \to \infty.
\end{equation}
In order to get such a sequence $\{t_j\}_{j = 1}^\infty$, one can use
the nontriviality of $|\cdot|$ on $k$, which implies that there are
elements of $k$ whose absolute value is comparable to any given
positive real number.

        Suppose that $\phi$ is a bounded linear mapping from $V_1$
into $V_2$, and let $\{v_j\}_{j = 1}^\infty$ be a sequence of elements
of $V_1$ that converges to $0$.  As in the preceding paragraph, there
is a sequence $\{t_j\}_{j = 1}^\infty$ of nonzero elements of $k$ that
converges to $0$ and satisfies (\ref{lim_{j to infty} t_j^{-1} v_j =
  0}).  This implies that
\begin{equation}
\label{E = {t_j^{-1} v_j : j in {bf Z}_+} cup {0}}
        E = \{t_j^{-1} \, v_j : j \in {\bf Z}_+\} \cup \{0\}
\end{equation}
is a compact subset of $V_1$, by standard arguments, and hence that
$E$ is bounded in $V_1$, as in Section \ref{bounded sets}.  It follows
that $\phi(E)$ is a bounded set in $V_2$, by hypothesis, so that
\begin{equation}
\label{phi(v_j) = t_j phi(t_j^{-1} v_j) to 0 as j to infty}
 \phi(v_j) = t_j \, \phi(t_j^{-1} \, v_j) \to 0 \quad\hbox{as } j \to \infty
\end{equation}
in $V_2$, because $\{t_j\}_{j = 1}^\infty$ converges to $0$ in $k$.
This criterion for convergence of sequences using bounded sets was
mentioned in Section \ref{bounded sets}.  This shows that $\phi$ is
sequentially continuous at $0$ under these conditions.  Thus $\phi$ is
continuous on $V_1$, since $V_1$ is supposed to have a local base for
its topology at $0$ with only finitely or countably many elements.

\section{The strong product topology}
\label{strong product topology}
\setcounter{equation}{0}

        Let $I$ be a nonempty set, and let $X_j$ be a topological
space for each $j \in I$.  As in Section \ref{cartesian products},
we let
\begin{equation}
\label{X = prod_{j in I} X_j, 2}
        X = \prod_{j \in I} X_j
\end{equation}
be the Cartesian product of the $X_j$'s, and we let $x_j \in X_j$ be
the $j$th coordinate of $x \in X$ for each $j \in I$.  A set $W
\subseteq X$ is said to be an open set with respect to the
\emph{strong product topology}\index{strong product topology} if for
each $x \in W$ and $j \in I$ there is an open set $U_j \subseteq X_j$
such that $x_j \in U_j$ and
\begin{equation}
\label{U = prod_{j in I} U_j, 3}
        U = \prod_{j \in I} U_j
\end{equation}
is contained in $W$.  It is well known and not difficult to check that
this defines a topology on $X$, which is the same as the product
topology on $X$ when $I$ has only finitely many elements.  If $I$ is
any nonempty set, then every open set in $X$ with respect to the
product topology is also an open set with respect to the strong
product topology.  Equivalently, if $U_j \subseteq X_j$ is an open set
for each $j \in I$, then (\ref{U = prod_{j in I} U_j, 3}) is an open
set in $X$ with respect to the strong product topology, and the
collection of these open sets forms a base for the strong product
topology on $X$.  If $X_j$ is equipped with the discrete topology for
each $j \in I$, then the strong product topology on $X$ is the same as
the discrete topology on $X$.

        Now let $A_j$ be a commutative topological group for each
$j \in I$, and let
\begin{equation}
\label{A = prod_{j in I} A_j, 3}
        A = \prod_{j \in I} A_j
\end{equation}
be their Cartesian product.  As in Section \ref{cartesian products},
$A$ is also a commutative group with respect to coordinatewise
addition.  One can check that $A$ is a topological group with respect
to the strong product topology as well.

        Let $k$ be a field with a $q$-absolute value function $|\cdot|$
for some $q > 0$, and let $V_j$ be a topological vector space over
$k$ for each $j \in I$.  As in Section \ref{cartesian products}, the
Cartesian product
\begin{equation}
\label{V = prod_{j in I} V_j, 3}
        V = \prod_{j \in I} V_j
\end{equation}
is a vector space over $k$ too, where the vector space operations are
defined coordinatewise.  In particular, $V_j$ is a commutative
topological group with respect to addition for each $j$, so that $V$
is also a commutative topological group with respect to addition, as
in the preceding paragraph.  It is easy to see that for each $t \in
k$, multiplication by $t$ is continuous on $V$ with respect to the
strong product topology, because of the analogous property of $V_j$
for each $j \in I$.  However, if $I$ has infinitely many elements,
then $V$ is not necessarily a topological vector space over $k$ with
respect to the strong product topology.  More precisely, if $v \in V$
has infinitely many nonzero coordinates, then $t \mapsto t \, v$ is
not necessarily continuous as a mapping from $k$ into $V$, with
respect to the strong product topology on $V$.  Of course, this is not
a problem when $|\cdot|$ is the trivial absolute value function on
$k$.

        If $|\cdot|$ is nontrivial on $k$, then for each $j \in I$,
there is a local base for the topology of $V_j$ at $0$ consisting of
balanced open sets in $V_j$, as in Section \ref{topological vector
  spaces}.  This leads to a local base for the strong product topology
on $V$ at $0$ consisting of balanced open sets in $V$, by taking
Cartesian products of balanced open subsets of the $V_j$'s.
Similarly, if $|\cdot|$ is the trivial absolute value function on $k$,
and if there a local base for the topology of $V_j$ at $0$ consisting
of balanced open sets for each $j \in I$, then there is a local base
for the strong product topology on $V$ at $0$ consisting of balanced
open sets, for the same reasons as before.

\section{Direct sums}
\label{direct sums}
\setcounter{equation}{0}

        Let $I$ be a nonempty set again, and let $A_j$ be a commutative
group for each $j \in I$.  As in Section \ref{cartesian products}, the
Cartesian product
\begin{equation}
\label{prod_{j in I} A_j}
        \prod_{j \in I} A_j
\end{equation}
is a commutative group with respect to coordinatewise addition, and is
known as the direct product of the $A_j$'s.  The corresponding
\emph{direct sum}\index{direct sums} is denoted
\begin{equation}
\label{sum_{j in I} A_j}
        \sum_{j \in I} A_j,
\end{equation}
and is the subgroup of the direct product consisting of elements $a$
of (\ref{prod_{j in I} A_j}) whose $j$th coordinate $a_j$ is equal to
$0$ for all but finitely many $j \in I$.  Of course, the direct sum of
the $A_j$'s is the same as the direct product when $I$ has only
finitely many elements.

        Suppose now that $A_j$ is a commutative topological group
for each $j \in I$.  Note that (\ref{sum_{j in I} A_j}) is dense in
(\ref{prod_{j in I} A_j}) with respect to the product topology on
(\ref{prod_{j in I} A_j}) associated to the given topology on $A_j$
for each $j$.  However, if $\{0\}$ is a closed set in $A_j$ for each
$j \in I$, then (\ref{sum_{j in I} A_j}) is a closed subset of
(\ref{prod_{j in I} A_j}) with respect to the strong product topology.
To see this, let $a$ be an element of (\ref{prod_{j in I} A_j}) which
is not an element of (\ref{sum_{j in I} A_j}), so that $a_j \ne 0$
for infinitely many $j \in I$.  Put
\begin{equation}
\label{U_j = A_j setminus {0}}
        U_j = A_j \setminus \{0\}
\end{equation}
for every $j \in I$ such that $a_j \ne 0$, and $U_j = A_j$ for every
$j \in I$ such that $a_j = 0$.  Because $\{0\}$ is a closed set in
$A_j$ for each $j \in I$, by hypothesis, $A_j \setminus \{0\}$ is an
open set in $A_j$ for every $j \in I$.  This implies that
\begin{equation}
\label{U = prod_{j in I} U_j, 4}
        U = \prod_{j \in I} U_j
\end{equation}
is an open set in (\ref{prod_{j in I} A_j}) with respect to the strong
product topology.  By construction, $a \in U$, and every element of
$U$ is not in the direct sum (\ref{sum_{j in I} A_j}).  It follows
that the complement of (\ref{sum_{j in I} A_j}) in (\ref{prod_{j in I}
  A_j}) is an open set with respect to the strong product topology, so
that (\ref{sum_{j in I} A_j}) is a closed set with respect to the
strong product topology.

        As in the previous section, the direct product
(\ref{prod_{j in I} A_j}) is a commutative topological group
with respect to the strong product topology, which implies that the
direct sum (\ref{sum_{j in I} A_j}) is a commutative topological group
with respect to the induced topology.  Let us continue to ask that
$\{0\}$ be a closed set in $A_j$ for each $j \in I$, and suppose that
$E$ is a totally bounded subset of (\ref{sum_{j in I} A_j}) with
respect to the topology induced on (\ref{sum_{j in I} A_j}) by the
strong product topology on (\ref{prod_{j in I} A_j}), as in Section
\ref{totally bounded sets, continued}.  Put
\begin{equation}
\label{I(E) = {j in I : there is an a in E such that a_j ne 0}}
        I(E) = \{j \in I : \hbox{ there is an } a \in E \hbox{ such that }
                                               a_j \ne 0\},
\end{equation}
and let us show that $I(E)$ has only finitely many elements under
these conditions.  Suppose for the sake of a contradiction that $I(E)$
has infinitely many elements, and for each $j \in I(E)$, let $a(j)$
be an element of $E$ such that
\begin{equation}
\label{a_j(j) ne 0}
        a_j(j) \ne 0.
\end{equation}
If we put
\begin{equation}
\label{U_j = A_j setminus {a_j(j)}}
        U_j = A_j \setminus \{a_j(j)\}
\end{equation}
for each $j \in I(E)$, and $U_j = A_j$ when $j \in I \setminus I(E)$,
then
\begin{equation}
\label{U = (sum_{j in I} A_j) cap (prod_{j in I} U_j)}
        U = \Big(\sum_{j \in I} A_j\Big) \cap \Big(\prod_{j \in I} U_j\Big)
\end{equation}
is an open subset of the direct sum (\ref{sum_{j in I} A_j}) with
respect to the topology induced by the strong product topology on the
direct product (\ref{prod_{j in I} A_j}).  This uses the hypothesis
that $\{0\}$ be a closed set in $A_j$ for each $j \in I$ to get that
(\ref{U_j = A_j setminus {a_j(j)}}) is an open set in $A_j$ for every
$j \in I(E)$.  By construction, $0 \in U_j$ for every $j \in I$,
which implies that $0 \in U$.

        If $E$ is totally bounded as a subset of the direct sum
(\ref{sum_{j in I} A_j}) with respect to the topology induced by the
strong product topology on the direct product (\ref{prod_{j in I} A_j}),
then there are finitely many elements $b(1), \ldots, b(n)$ of
(\ref{sum_{j in I} A_j}) such that
\begin{equation}
\label{E subseteq bigcup_{l = 1}^n (b(l) + U)}
        E \subseteq \bigcup_{l = 1}^n (b(l) + U).
\end{equation}
Of course, for each $l = 1, \ldots, n$, we have that $b_j(l) = 0$ for
all but finitely many $j \in I$, because $b(l)$ is an element of
(\ref{sum_{j in I} A_j}).  If $I(E)$ has infinitely many elements,
then there is a $j_0 \in I(E)$ such that
\begin{equation}
\label{b_{j_0}(l) = 0}
        b_{j_0}(l) = 0
\end{equation}
for each $l = 1, \ldots, n$.  Let $a(j_0) \in E$ be as in the previous
paragraph, so that $a_{j_0}(j_0) \not\in U_{j_0}$ by the definition
(\ref{U_j = A_j setminus {a_j(j)}}) of $U_{j_0}$.  This implies that
\begin{equation}
\label{a_{j_0}(j_0) not in b_{j_0}(l) + U_{j_0}}
        a_{j_0}(j_0) \not\in b_{j_0}(l) + U_{j_0}
\end{equation}
for each $l = 1, \ldots, n$, because of (\ref{b_{j_0}(l) = 0}).  It
follows that
\begin{equation}
\label{a(j_0) not in b(l) + U}
        a(j_0) \not\in b(l) + U
\end{equation}
for each $l = 1, \ldots, n$, by the definition (\ref{U = (sum_{j in I}
  A_j) cap (prod_{j in I} U_j)}) of $U$.  This contradicts (\ref{E
  subseteq bigcup_{l = 1}^n (b(l) + U)}), since $a(j_0) \in E$. Thus
$I(E)$ has only finitely many elements in this situation, as desired.

\section{Direct sums, continued}
\label{direct sums, continued}
\setcounter{equation}{0}

        Let $k$ be a field, let $I$ be a nonempty set, and let
$V_j$ be a vector space over $k$ for each $j \in I$.  As in
Section \ref{cartesian products}, the Cartesian product
\begin{equation}
\label{prod_{j in I} V_j}
        \prod_{j \in I} V_j
\end{equation}
is a vector space over $k$ with respect to coordinatewise addition and
scalar multiplication, and is known as the direct product of the
$V_j$'s.  The \emph{direct sum}\index{direct sums}
\begin{equation}
\label{sum_{j in I} V_j}
        \sum_{j \in I} V_j
\end{equation}
of the $V_j$'s is the linear subspace of the direct product
(\ref{prod_{j in I} V_j}) consistinng of the vectors $v$ in
(\ref{prod_{j in I} V_j}) whose $j$th coordinate $v_j$ is equal to $0$
for all but finitely many $j \in I$.  Of course, a vector space over
$k$ is a commutative group with respect to addition in particular, and
the commutative groups corresponding to the direct sum or product of
the $V_j$'s as vector spaces over $k$ are the same as the direct sum
or product of the $V_j$'s as commutative groups, respectively.  As
before, the direct sum of the $V_j$'s is the same as the direct
product when $I$ has only finitely many elements.

        Suppose from now on in this section that $|\cdot|$ is a
$q$-absolute value function on $k$ for some $q > 0$, and that
$V_j$ is a topological vector space with respect to $|\cdot|$
on $k$ for each $j \in I$.  Under these conditions, one can check
that the direct sum (\ref{sum_{j in I} V_j}) is a topological
vector space with respect to $|\cdot|$ on $k$ as well, using the
topology induced on the direct sum by the strong product topology
on the direct product (\ref{prod_{j in I} V_j}).  The main point
is that if $v$ is an element of (\ref{sum_{j in I} V_j}), then
$t \mapsto t \, v$ defines a continuous mapping from $k$ into
(\ref{sum_{j in I} V_j}), with respect to the topology induced
on (\ref{sum_{j in I} V_j}) by the strong product topology
on (\ref{prod_{j in I} V_j}).  As in Section \ref{strong product
topology}, the analogous statement for $v$ in the direct product
(\ref{prod_{j in I} V_j}) does not necessarily hold when $I$
has infinitely any elements.  However, some other properties of
the direct sum can be derived from analogous statements for the
direct product, as in Section \ref{strong product topology}.

        Suppose now in addition that $|\cdot|$ is not the trivial
absolute value function on $k$, and that $\{0\}$ is a closed set in
$V_j$ for each $j \in I$.  Let $E$ be a bounded subset of the direct
sum (\ref{sum_{j in I} V_j}), with respect to the topology on
(\ref{sum_{j in I} V_j}) induced by the strong product topology on the
direct product (\ref{prod_{j in I} V_j}), as in Section \ref{bounded
  sets}.  As in the previous section, we put
\begin{equation}
\label{I(E) = {j in I : there is a v in E such that v_j ne 0}}
 I(E) = \{j \in I : \hbox{ there is a } v \in E \hbox{ such that } v_j \ne 0\},
\end{equation}
and we would like to show that $I(E)$ has only finitely many elements
under these conditions.  Suppose for the sake of a contradiction that
$I(E)$ has infinitely many elements, and let $\{j_l\}_{l = 1}^\infty$
be an infinite sequence of distinct elements of $I(E)$.  Thus for
each $l \in {\bf Z}_+$ there is a vector $v(j_l) \in E$ such that
$v_{j_l}(j_l) \ne 0$.  Let $t_0$ be an element of $k$ such that
$|t_0| > 1$, which exists because $|\cdot|$ is nontrivial on $k$.
Observe that
\begin{equation}
\label{V_{j_l} setminus {t_0^{-l} v_{j_l}(j_l)}}
        V_{j_l} \setminus \{t_0^{-l} \, v_{j_l}(j_l)\}
\end{equation}
is an open set in $V_{j_l}$ that contains $0$ for each $l \in {\bf
  Z}_+$, because $\{0\}$ is a closed set in $V_j$ for each $j \in I$.
It follows that there is a nonempty balanced open set $U_{j_l}
\subseteq V_{j_l}$ contained in (\ref{V_{j_l} setminus {t_0^{-l}
    v_{j_l}(j_l)}}) for each $l \in {\bf Z}_+$, as in Section
\ref{topological vector spaces}, using the nontriviality of $|\cdot|$
on $k$ again.  Put $U_j = V_j$ for every $j \in I$ such that $j \ne
j_l$ for each $l \in {\bf Z}_+$, so that $U_j$ is an open subset of
$V_j$ that contains $0$ for every $j \in I$.  This implies that
\begin{equation}
\label{U = (sum_{j in I} V_j) cap (prod_{j in I} U_j)}
        U = \Big(\sum_{j \in I} V_j\Big) \cap \Big(\prod_{j \in I} U_j\Big)
\end{equation}
is an open subset of the direct sum (\ref{sum_{j in I} V_j}) with
respect to the topology induced by the strong product topology on the
direct product (\ref{prod_{j in I} V_j}), and that $0 \in U$.

        If $E$ is a bounded subset of the direct sum
(\ref{sum_{j in I} V_j}) with respect to this topology, then it follows
that $E \subseteq t \, U$ for every $t \in k$ such that $|t|$ is
sufficiently large.  In particular, if $|t|$ is sufficiently large,
then
\begin{equation}
\label{v_{j_l} in t U}
        v_{j_l} \in t \, U
\end{equation}
for every $l \in {\bf Z}_+$, because $v(j_l) \in E$.  Using the
definition (\ref{U = (sum_{j in I} V_j) cap (prod_{j in I} U_j)}) of
$U$, we get that if $|t|$ is sufficiently large, then
\begin{equation}
\label{v_{j_l}(j_l) in t U_{j_l}}
        v_{j_l}(j_l) \in t \, U_{j_l}
\end{equation}
for every $l \in {\bf Z}_+$.  However, we also have that
\begin{equation}
\label{v_{j_l}(j_l) not in t_0^l U_{j_l}}
        v_{j_l}(j_l) \not\in t_0^l \, U_{j_l}
\end{equation}
for every $l \in {\bf Z}_+$, because $U_{j_l}$ is contained in
(\ref{V_{j_l} setminus {t_0^{-l} v_{j_l}(j_l)}}) for every $l$, by
construction.  This implies that
\begin{equation}
\label{v_{j_l}(j_l) not in t U_{j_l}}
        v_{j_l}(j_l) \not\in t \, U_{j_l}
\end{equation}
when $|t| \le |t_0|^l$, because $U_{j_l}$ is supposed to be balanced
in $V_{j_l}$ for every $l$.  If $t \in k$ is given, then $|t| \le
|t_0|^l$ for all but finitely many $l \in {\bf Z}_+$, because $|t_0| >
1$, as in the preceding paragraph.  This leads to a contradiction, so
that $I(E)$ should have only finitely many elements, as desired.

\section{Combining semimetrics}
\label{combining semimetrics}
\setcounter{equation}{0}

        Let $X$ be a set, let $q$ be a positive real number, and
let $d_1, \ldots, d_l$ be finitely many $q$-semimetrics on $X$.  If $q
\le r < \infty$, then
\begin{equation}
\label{(sum_{j = 1}^l d_j(x, y)^r)^{1/r}}
        \Big(\sum_{j = 1}^l d_j(x, y)^r\Big)^{1/r}
\end{equation}
is a $q$-semimetric on $X$ as well.  To see this, the main point is to
verify that (\ref{(sum_{j = 1}^l d_j(x, y)^r)^{1/r}}) satisfies the
$q$-semimetric version of the triangle inequality.  This is easy to do
when $r = q$, and otherwise one can use Minkowski's inequality, which
is the triangle inequality for the $\ell^{r/q}$ norm when $r/q \ge 1$.
There is an analogous statement for infinite families of
$q$-semimetrics on $X$, as long as the infinite sum corresponding to
the finite sum in (\ref{(sum_{j = 1}^l d_j(x, y)^r)^{1/r}}) is finite
for every $x, y \in X$.

        The analogue of (\ref{(sum_{j = 1}^l d_j(x, y)^r)^{1/r}}) for
$r = \infty$ is
\begin{equation}
\label{max_{1 le j le l} d_j(x, y)}
        \max_{1 \le j \le l} d_j(x, y).
\end{equation}
It is easy to check directly that (\ref{max_{1 le j le l} d_j(x, y)})
satisfies the $q$-semimetric version of the triangle inequality, and
hence defines a $q$-semimetric on $X$.  This also works when $q =
\infty$, which is to say that if $d_1, \ldots, d_l$ are
semi-ultrametrics on $X$, then (\ref{max_{1 le j le l} d_j(x, y)}) is
a semi-ultrametric on $X$ too.  Similarly, the supremum of an infinite
family of $q$-semimetrics on $X$ is a $q$-semimetric on $X$, as long
as the supremum is finite for every $x, y \in X$.  As before, this
works when $q = \infty$, so that the supremum of an infinite family of
semi-ultrametrics on $X$ is a semi-ultrametric on $X$, as long as the
supremum is finite for every $x, y \in X$.

        Now let $I$ be a nonempty set, and let $X_j$ be a set for each
$j \in I$.  Also let $q$ be a positive real number, and let $\mathcal{M}_j$
be a nonempty collection of $q$-semimetrics on $X_j$ for each $j \in
I$.  As in the previous paragraphs, we ask here that these
$q$-semimetrics use the same $q$, which can always be arranged as in
Section \ref{q-semimetrics}.  Let us also ask that for each $j \in I$,
the collection of open balls in $X_j$ associated to elements of
$\mathcal{M}_j$ form a base for the topology on $X_j$ determined by
$\mathcal{M}_j$, and not just a sub-base, as in Section
\ref{collections of semimetrics}.  In particular, this holds when the
maximum of any finite number of elements of $\mathcal{M}_j$ is an
element of $\mathcal{M}_j$, which can easily be arranged by adding
these maxima to $\mathcal{M}_j$, if necessary.

        Let $d_j$ be an element of $\mathcal{M}_j$ for some $j \in I$,
and let $a_j$ be a positive real number, so that $a_j \, d_j(\cdot, \cdot)$
is a $q$-semimetric on $X_j$ too.  Note that
\begin{equation}
\label{B_{a_j d_j}(x_j, r) = B_{d_j}(x_j, r / a_j)}
        B_{a_j \, d_j}(x_j, r) = B_{d_j}(x_j, r / a_j)
\end{equation}
for every $x_j \in X_j$ and $r > 0$, where these open balls in $X_j$
centered at $x_j$ associated to $d_j$ and $a_j \, d_j$ are defined as
in (\ref{B(x, r) = B_d(x, r) = {y in X : d(x, y) < r}}), as usual.
Put
\begin{equation}
\label{d_j'(x_j, y_j) = min(a_j d_j(x_j, y_j), 1)}
        d_j'(x_j, y_j) = \min(a_j \, d_j(x_j, y_j), 1)
\end{equation}
for every $x_j, y_j \in X_j$, which defines a $q$-semimetric on $X_j$
as well, as in Section \ref{bounded semimetrics}.  If $0 < r \le 1$,
then
\begin{equation}
\label{B_{d_j'}(x_j, r) = B_{a_j d_j}(x_j, r) = B_{d_j}(x_j, r / a_j)}
        B_{d_j'}(x_j, r) = B_{a_j \, d_j}(x_j, r) = B_{d_j}(x_j, r / a_j)
\end{equation}
for every $x_j \in X_j$, where $B_{d_j'}(x_j, r)$ is the open ball in
$X_j$ centered at $x_j$ with radius $r$ associated to $d_j'$.  The
first equality in (\ref{B_{d_j'}(x_j, r) = B_{a_j d_j}(x_j, r) =
  B_{d_j}(x_j, r / a_j)}) corresponds to (\ref{B_{d'}(x, r) = B_d(x,
  r)}), and we have that
\begin{equation}
\label{B_{d_j'}(x_j, r) = X_j}
        B_{d_j'}(x_j, r) = X_j
\end{equation}
for every $x_j \in X_j$ when $r > 1$, as in (\ref{B_{d'}(x, r) = X}).

        As in Section \ref{cartesian products}, let
\begin{equation}
\label{X = prod_{j in I} X_j, 3}
        X = \prod_{j \in I} X_j
\end{equation}
be the Cartesian product of the $X_j$'s, and let $x_j \in X_j$ be the
$j$th component of $x \in X$ for each $j \in I$.  Also let $X$ be
equipped with the strong product topology corresponding to the
topology on $X_j$ determined by $\mathcal{M}_j$ for each $j \in I$.
Suppose that $d_j \in \mathcal{M}_j$ and $a_j > 0$ are given as in the
preceding paragraph for each $j \in I$, and put
\begin{equation}
\label{d(x, y) = sup_{j in I} d_j'(x_j, y_j)}
        d(x, y) = \sup_{j \in I} d_j'(x_j, y_j)
\end{equation}
for every $x, y \in X$, where $d_j'$ is associated to $d_j$ and $a_j$
as in (\ref{d_j'(x_j, y_j) = min(a_j d_j(x_j, y_j), 1)}).  Of course,
\begin{equation}
        d(x, y) \le 1
\end{equation}
for every $x, y \in X$, because $d_j'(x_j, y_j) \le 1$ for each $j \in
I$, by construction.  Thus (\ref{d(x, y) = sup_{j in I} d_j'(x_j,
  y_j)}) is finite, and hence defines a $q$-semimetric on $X$ under
these conditions, as mentioned earlier.

        Observe that
\begin{equation}
\label{B_d(x, r) subseteq prod_{j in I} B_{d_j'}(x_j, r)}
        B_d(x, r) \subseteq \prod_{j \in I} B_{d_j'}(x_j, r)
\end{equation}
for every $x \in X$ and $r > 0$, where $B_d(x, r)$ is the open ball in
$X$ centered at $x$ with radius $r$ associated to $d$ as in (\ref{B(x,
  r) = B_d(x, r) = {y in X : d(x, y) < r}}), and $B_{d_j'}(x_j, r)$ is
the open ball in $X_j$ centered at $x_j$ with radius $r$ associated to
$d_j'$ for each $j \in I$, as before.  More precisely, one can check
that
\begin{equation}
\label{B_d(x, r) = bigcup_{0 < widetilde{r} < r} ...}
        B_d(x, r) = \bigcup_{0 < \widetilde{r} < r}
                     \, \prod_{j \in I} B_{d_j'}(x_j, \widetilde{r})
\end{equation}
for every $x \in X$ and $r > 0$.  This implies that $B_d(x, r)$ is an
open set in $X$ with respect to the strong product topology for every
$x \in X$ and $r > 0$.  This uses the fact that $B_{d_j'}(x_j,
\widetilde{r})$ is an open set in $X_j$ for each $j \in I$, by
(\ref{B_{d_j'}(x_j, r) = B_{a_j d_j}(x_j, r) = B_{d_j}(x_j, r / a_j)})
and (\ref{B_{d_j'}(x_j, r) = X_j}), so that
\begin{equation}
\label{prod_{j in I} B_{d_j'}(x_j, widetilde{r})}
        \prod_{j \in I} B_{d_j'}(x_j, \widetilde{r})
\end{equation}
is an open set in $X$ with respect to the strong product topology.
One can also verify that the topology on $X$ determined by the
collection of $q$-semimetrics on $X$ of the form (\ref{d(x, y) =
  sup_{j in I} d_j'(x_j, y_j)}) with $d_j \in \mathcal{M}_j$ and $a_j
> 0$ is the same as the strong product topology, using (\ref{B_d(x, r)
  subseteq prod_{j in I} B_{d_j'}(x_j, r)}).  It suffices for this to
use any collection of positive real numbers $a_j$ that can be
arbitrarily large, instead of all $a_j > 0$.  It is here that we use
the hypothesis that the open balls in $X_j$ associated to elements of
$\mathcal{M}_j$ form a base for the topology of $X_j$ for each $j \in
I$.

        Suppose that $X_j$ is a commutative group for each $j \in I$,
so that $X$ is also a commutative group with respect to coordinatewise
addition.  If $d_j$ is a translation-invariant $q$-semimetric on $X_j$
for some $j \in I$, then (\ref{d_j'(x_j, y_j) = min(a_j d_j(x_j, y_j),
  1)}) is invariant under translations on $X_j$ for every $a_j > 0$.
If $d_j$ is a translation-invariant $q$-semimetric on $X_j$ and $a_j >
0$ for each $j \in I$, then it follows that (\ref{d(x, y) = sup_{j in
    I} d_j'(x_j, y_j)}) is invariant under translations on $X$.  If
every element of $\mathcal{M}_j$ is invariant under translations on
$X_j$ for each $j \in I$, then we get a collection of
translation-invariant $q$-semimetrics on $X$ for which the
corresponding topology is the same as the strong product topology, as
in the previous paragraph.

\section{Combining seminorms}
\label{combining seminorms}
\setcounter{equation}{0}

        Let $k$ be a field with a $q$-absolute value function $|\cdot|$
for some $q > 0$, and let $V$ be a vector space over $k$.  If $N_1,
\ldots, N_l$ are finitely many $q$-seminorms on $V$, and $q \le r <
\infty$, then
\begin{equation}
\label{(sum_{j = 1}^l N_j(v)^r)^{1/r}}
        \Big(\sum_{j = 1}^l N_j(v)^r\Big)^{1/r}
\end{equation}
is a $q$-seminorm on $V$ too.  As in the previous section, the main
point is to check that (\ref{(sum_{j = 1}^l N_j(v)^r)^{1/r}})
satisfies the $q$-seminorm version of the triangle inequality, which
is essentially the same as for the corresponding $q$-semimetrics.  As
before, the analogue of (\ref{(sum_{j = 1}^l N_j(v)^r)^{1/r}}) for $r
= \infty$ is
\begin{equation}
\label{max_{1 le j le l} N_j(v)}
        \max_{1 \le j \le l} N_j(v).
\end{equation}
Similarly, if the $N_j$'s are semi-ultranorms on $V$, then
(\ref{max_{1 le j le l} N_j(v)}) is a semi-ultranorm on $V$ as well.

        As in the previous section again, there are analogous statements
for infinite families of $q$-seminorms, as long as the relevant
quantities are finite.  In particular, this can be applied to direct
sums of vector spaces.  Let $I$ be a nonempty set, let $V_j$ be a
vector space over $k$ for each $j \in I$, and let
\begin{equation}
\label{V = sum_{j in I} V_j}
        V = \sum_{j \in I} V_j
\end{equation}
be the corresponding direct sum, as in Section \ref{direct sums,
  continued}.  If $N_j$ is a $q$-seminorm on $V_j$ for each $j \in I$,
and if $q \le r < \infty$, then
\begin{equation}
\label{(sum_{j in I} N_j(v_j)^r)^{1/r}}
        \Big(\sum_{j \in I} N_j(v_j)^r\Big)^{1/r}
\end{equation}
defines a $q$-seminorm on $V$.  Here $v_j \in V_j$ is the $j$th
component of $v \in V$ for each $j \in I$, as usual, which is equal to
$0$ for all but finitely many $j \in I$, by the definition of the
direct sum.  This implies that $N_j(v_j) = 0$ for all but finitely
many $j \in I$, so that the sum in (\ref{(sum_{j in I}
  N_j(v_j)^r)^{1/r}}) is finite.  Similarly, the maximum
\begin{equation}
\label{max_{j in I} N_j(v_j)}
        \max_{j \in I} N_j(v_j)
\end{equation}
is attained for every $v \in V$, and defines a $q$-seminorm on $V$
too.  If $N_j$ is a semi-ultranorm on $V_j$ for each $j \in I$, then
(\ref{max_{j in I} N_j(v_j)}) defines a semi-ultranorm on $V$.

        Suppose now that $\mathcal{N}_j$ is a nonempty collection
of $q$-seminorms on $V_j$ for each $j \in I$, where more precisely the
same $q$ should be used for every element of $\mathcal{N}_j$ for each
$j \in I$.  Let $V_j$ be equipped with the topology determined by the
collection of $q$-semimetrics associated to the elements of
$\mathcal{N}_j$ for each $j \in I$, and let $V$ be equipped with the
topology induced by the corresponding strong product topology on the
Cartesian product of the $V_j$'s.  As in the preceding section, we ask
that the open balls in $V_j$ centered at $0$ associated to elements of
$\mathcal{N}_j$ form a local base for the topology of $V_j$ at $0$ for
each $j \in I$, and not just a sub-base.  As before, this holds
automatically when the maximum of any finite number of elements of
$\mathcal{N}_j$ is an element of $\mathcal{N}_j$ for each $j \in I$,
which can always be arranged by adding these maxima to
$\mathcal{N}_j$.

        Let $N_j$ be an element of $\mathcal{N}_j$ for each $j \in I$,
and let $a_j$ be a positive real number for each $j \in I$.  Thus
$a_j \, N_j$ is a $q$-seminorm on $V_j$ for each $j \in I$, so that
\begin{equation}
\label{N(v) = max_{j in I} (a_j N_j(v_j))}
        N(v) = \max_{j \in I} (a_j \, N_j(v_j))
\end{equation}
defines a $q$-seminorm on $V$, as in (\ref{max_{j in I} N_j(v_j)}).
As in (\ref{B_{a_j d_j}(x_j, r) = B_{d_j}(x_j, r / a_j)}),
\begin{equation}
\label{B_{a_j N_j}(0, r) = B_{N_j}(0, r / a_j)}
        B_{a_j \, N_j}(0, r) = B_{N_j}(0, r / a_j)
\end{equation}
for every $j \in I$ and $r > 0$, where these open balls in $V_j$
associated to $N_j$ and $a_j \, N_j$ are defined as in (\ref{B_N(0, r)
  = {v in V : N(v) < r}}).  It is easy to see that
\begin{equation}
\label{B_N(0, r) = V cap (prod_{j in I} B_{a_j N_j}(0, r))}
        B_N(0, r) = V \cap \Big(\prod_{j \in I} B_{a_j \, N_j}(0, r)\Big)
\end{equation}
for every $r > 0$, where the left side of (\ref{B_N(0, r) = V cap
  (prod_{j in I} B_{a_j N_j}(0, r))}) is the open ball in $V$
associated to $N$ as in (\ref{B_N(0, r) = {v in V : N(v) < r}}).
Remember that the direct sum $V$ is contained in the Cartesian product
of the $V_j$'s, which also contains the product of open balls in the
$V_j$'s on the right side of (\ref{B_N(0, r) = V cap (prod_{j in I}
  B_{a_j N_j}(0, r))}).  This step is a bit simpler than its analogue
in the previous section, because $N$ is defined in (\ref{N(v) = max_{j
    in I} (a_j N_j(v_j))}) as a maximum, instead of a supremum.  It
follows from (\ref{B_{a_j N_j}(0, r) = B_{N_j}(0, r / a_j)}) and
(\ref{B_N(0, r) = V cap (prod_{j in I} B_{a_j N_j}(0, r))}) that
\begin{equation}
\label{B_N(0, r) = V cap (prod_{j in I} B_{N_j}(0, r / a_j))}
        B_N(0, r) = V \cap \Big(\prod_{j \in I} B_{N_j}(0, r / a_j)\Big)
\end{equation}
for every $r > 0$.  Of course,
\begin{equation}
\label{prod_{j in I} B_{N_j}(0, r / a_j)}
        \prod_{j \in I} B_{N_j}(0, r / a_j)
\end{equation}
is an open set in the Cartesian product of the $V_j$'s with respect to
the strong product topology for every $r > 0$, since open balls in
$V_j$ with respect to $N_j$ are open subsets of $V_j$ for each $j \in
I$ by construction.  This implies that (\ref{B_N(0, r) = V cap
  (prod_{j in I} B_{N_j}(0, r / a_j))}) is an open set in $V$ with
respect to the topology induced by the strong product topology on the
Cartesian product of the $V_j$'s.  Using (\ref{B_N(0, r) = V cap
  (prod_{j in I} B_{N_j}(0, r / a_j))}), one can also check that the
topology induced on $V$ by the strong product topology on the
Cartesian product of the $V_j$'s is the same as the topology
determined by the collection of $q$-seminorms of the form (\ref{N(v) =
  max_{j in I} (a_j N_j(v_j))}), where $N_j \in \mathcal{N}_j$ and
$a_j > 0$ for each $j \in I$.  As in the previous section, it suffices
to use any collection of positive real numbers $a_j$ that can be
arbitrarily large here, instead of all $a_j > 0$.  It is also here
that we use the hypothesis that open balls in $V_j$ centered at $0$
associated to elements of $\mathcal{N}_j$ form a local base for the
topology of $V_j$ at $0$ for each $j \in I$, as before.

\section{Continuous linear mappings, continued}
\label{continuous linear mappings, continued}
\setcounter{equation}{0}

        Let $k$ be a field with a $q$-absolute value function $|\cdot|$
for some $q > 0$, and let $V$, $W$ be vector spaces over $k$.  Also
let $\mathcal{N}_V$, $\mathcal{N}_W$ be nonempty collections of
$q$-seminorms on $V$, $W$, respectively.  As before, one can let
$q > 0$ depend on the elements of $\mathcal{N}_V$ and $\mathcal{N}_W$,
as long as $|\cdot|$ is a $q$-absolute value function on $k$ for that
choice of $q$.  Note that for any finite collection of elements of
$\mathcal{N}_V$ or $\mathcal{N}_W$, there is a $q > 0$ that works for
each element of the finite collection, by taking the minimum of the
corresponding finitely many $q$'s.  Let $V$, $W$ be equipped with the
topologies determined by $\mathcal{N}_V$, $\mathcal{N}_W$, respectively,
so that $V$ and $W$ are topological vector spaces.

        Let $\phi$ be a linear mapping from $V$ into $W$.  It is
easy to see that $\phi$ is continuous at $0$ if and only if
for every $N_W \in \mathcal{N}_W$ and $r > 0$ there are finitely
many elements $N_{V, 1}, \ldots, N_{V, l}$ of $\mathcal{N}_V$ and
finitely many positive real numbers $r_1, \ldots, r_l$ such that
\begin{equation}
\label{phi(bigcap_{j = 1}^l B_{N_{V, j}}(0, r_j)) subseteq B_{N_W}(0, r)}
 \phi\Big(\bigcap_{j = 1}^l B_{N_{V, j}}(0, r_j)\Big) \subseteq B_{N_W}(0, r).
\end{equation}
Here $B_{N_W}(0, r)$ is the open ball in $W$ associated to $N_W$ and
$r$ as in (\ref{B_N(0, r) = {v in V : N(v) < r}}), and similarly
$B_{N_{V, j}}(0, r_j)$ is the open ball in $V$ associated to $N_{V,
  j}$ and $r_j$ for $j = 1, \ldots, l$.  In this case, we get that
\begin{equation}
\label{phi(bigcap B_{N_{V, j}}(0, |t| r_j)) subseteq B_{N_W}(0, |t| r)}
        \phi\Big(\bigcap_{j = 1}^l B_{N_{V, j}}(0, |t| \, r_j)\Big)
                                     \subseteq B_{N_W}(0, |t| \, r)
\end{equation}
for every $t \in k$ with $t \ne 0$, using (\ref{t B_N(0, r) = B_N(0,
  |t| r)}).  If $|\cdot|$ is nontrivial on $k$, then it follows that a
condition like (\ref{phi(bigcap_{j = 1}^l B_{N_{V, j}}(0, r_j))
  subseteq B_{N_W}(0, r)}) for a single $r > 0$ implies an analogous
condition for every $r > 0$.

        Suppose that there is a positive real number $C$ such that
\begin{equation}
\label{N_W(phi(v)) le C max_{1 le j le l} N_{V, j}(v)}
        N_W(\phi(v)) \le C \, \max_{1 \le j \le l} N_{V, j}(v)
\end{equation}
for every $v \in V$.  This implies that for each $r > 0$,
(\ref{phi(bigcap_{j = 1}^l B_{N_{V, j}}(0, r_j)) subseteq B_{N_W}(0,
  r)}) holds with $r_j = r / C$ for $j = 1, \ldots, l$.  As a partial
converse, if (\ref{phi(bigcap_{j = 1}^l B_{N_{V, j}}(0, r_j)) subseteq
  B_{N_W}(0, r)}) holds for some positive real numbers $r, r_1,
\ldots, r_l$, and if $|\cdot|$ is not trivial on $k$, then one can
check that there is a $C > 0$ such that (\ref{N_W(phi(v)) le C max_{1
    le j le l} N_{V, j}(v)}) holds for every $v \in V$.  This is a bit
simpler when $|\cdot|$ is not discrete on $k$, so that (\ref{{|x| : x
    in k, x ne 0}}) is dense in ${\bf R}_+$ with respect to the
standard Euclidean topology on ${\bf R}$.  Otherwise, if $|\cdot|$ is
nontrivial and discrete on $k$, then the constant $C$ just mentioned
also depends on the quantity (\ref{rho_1 = sup {|x| : x in k, |x| <
    1}}) associated to $|\cdot|$ on $k$.

        Suppose now that $|\cdot|$ is the trivial absolute value
function on $k$, and that $\mathcal{N}_V$ consists of only the trivial
ultranorm on $V$.  This implies that the topology on $V$ determined by
$\mathcal{N}_V$ is the same as the discrete topology, so that any
mapping from $V$ into any topological space is continuous.  In this
situation, (\ref{N_W(phi(v)) le C max_{1 le j le l} N_{V, j}(v)})
would say that
\begin{equation}
\label{N_W(phi(v)) le C}
        N_W(\phi(v)) \le C
\end{equation}
for every $v \in V$ with $v \ne 0$.  If $V$ is finite-dimensional over
$k$, then there is always a $C > 0$ so that this holds.  Otherwise, if
$V$ is infinite-dimensional over $k$, then one can give examples where
this does not work, with $W = V$ and $\phi$ equal to the identity
mapping.

\section{Direct sums, revisited}
\label{direct sums, revisited}
\setcounter{equation}{0}

        Let $I$ be a nonempty set, let $A_j$ be a commutative group
for each $j \in I$, and let
\begin{equation}
\label{A = sum_{j in I} A_j}
        A = \sum_{j \in I} A_j
\end{equation}
be the corresponding direct sum, as in Section \ref{direct sums}.
Note that for each $l \in I$, there is a natural embedding of $A_l$
into $A$, which sends each element of $A_l$ to the element of $A$
whose $l$th coordinate is that element of $A_l$, and whose $j$th
coordinate is equal to $0$ for every $j \in I$ with $j \ne l$.  If
$\phi$ is a homomorphism from $A$ into another commutative group $B$,
then one gets a homomorphism $\phi_l$ from $A_l$ into $B$ for each $l
\in I$, by composing the natural embedding of $A_l$ into $A$ with
$\phi$.  In the other direction, if $\phi_j$ is a homomorphism from
$A_j$ into $B$ for each $j \in I$, then
\begin{equation}
\label{phi(a) = sum_{j in I} phi_j(a_j)}
        \phi(a) = \sum_{j \in I} \phi_j(a_j)
\end{equation}
defines a homomorphism from $A$ into $B$.  Remember that for each $a
\in A$, $a_j = 0$ for all but finitely $j \in I$, which implies that
$\phi_j(a_j) = 0$ for all but finitely many $j \in I$, so that the sum
in (\ref{phi(a) = sum_{j in I} phi_j(a_j)}) is makes sense.  If the
$A_j$'s are vector spaces over a field $k$, then $A$ is also a vector
space over $k$, and the natural embedding of $A_l$ into $A$ is linear
for each $l \in I$.  If $B$ is another vector space over $k$, then
linear mappings from $A$ into $B$ correspond to families of linear
mappings from the $A_j$'s into $B$ as before.

        Suppose now that $A_j$ is a commutative topological group
for each $j \in I$, and let $A$ be equipped with the topology
induced by the corresponding strong product topology on the Cartesian
product of the $A_j$'s, as in Section \ref{direct sums}.  It is easy
to see that for each $l \in I$, the natural embedding from $A_l$
into $A$ is a homeomorphism onto its image with respect to the induced
topology.  Let $B$ be another commutative topological group, and
suppose that $\phi$ is a continuous homomorphism from $A$ into $B$.
If $\phi_j$ is the corresponding homomorphism from $A_j$ into $B$
for each $j \in I$, as in the preceding paragraph, then $\phi_j$ is
also continuous for each $j$, since it is the composition of two
continuous mappings.

        Now let $\phi_j$ be a continuous homomorphism from $A_j$
into $B$ for each $j \in I$, and let $\phi$ be the corresponding
homomorphism from $A$ into $B$, as in (\ref{phi(a) = sum_{j in I}
  phi_j(a_j)}).  If $I$ has only finitely many elements, then one
can check that $\phi$ is also continuous.  More precisely, put
\begin{equation}
\label{widetilde{phi}_j(a) = phi_j(a_j)}
        \widetilde{\phi}_j(a) = \phi_j(a_j)
\end{equation}
for every $a \in A$ and $j \in I$, which is the same as the
composition of $\phi_j$ with the natural coordinate projection from
$A$ on $A_j$.  Thus $\widetilde{\phi}_j$ is continuous as a mapping
from $A$ into $B$ for each $j \in I$, because it is the composition of
two continuous mappings.  If $I$ has only finitely many elements, then
it follows that $\phi$ is continuous as a mapping from $A$ into $B$ as
well, since it is the sum of finitely many continuous mappings from
$A$ into $B$.

        In order to prove an analogous statement when $I$ is
countably infinite, let us begin with some consequences of continuity
of addition on $B$.  Let $W$ be an open subset of $B$ that contains
$0$, and let $W_1 \subseteq B$ be an open set such that $0 \in W_1$
and
\begin{equation}
\label{W_1 + W_1 subseteq W}
        W_1 + W_1 \subseteq W.
\end{equation}
Continuing in this way, we get for each integer $n \ge 2$ an open set
$W_n \subseteq B$ such that $0 \in W_n$ and
\begin{equation}
\label{W_n + W_n  subseteq W_{n - 1}}
        W_n + W_n \subseteq W_{n - 1}.
\end{equation}
This implies that
\begin{equation}
\label{W_1 + W_2 + cdots + W_{n - 1} + W_n + W_n subseteq ...}
        \qquad W_1 + W_2 + \cdots + W_{n - 1} + W_n + W_n
         \subseteq W_1 + W_2 + \cdots + W_{n - 1} + W_{n - 1}
\end{equation}
when $n \ge 2$, and hence that
\begin{equation}
\label{W_1 + cdots + W_n + W_n subseteq W}
        W_1 + \cdots + W_n + W_n \subseteq W
\end{equation}
for every $n \ge 1$.  In particular,
\begin{equation}
\label{W_1 + cdots + W_n subseteq W}
        W_1 + \cdots + W_n \subseteq W
\end{equation}
for each $n$, since $0 \in W_n$.

        To deal with the case where $I$ is countably infinite, we may
as well suppose that $I = {\bf Z}_+$.  Let $W$ be given as in the
previous paragraph, and let $W_n$ be chosen for $n \in {\bf Z}_+$
as before.  Suppose that $\phi_j$ is a continuous homomorphism from
$A_j$ into $B$ for each $j \in I = {\bf Z}_+$, and let $U_j$ be an
open subset of $A_j$ such that $0 \in U_j$ and
\begin{equation}
\label{phi_j(U_j) subseteq W_j}
        \phi_j(U_j) \subseteq W_j
\end{equation}
for each $j$.  Put
\begin{equation}
\label{U = A cap (prod_{j = 1}^infty U_j)}
        U = A \cap \Big(\prod_{j = 1}^\infty U_j\Big),
\end{equation}
which is an open subset of $A$ with respect to the topology induced by
the strong product topology on $\prod_{j = 1}^\infty A_j$.  Note that
$0 \in U$, since $0 \in U_j$ for each $j$.  One can also check that
\begin{equation}
\label{phi(U) subseteq W, 2}
        \phi(U) \subseteq W,
\end{equation}
where $\phi$ is as in (\ref{widetilde{phi}_j(a) = phi_j(a_j)}), using
(\ref{W_1 + cdots + W_n subseteq W}) and (\ref{phi_j(U_j) subseteq
  W_j}).  This implies that $\phi$ is continuous at $0$ with respect
to the topology on $A$ induced by the strong product topology on
$\prod_{j = 1}^\infty A_j$, and hence that $\phi$ is continuous
everywhere on $A$ with respect to this topology, because $\phi$ is a
homomorphism.

        Here is another version of this type of argument.  Suppose
that $d_B(\cdot, \cdot)$ is a $q$-semimetric on $B$ for some $q > 0$
that is invariant under translations, and compatible with the given
topology on $B$.  Let $\epsilon > 0$ be given.  If $I$ has only
finitely or countably many elements, then one can choose $\epsilon_j
> 0$ for every $j \in I$ so that
\begin{equation}
\label{sum_{j in I} epsilon_j^q le epsilon^q}
        \sum_{j \in I} \epsilon_j^q \le \epsilon^q,
\end{equation}
which is the same as saying that the sum of $\epsilon_j^q$ over any
finite subset of $I$ is less than or equal to $\epsilon^q$.  Suppose
that $\phi_j$ is a continuous homomorphism from $A_j$ into $B$ for
each $j \in I$, so that $\phi_j$ is continuous at $0$ with respect to
$d_B(\cdot, \cdot)$ on $B$ in particular.  This implies that for each
$j \in I$, there is an open set $U_j \subseteq A_j$ such that $0 \in
A_j$ and
\begin{equation}
\label{phi_j(U_j) subseteq B_{d_B}(0, epsilon_j)}
        \phi_j(U_j) \subseteq B_{d_B}(0, \epsilon_j),
\end{equation}
where the right side of (\ref{phi_j(U_j) subseteq B_{d_B}(0,
  epsilon_j)}) is the open ball in $B$ with respect to $d_B(\cdot,
\cdot)$ centered at $0$ with radius $\epsilon_j$, as in (\ref{B(x, r)
  = B_d(x, r) = {y in X : d(x, y) < r}}).  As before,
\begin{equation}
\label{U = A cap (prod_{j in I} U_j)}
        U = A \cap \Big(\prod_{j \in I} U_j\Big)
\end{equation}
is an open subset of $A$ with respect to the topology induced by the
strong product topology on $\prod_{j \in I} A_j$, and $0 \in U$.
Using (\ref{sum_{j in I} epsilon_j^q le epsilon^q}), (\ref{phi_j(U_j)
  subseteq B_{d_B}(0, epsilon_j)}), and the $q$-semimetric version of
the triangle inequality, one can verify that
\begin{equation}
\label{phi(U) subseteq B_{d_B}(0, epsilon)}
        \phi(U) \subseteq B_{d_B}(0, \epsilon),
\end{equation}
where the right side of (\ref{phi(U) subseteq B_{d_B}(0, epsilon)}) is
the open ball in $B$ with respect to $d_B(\cdot, \cdot)$ centered at
$0$ with radius $\epsilon$.  This means that $\phi$ is continuous at
$0$ with respect to $d_B(\cdot, \cdot)$ on $B$ and the topology on $A$
induced by the strong product topology on $\prod_{j \in I} A_j$.

        Suppose now that $d_B(\cdot, \cdot)$ is a semi-ultrametric
on $B$ that is invariant under translations and compatible with the
given topology on $B$, which corresponds to $q = \infty$ in the
preceding paragraph.  In this case, the analogous argument works for
any nonempty set $I$, with $\epsilon_j = \epsilon$ for every $j \in I$.

\section{Combining seminorms, continued}
\label{combining seminorms, continued}
\setcounter{equation}{0}

        Let $k$ be a field with a $q$-absolute value function $|\cdot|$
for some $q > 0$, and let $I$ be a nonempty set.  Also let $V_j$
be a vector space over $k$ for each $j \in I$, and let
\begin{equation}
\label{V = sum_{j in I} V_j, 2}
        V = \sum_{j \in I} V_j
\end{equation}
be the corresponding direct sum, as in Section \ref{direct sums,
  continued}.  If $N_j$ is a $q$-seminorm on $V_j$ for each $j \in
I$, then
\begin{equation}
\label{(sum_{j in I} N_j(v_j)^r)^{1/r}, 2}
        \Big(\sum_{j \in I} N_j(v_j)^r\Big)^{1/r}
\end{equation}
defines a $q$-seminorm on $V$ when $q \le r < \infty$, as in Section
\ref{combining seminorms}.  Similarly,
\begin{equation}
\label{max_{j in I} N_j(v_j), 2}
        \max_{j \in I} N_j(v_j)
\end{equation}
defines a $q$-seminorm on $V$ too, as before.  If $I$ has only
finitely many elements, then (\ref{(sum_{j in I} N_j(v_j)^r)^{1/r},
  2}) is bounded by (\ref{max_{j in I} N_j(v_j), 2}) times the
$r$th root of the number of elements of $I$.

        Suppose now that  $I$ is countably infinite, and let $a_j$ be 
a positive real number for each $j \in I$ such that
\begin{equation}
\label{sum_{j in I} a_j^{-r} < infty}
        \sum_{j \in I} a_j^{-r} < \infty,
\end{equation}
where the sum is defined as the supremum of the corresponding finite
subsums.  Observe that
\begin{eqnarray}
        \Big(\sum_{j \in I} N_j(v_j)^r\Big)^{1/r}
  & = & \Big(\sum_{j \in I} a_j^{-r} \, (a_j \, N_j(v_j))^r\Big)^{1/r} \\
 & \le & \Big(\sum_{j \in I} a_j^{-r}\Big)^{1/r} \, 
                         \max_{j \in I} (a_j \, N_j(v_j)) \nonumber
\end{eqnarray}
for every $v \in V$.  Of course,
\begin{equation}
\label{max_{j in I} (a_j N_j(v_j))}
        \max_{j \in I} (a_j \, N_j(v_j))
\end{equation}
also defines a $q$-seminorm on $V$, as in (\ref{max_{j in I} N_j(v_j),
  2}).

        Let $I$ be any nonempty set again, and for each $l \in I$,
let $\eta_l$ be the obvious inclusion mapping from $V_l$ into $V$.
Thus for each $l \in I$ and $v_l \in V_l$, $\eta_l(v_l)$ is the
element of $V$ whose $l$th coordinate is equal to $v_l$, and whose
$j$th coordinate is equal to $0$ for every $j \in I$ with $j \ne l$.
Let $N$ be a $q$-seminorm on $V$, and suppose that for each $l \in I$
there is a nonnegative real number $C_l$ such that
\begin{equation}
\label{N(eta_l(v_l)) le C_l N_l(v_l)}
        N(\eta_l(v_l)) \le C_l \, N_l(v_l)
\end{equation}
for every $v_l \in V_l$.  Using the $q$-semimetric version of the
triangle inequality, we get that
\begin{equation}
\label{N(v) le (sum_{j in I} C_l^q N_l(v_l)^q)^{1/q}}
        N(v) \le \Big(\sum_{j \in I} C_l^q \, N_l(v_l)^q\Big)^{1/q}
\end{equation}
for every $v \in V$.  Note that the right side of (\ref{N(v) le
  (sum_{j in I} C_l^q N_l(v_l)^q)^{1/q}}) is a $q$-semimetric on $V$,
as in (\ref{(sum_{j in I} N_j(v_j)^r)^{1/r}, 2}).

        In the analogous situation for $q = \infty$, we ask that $N_l$
be a semi-ultranorm on $V_l$ for each $l \in I$, and that $N$ be a
semi-ultranorm on $V$.  In this case, (\ref{N(eta_l(v_l)) le C_l
  N_l(v_l)}) implies that
\begin{equation}
\label{N(v) le max_{j in I} (C_l N_l(v_l))}
        N(v) \le \max_{j \in I} (C_l \, N_l(v_l))
\end{equation}
for every $v \in V$.  As before, the right side of (\ref{N(v) le
  max_{j in I} (C_l N_l(v_l))}) is a semi-ultranorm on $V$ too.

\section{Completeness}
\label{completeness}
\setcounter{equation}{0}

        Let $X$ be a set, and let $d(x, y)$ be a $q$-semimetric on $X$
for some $q > 0$.  As usual, a sequence $\{x_j\}_{j = 1}^\infty$ of
elements of $X$ is said to be a \emph{Cauchy sequence}\index{Cauchy
sequences} with respect to $d(\cdot, \cdot)$ if
\begin{equation}
\label{lim_{j, l to infty} d(x_j, x_l) = 0}
        \lim_{j, l \to \infty} d(x_j, x_l) = 0.
\end{equation}
Of course, this happens if and only if
\begin{equation}
\label{lim_{j, l to infty} d(x_j, x_l)^q = 0}
        \lim_{j, l \to \infty} d(x_j, x_l)^q = 0,
\end{equation}
so that $\{x_j\}_{j = 1}^\infty$ is a Cauchy sequence with respect to
$d(\cdot, \cdot)^q$ as an ordinary semimetric on $X$.  If $\{x_j\}_{j
  = 1}^\infty$ converges to an element of $X$ with respect to the
topology determined by $d(\cdot, \cdot)$, then it is easy to see that
$\{x_j\}_{j = 1}^\infty$ is a Cauchy sequence with respect to
$d(\cdot, \cdot)$, using the $q$-semimetric version of the triangle
inequality.  Conversely, if every Cauchy sequence of elements of $X$
with respect to $d(\cdot, \cdot)$ converges to an element of $X$ with
respect to the topology determined by $d(\cdot, \cdot)$, then $X$ is
said to be \emph{complete}\index{completeness} with respect to
$d(\cdot, \cdot)$.

        Now let $A$ be a commutative topological group.  A sequence
$\{x_j\}_{j = 1}^\infty$ of elements of $A$ is said to be a
\emph{Cauchy sequence}\index{Cauchy sequences} in $A$ if
\begin{equation}
\label{lim_{j, l to infty} (x_j - x_l) = 0}
        \lim_{j, l \to \infty} (x_j - x_l) = 0.
\end{equation}
Equivalently, this means that for each open set $U \subseteq A$ that
contains $0$, there is a positive integer $L$ such that
\begin{equation}
\label{x_j - x_l in U}
        x_j - x_l \in U
\end{equation}
for every $j, l \ge 1$.  As before, one can check that if $\{x_j\}_{j
  = 1}^\infty$ converges to an element of $A$, then $\{x_j\}_{j =
  1}^\infty$ is a Cauchy sequence in $A$.  This uses the continuity of
addition on $A$ at $0$.

        Let $d(x, y)$ be a translation-invariant $q$-semimetric on $A$
for some $q > 0$.  Suppose that $d(x, y)$ is compatible with the given
topology on $A$, in the sense that $d(x, 0)$ is continuous as a
real-valued function on $A$ at $x = 0$ with respect to the standard
topology on ${\bf R}$, as in Sections \ref{continuity of semimetrics}
and \ref{translation-invariant semimetrics}.  If $\{x_j\}_{j =
  1}^\infty$ is a Cauchy sequence in $A$ as a commutative topological
group, as in the preceding paragraph, then $\{x_j\}_{j = 1}^\infty$ is
also a Cauchy sequence with respect to $d(\cdot, \cdot)$.  Suppose for
the moment that the topology on $A$ is determined by a nonempty
collection $\mathcal{M}$ of translation-invariant $q$-semimetrics on
$A$, where $q > 0$ is allowed to depend on the element of
$\mathcal{M}$.  In this case, a sequence $\{x_j\}_{j = 1}^\infty$ of
elements of $A$ is a Cauchy sequence in $A$ as a commutative
topological group if and only if $\{x_j\}_{j = 1}^\infty$ is a Cauchy
sequence with respect to each element of $\mathcal{M}$.

        If every Cauchy sequence of elements of $A$ as a commutative
topological group converges to an element of $A$, then $A$ is said to
be \emph{sequentially complete}.\index{sequential completeness} If the
topology on $A$ is determined by a single translation-invariant
$q$-semimetric $d(x, y)$ for some $q > 0$, then $A$ is sequentially
complete as a commutative topological group if and only if $A$ is
complete with respect to $d(x, y)$.  Remember that if there is a local
base for the topology of $A$ at $0$ with only finitely or countably
many elements, then there is a translation-invariant semimetric on $A$
that determines the same topology on $A$, as in Section
\ref{metrization}.  In this case, it is reasonable to simply say that
$A$ is complete\index{completeness} as a commutative topological group
when $A$ is sequentially complete.  Otherwise, one would normally
define completeness of $A$ in terms of convergence of Cauchy nets in $A$,
or Cauchy filters on $A$.

\section{Equicontinuity}
\label{equicontinuity}
\setcounter{equation}{0}

        Let $A_1$, $A_2$ be commutative topological groups, and
let $\mathcal{L}$ be a collection of group homomorphisms from $A_1$
into $A_2$.  We say that $\mathcal{L}$ is
\emph{equicontinuous}\index{equicontinuity} at $0$ if for each open
set $U_2 \subseteq A_2$ that contains $0$ there is an open set
$U_1 \subseteq A_1$ such that $0 \in U_1$ and
\begin{equation}
\label{phi(U_1) subseteq U_2}
        \phi(U_1) \subseteq U_2
\end{equation}
for every $\phi \in \mathcal{L}$.  Of course, this implies that each
element of $\mathcal{L}$ is continuous at $0$.  If $\mathcal{L}$ has
only finitely many elements, and each element of $\mathcal{L}$ is
continuous at $0$, then $\mathcal{L}$ is equicontinuous at $0$.  If
$A_1$ is equipped with the discrete topology, then the collection of
all group homomorphisms from $A_1$ into $A_2$ is equicontinuous at
$0$.

        Now let $k$ be a field with a $q$-absolute value function $|\cdot|$
for some $q > 0$, and let $V$ be a topological vector space over $k$.
Put
\begin{equation}
\label{delta_a(v) = a v}
        \delta_a(v) = a \, v
\end{equation}
for each $a \in k$ and $v \in V$, so that $\delta_a$ is a linear
mapping from $V$ into itself for each $a \in k$.  If $|\cdot|$ is not
the trivial absolute value function on $k$, then
\begin{equation}
\label{{delta_a : a in k, |a| le 1}}
        \{\delta_a : a \in k, \, |a| \le 1\}
\end{equation}
is equicontinuous at $0$ on $V$.  This uses the fact that balanced
nonempty open subsets of $V$ form a local base for the topology of $V$
at $0$ when $|\cdot|$ is nontrivial on $k$, as in Section
\ref{topological vector spaces}.  Under these conditions, it follows
that
\begin{equation}
\label{{delta_a : a in k, |a| le |b|}}
        \{\delta_a : a \in k, \, |a| \le |b|\}
\end{equation}
is equicontinuous at $0$ for each $b \in k$.  This is the same as
saying that
\begin{equation}
\label{{delta_a : a in k, |a| le r}}
        \{\delta_a : a \in k, \, |a| \le r\}
\end{equation}
is equicontinuous at $0$ on $V$ for every positive real number $r$.
More precisely, if $|\cdot|$ is nontrivial on $k$, then there is a
$b_0 \in k$ such that $|b_0| > 1$.  Thus for each $r > 0$ there is
an $n \in {\bf Z}_+$ such that
\begin{equation}
\label{r le |b_0|^n = |b_0^n|}
        r \le |b_0|^n = |b_0^n|,
\end{equation}
so that (\ref{{delta_a : a in k, |a| le r}}) is contained in
(\ref{{delta_a : a in k, |a| le |b|}}) with $b = b_0^n$.

        If $|\cdot|$ is the trivial absolute value function on $k$,
then the equicontinuity of
\begin{equation}
\label{{delta_a : a in k}}
        \{\delta_a : a \in k\}
\end{equation}
at $0$ corresponds to (\ref{t widetilde{U} subseteq W}).  As in
Section \ref{topological vector spaces}, this is equivalent to the
condition that nonempty balanced open subsets of $V$ form a local base
for the topology of $V$ at $0$, which is not automatic in this case.

        Suppose again that $|\cdot|$ is nontrivial on $k$, and let $V_1$,
$V_2$ be topological vector spaces over $k$.  Thus the notion of
bounded subsets of $V_1$, $V_2$ can be defined as in Section
\ref{bounded sets}.  Let $\mathcal{L}$ be a collection of linear
mappings from $V_1$ into $V_2$.  If $\mathcal{L}$ is equicontinuous at
$0$, and if $E_1$ is a bounded subset of $V_1$, then one can check
that
\begin{equation}
\label{bigcup_{phi in mathcal{L}} phi(E_1)}
        \bigcup_{\phi \in \mathcal{L}} \phi(E_1)
\end{equation}
is a bounded subset of $V_2$.  If $E_1$ is an open subset of $V_1$
that contains $0$, and if (\ref{bigcup_{phi in mathcal{L}} phi(E_1)})
is a bounded subset of $V_2$, then it is easy to see that
$\mathcal{L}$ is equicontinuous at $0$.

\part{Topological dimension $0$}
\label{topological dimension 0}

\section{Separated sets}
\label{separated sets}
\setcounter{equation}{0}

        As usual, a pair $A$, $B$ of subsets of a topological space $X$
are said to be \emph{separated}\index{separated sets} in $X$ if
\begin{equation}
\label{overline{A} cap B = A cap overline{B} = emptyset}
        \overline{A} \cap B = A \cap \overline{B} = \emptyset,
\end{equation}
where $\overline{A}$, $\overline{B}$ are the closures of $A$, $B$ in
$X$, respectively.  In particular, disjoint closed subsets of $X$ are
separated, and it is easy to see that disjoint open sets in $X$ are
separated too.  If $A, B \subseteq X$ are separated and $A \cup B =
X$, then both $A$ and $B$ are both open and closed in $X$.  If $Y
\subseteq X$, then $Y$ may also be considered as a topological space,
with respect to the topology induced by the one on $X$.  It is well
known and easy to check that $A, B \subseteq Y$ are separated with
respect to the induced topology on $Y$ if and only if $A$, $B$ are
separated as subsets of $X$.  A set $E \subseteq X$ is said to be
\emph{connected}\index{connected sets} if it cannot be expressed as
the union of two nonempty separated sets.  Thus $X$ is connected if it
cannot be expressed as the union of two nonempty disjoint closed sets,
or equivalently as the union of two nonempty disjoint open sets.  If
$E \subseteq Y \subseteq X$, then $E$ is connected with respect to the
induced topology on $Y$ if and only if $E$ is connected as a subset of
$X$.  One can check that the closure of a connected subset of $X$ is
also connected in $X$.  A subset of $X$ is said to be \emph{totally
  disconnected}\index{totally disconnected sets} if it does not
contain any connected set with at least two elements.

        Now let $A$ be a commutative topological group, and let
$U \subseteq A$ be an open set that contains $0$.  It will be convenient
to ask that $U$ also be symmetric about $0$ in $A$, in the sense that
\begin{equation}
\label{-U = U}
        -U = U,
\end{equation}
which can always be arranged by replacing $U$ with $U \cap (-U)$.  Let
us say that $B, C \subseteq A$ are
\emph{$U$-separated}\index{separated sets} in $A$ if
\begin{equation}
\label{(B + U) cap C = emptyset}
        (B + U) \cap C = \emptyset,
\end{equation}
which is equivalent to asking that
\begin{equation}
\label{B cap (C + U) = emptyset}
        B \cap (C + U) = \emptyset,
\end{equation}
because of (\ref{-U = U}).  Using the continuity of the group
operations on $A$ at $0$, it is easy to see that there is an open set
$U_1 \subseteq A$ that contains $0$, is symmetric about $0$, and satisfies
\begin{equation}
\label{U_1 + U_1 subseteq U}
        U_1 + U_1 \subseteq U,
\end{equation}
as in (\ref{U + V subseteq W}).  Combining this with (\ref{(B + U) cap
  C = emptyset}) or (\ref{B cap (C + U) = emptyset}), one can check
that
\begin{equation}
\label{(B + U_1) cap (C + U_1) = emptyset}
        (B + U_1) \cap (C + U_1) = \emptyset.
\end{equation}
Remember that the closures of $B$, $C$ in $A$ are contained in $B +
U_1$, $C + U_1$, respectively, as in (\ref{overline{E} subseteq E +
  V}).  Thus (\ref{(B + U_1) cap (C + U_1) = emptyset}) implies in
particular that the closures of $B$ and $C$ in $A$ are disjoint.

        Suppose that $C \subseteq A$ is compact, $W \subseteq A$
is an open set, and that $C \subseteq W$.  If $x \in C$, then there
is an open set $U(x) \subseteq A$ that contains $0$ and satisfies
\begin{equation}
\label{x + U(x) + U(x) subseteq W}
        x + U(x) + U(x) \subseteq W,
\end{equation}
by continuity of addition on $A$.  The collection of open sets $x +
U(x)$ with $x \in C$ covers $C$, and so there are finitely many
elements $x_1, \ldots, x_n$ of $C$ such that
\begin{equation}
\label{C subseteq bigcup_{j = 1}^n (x_j + U(x_j))}
        C \subseteq \bigcup_{j = 1}^n (x_j + U(x_j)),
\end{equation}
by compactness.  Put
\begin{equation}
\label{U = bigcap_{j = 1}^n (U(x_j) cap (-U(x_j)))}
        U = \bigcap_{j = 1}^n (U(x_j) \cap (-U(x_j))),
\end{equation}
which is an open subset of $A$ that contains $0$ and is symmetric
about $0$.  Using (\ref{x + U(x) + U(x) subseteq W}) and (\ref{C
  subseteq bigcup_{j = 1}^n (x_j + U(x_j))}), we get that
\begin{equation}
\label{C + U subseteq ... subseteq W}
        C + U \subseteq \bigcup_{j = 1}^n (x_j + U(x_j) + U)
               \subseteq \bigcup_{j = 1}^n (x_j + U(x_j) + U(x_j)) \subseteq W.
\end{equation}
If $E \subseteq A$ is a closed set that is disjoint from $C$, then we
can apply the previous argument to $W = A \setminus E$.  In this case,
(\ref{C + U subseteq ... subseteq W}) is equivalent to saying that
\begin{equation}
\label{(C + U) cap E = emptyset}
        (C + U) \cap E = \emptyset,
\end{equation}
so that $C$ and $E$ are $U$-separated in $A$.

        Let $X$ be a topological space that is regular in the strict
sense, without including the first or $0$th separation condition.  If
$K \subseteq X$ is compact, $W \subseteq X$ is an open set, and $K
\subseteq W$, then there is an open set $U \subseteq X$ such that $K
\subseteq U$ and $\overline{U} \subseteq W$.  More precisely, for each
$x \in K$, there is an open set $U(x) \subseteq X$ such that $x \in
U(x)$ and $\overline{U(x)} \subseteq W$, because $X$ is regular.
Because $K$ is compact, $K$ can be covered by finitely many $U(x)$'s,
whose union $U$ has the desired properties.  In particular, if $K
\subseteq X$ is compact and open, then one can apply the previous
argument with $K = W$, to conclude that $K$ is a closed set in $X$.
Of course, if $X$ is Hausdorff, then every compact subset of $X$
is closed.  However, if $X$ is regular in the strict sense, then
$X$ need not be Hausdorff, and arbitrary compact subsets of $X$
need not be closed.  If $X$ is any set equipped with the indiscrete
topology, for instance, then $X$ is regular in the strict sense,
and every subset of $X$ is compact.

\section{Open subgroups}
\label{open subgroups}
\setcounter{equation}{0}

        Let $A$ be a commutative topological group again, and let $B$
be a subgroup of $A$ which is also an open set in $A$.  It is well
known that the complement of $B$ in $A$ can be expressed as a union
of cosets of $B$ in $A$, which are translates of $B$ in $A$.  Each
translate of $B$ is also an open set in $A$, and hence any union of
translates of $B$ in $A$ is an open set in $A$ too.  This implies
that the complement of $B$ in $A$ is an open set in $A$, which means
that $B$ is a closed set in $A$ as well.  In particular, if $A$ is
connected as a topological space, then $A$ has no proper open
subgroups.

        If $B$ is any subgroup of $A$, then
\begin{eqnarray}
\label{0 in B}
        0 \in B, \\
\label{-B = B}
        -B = B,
\end{eqnarray}
and
\begin{equation}
\label{B + B = B}
        B + B = B.
\end{equation}
Conversely, if $B$ is any subset of $A$ that satisfies (\ref{0 in B}),
(\ref{-B = B}), and (\ref{B + B = B}), then $B$ is a subgroup of $A$.
If $B$ is an open subgroup of $A$, then (\ref{B + B = B}) implies that
$B$ is $B$-separated from its complement in $A$, as discussed in the
previous section.

        Suppose that $U \subseteq A$ contains $0$ and is symmetric
about $0$.  Put $U_1 = U$, and define $U_j$ recursively for
$j \in {\bf Z}_+$ by putting
\begin{equation}
\label{U_{j + 1} = U_j + U}
        U_{j + 1} = U_j + U
\end{equation}
for each $j \ge 1$.  Equivalently, $U_j$ is the subset of $A$
consisting of sums of exactly $j$ elements of $U$ for each $j$.  Note
that $0 \in U_j$, $U_j$ is symmetric about $0$, and $U_j \subseteq
U_{j + 1}$ for each $j$, and that
\begin{equation}
\label{U_{j + l} = U_j + U_l}
        U_{j + l} = U_j + U_l
\end{equation}
for each $j, l \ge 1$.  If we put
\begin{equation}
\label{B = bigcup_{j = 1}^infty U_j}
        B = \bigcup_{j = 1}^\infty U_j,
\end{equation}
then $B$ satisfies (\ref{0 in B}), (\ref{-B = B}), and (\ref{B + B =
  B}).  Thus $B$ is a subgroup of $A$.  If $U$ is an open set in $A$,
then $U_j$ is an open set in $A$ for each $j$, and hence $B$ is an
open set in $A$ too.

        Let $U \subseteq A$ be an open set that contains $0$ and is
symmetric about $0$, and suppose that $E$ is a nonempty subset of
$A$ that is $U$-separated from its complement in $A$.  This means
that
\begin{equation}
\label{(E + U) cap (A setminus E) = emptyset}
        (E + U) \cap (A \setminus E) = \emptyset,
\end{equation}
as in the previous section, which is the same as saying that
\begin{equation}
\label{E + U subseteq E}
        E + U \subseteq E.
\end{equation}
If $U_j$ is defined as in the preceding paragraph for each $j \in {\bf
  Z}_+$, then it follows that
\begin{equation}
\label{E + U_j subseteq E}
        E + U_j \subseteq E
\end{equation}
for each $j \ge 1$.  This implies that
\begin{equation}
\label{E + B subseteq E}
        E + B \subseteq E,
\end{equation}
where $B$ is as in (\ref{B = bigcup_{j = 1}^infty U_j}), by taking the
union over $j \in {\bf Z}_+$ in (\ref{E + U_j subseteq E}).  Of
course, we actually have equality in (\ref{E + U subseteq E}), (\ref{E
  + U_j subseteq E}), and (\ref{E + B subseteq E}), because $0$ is an
element of $U$, and hence of $U_j$ for each $j$, as well as $B$.

\section{Semimetrics and partitions}
\label{semimetrics, partitions}
\setcounter{equation}{0}

        Let $X$ be a set, and let $\mathcal{P}$ be a partition of
$X$.\index{partitions}  Thus $\mathcal{P}$ is a collection of
pairwise-disjoint nonempty subsets of $X$ whose union is equal to
$X$.  Define
\begin{equation}
\label{d_mathcal{P}(x, y)}
        d_\mathcal{P}(x, y)
\end{equation}
for $x, y \in X$ by putting (\ref{d_mathcal{P}(x, y)}) equal to $0$
when $x$ and $y$ are contained in the same element of $\mathcal{P}$,
and equal to $1$ otherwise.  It is easy to see that this defines a
semi-ultrametric on $X$, which one might describe as the
\emph{discrete semimetric}\index{discrete semimetrics} associated to
$\mathcal{P}$.  If $\mathcal{P}$ is the partition of $X$ consisting of
all subsets of $X$ with exactly one element, then
(\ref{d_mathcal{P}(x, y)}) is the same as the discrete metric on $X$.
If $\mathcal{P}$ consists of only $X$ itself, then let us call
$\mathcal{P}$ the \emph{trivial partition}\index{trivial partition} of
$X$.  In this case, (\ref{d_mathcal{P}(x, y)}) is equal to $0$ for
every $x, y \in X$.

        Let $Z$ be another set, and let $\phi$ be a mapping from
$X$ into $Z$.  If $d_Z(\cdot, \cdot)$ is a $q$-semimetric on $Z$
for some $q > 0$, then
\begin{equation}
\label{d_Z(phi(x), phi(y))}
        d_Z(\phi(x), \phi(y))
\end{equation}
defines a $q$-semimetric on $X$.  Note that
\begin{equation}
\label{{phi^{-1}(z) : z in phi(X)}}
        \{\phi^{-1}(z) : z \in \phi(X)\}
\end{equation}
is a partition of $X$ under these conditions.  If $d_Z(\cdot, \cdot)$
is the discrete metric on $Z$, then (\ref{d_Z(phi(x), phi(y))}) is the
same as the discrete semimetric associated to the partition
(\ref{{phi^{-1}(z) : z in phi(X)}}), as in the preceding paragraph.
Of course, if $\mathcal{P}$ is any partition of $X$, then
$\mathcal{P}$ is of the form (\ref{{phi^{-1}(z) : z in phi(X)}}),
where $Z = \mathcal{P}$ and $\phi$ is the mapping that sends $x \in X$
to the element of $\mathcal{P}$ that contains $x$.

        If $d(x, y)$ is a $q$-semimetric on $X$ for any $q > 0$, then
\begin{equation}
\label{d(x, y) = 0}
        d(x, y) = 0
\end{equation}
defines an equivalence relation on $X$, whose equivalence classes
determine a partition of $X$.  These equivaence classes are the same
as the closed balls in $X$ of radius $0$ with respect to $d(\cdot,
\cdot)$.  Similarly, if $d(x, y)$ is a semi-ultrametric on $X$, then
for each $r > 0$, $X$ is partitioned by the open balls in $X$ with
respect to $d(\cdot, \cdot)$ of radius $r$, and by the closed balls in
$X$ of radius $r$.  In the case of (\ref{d_mathcal{P}(x, y)}), the
elements of $\mathcal{P}$ are the same as the open balls of radius $0
< r \le 1$, and the closed balls of radius $0 \le r < 1$.  If $d(x,
y)$ is a $q$-semimetric on $X$ for some $q > 0$ that takes values in
the set $\{0, 1\}$, then $d(x, y)$ is the same as the discrete
semimetric associated to the partition of $X$ consisting of the closed
balls of radius $0$.

        Let $X$ be a topological space, and let $\mathcal{P}$ be a
partition of $X$.  If every element of $\mathcal{P}$ is an open subset
of $X$, then every element of $\mathcal{P}$ is a closed set in $X$
too, because the complement of every element of $\mathcal{P}$ in $X$
is equal to the union of the other elements of $\mathcal{P}$, and
hence is an open subset of $X$.  In particular, if $X$ is connected,
then $\mathcal{P}$ is the trivial partition of $X$ under these
conditions.  Note that every element of $\mathcal{P}$ is an open
subset of $X$ exactly when the corresponding discrete semimetric
(\ref{d_mathcal{P}(x, y)}) is compatible with the topology on $X$, in
the sense that open subsets of $X$ with respect to
(\ref{d_mathcal{P}(x, y)}) are also open with respect to the given
topology on $X$.  If $Z$ is any set equipped with the discrete
topology, then a mapping $\phi$ from $X$ into $Z$ is continuous if and
only if (\ref{{phi^{-1}(z) : z in phi(X)}}) is an open set in $X$ for
every $z \in Z$.

\section{Dimension $0$}
\label{dimension 0}
\setcounter{equation}{0}

        A topological space $X$ is said to have \emph{topological
dimension $0$}\index{topological dimension 0@topological dimension $0$}
at a point $x \in X$ if for each open set $W \subseteq X$ that
contains $x$ there is an open set $W_1 \subseteq X$ such that $x \in
W_1$, $W_1 \subseteq W$, and $W_1$ is also a closed set in $X$.
Equivalently, this means that there is a local base for the topology
of $X$ at $x$ consisting of subsets of $X$ that are both open and
closed.  If $X$ has topological dimension $0$ at every point $x \in
X$, then $X$ is said to have topological dimension $0$ as a
topological space.  As before, this is the same as saying that there
is a base for the topology of $X$ consisting of subsets of $X$ that
are both open and closed.  Sometimes $X$ would also be required to be
nonempty to have topological dimension $0$, and the empty set is
defined to have topological dimension $-1$, in order to define
topological dimension inductively when it is positive.

        Let $X$ be a topological space with topological dimension $0$.
Observe that every subset $Y$ of $X$ has topological dimension $0$
with respect to the induced topology.  Clearly $X$ is regular as a
topological space in the strict sense, without including the first or
$0$th separation condition.  One may wish to include the first or
$0$th separation condition as part of the definition of topological
dimension $0$, in which case the space is Hausdorff as well.  If $X$
has at least two elements, and if $X$ satisfies the first or $0$th
separation condition, then it is easy to see that $X$ is not
connected.  It follows that $X$ is totally disconnected when $X$ has
topological dimension $0$ and $X$ satisfies the first or $0$th
separation condition, since every subset of $X$ has the same
properties with respect to the induced topology.  Otherwise, any set
equipped with the indiscrete topology has topological dimension $0$ in
the strict sense, without including the first or $0$th separation
condition.

        Suppose that $d(x, y)$ is a semi-ultrametric on a set $X$.
As in Section \ref{semi-ultrametrics}, both open and closed balls
in $X$ with respect to $d$ of positive radius are both open and
closed as subsets of $X$, with respect to the topology on $X$
determined by $d$.  This implies that $X$ has topological dimension
$0$ in the strict sense, without including the first or $0$th separation
condition, with respect to the topology determined by $d$.  Similarly,
if $\mathcal{M}$ is a nonempty collection of semi-ultrametrics on $X$,
then $X$ has topological dimension $0$ in the strict sense with respect
to the topology determined by $\mathcal{M}$.  Of course, if $\mathcal{M}$
is nondegenerate on $X$, then $X$ is Hausdorff with respect to the
topology determined by $\mathcal{M}$.

        Let $E$ be a subset of a set $X$, and put
\begin{equation}
\label{d_E(x, y) = 0}
        d_E(x, y) = 0
\end{equation}
when $x, y \in E$ and when $x, y \in X \setminus E$, and put
\begin{equation}
\label{d_E(x, y) = 1}
        d_E(x, y) = 1
\end{equation}
when $x \in E$ and $y \in X \setminus E$, and when $x \in X \setminus
E$ and $y \in E$.  It is easy to see that this defines a
semi-ultrametric on $X$.  More precisely, if $E$ is a nonempty
proper subset of $X$, then
\begin{equation}
\label{{E, X setminus E}}
        \{E, X \setminus E\}
\end{equation}
is a partition of $X$, and $d_E(x, y)$ is the same as the discrete
semimetric associated to this partition, as in the preceding section.
Otherwise, if $E = \emptyset$ or $E = X$, then (\ref{d_E(x, y) = 0})
holds for every $x, y \in X$, and $d_E(x, y)$ is the same as the
discrete semimetric associated to the trivial partition of $X$.
Suppose now that $X$ is a topological space.  If $E \subseteq X$ is
both open and closed, then $d_E$ is compatible with the topology on
$X$, in the sense that every open set in $X$ with respect to $d_E$ is
also an open set with respect to the given topology on $X$.  If
\begin{equation}
\label{mathcal{M}_0 = {d_E : E subseteq X is both open and closed}}
 \mathcal{M}_0 = \{d_E : E \subseteq X \hbox{ is both open and closed}\},
\end{equation}
is the collection of all such semi-ultrametrics on $X$, then it
follows that every open set in $X$ with respect to the topology
determined by $\mathcal{M}_0$ is an open set with respect to the given
topology on $X$ too.  If $X$ has topological dimension $0$ in the
strict sense, then one can check that the topology on $X$ determined
by $\mathcal{M}_0$ is the same as the given topology on $X$.  If $X$
also satisfies the first or $0$th separation condition, then
$\mathcal{M}_0$ is nondegenerate on $X$.

\section{Totally separated topological spaces}
\label{totally separated topological spaces}
\setcounter{equation}{0}

        A topological space $X$ is said to be \emph{totally
separated}\index{totally separated topological spaces}
if for every pair of points $x, y \in X$ with $x \ne y$ there is an
open set $U \subseteq X$ such that $x \in U$, $y \in X \setminus U$,
and $U$ is also closed in $X$.  Note that this is symmetric in $x$ and
$y$, since $X \setminus U$ is both open and closed in $X$ too.
Suppose that $X$ is totally separated, which implies that $X$ is
Hausdorff.  It is easy to see that every subset of $X$ is totally
separated with respect to the induced topology.  If $X$ has at least
two elements, then $X$ is not connected.  This implies that $X$ is
totally disconnected, because subspaces of $X$ are totally separated.
If $\widetilde{\tau}$ is another topology on $X$ which contains the
given topology on $X$, then $X$ is totally separated with respect to
$\widetilde{\tau}$ too.

        Let $(X, \tau)$ be a topological space, and let $\tau_0$ be
the collection of subsets $W$ of $X$ with the property that for each
$x \in W$ there is an open set $U \subseteq X$ with respect to $\tau$
such that $x \in U$, $U \subseteq W$, and $U$ is closed in $X$ with
respect to $\tau$.  Equivalently, this means that $W$ can be expressed
as the union of a family of subsets of $X$ that are both open and
closed with respect to $\tau$.  In particular, this implies that $W$
is open with respect to $\tau$, so that
\begin{equation}
\label{tau_0 subseteq tau}
        \tau_0 \subseteq \tau.
\end{equation}
It is easy to see that $\tau_0$ is a topology on $X$, which is the
same as the topology determined by the collection (\ref{mathcal{M}_0 =
  {d_E : E subseteq X is both open and closed}}) of semi-ultrametrics
on $X$ associated to $\tau$.  By construction, if $U \subseteq X$ is
both open and closed with respect to $\tau$, then $U$ is an open set
with respect to $\tau_0$.  In this case, $X \setminus U$ is both open
and closed with respect to $\tau$ too, so that $X \setminus U$ is an
open set with respect to $\tau_0$, which means that $U$ is a closed
set with respect to $\tau_0$.  Conversely, if $U \subseteq X$ is both
open and closed with respect to $\tau_0$, then $U$ is both open and
closed with respect to $\tau$, because of (\ref{tau_0 subseteq tau}).
The collection of these sets forms a base for $\tau_0$, by definition
of $\tau_0$.  This implies that $X$ automatically has topological
dimension $0$ with respect to $\tau_0$, which also follows from the
characterization of $\tau_0$ in terms of (\ref{mathcal{M}_0 = {d_E : E
    subseteq X is both open and closed}}).  One can check that $X$ is
totally separated with respect to $\tau$ if and only if $X$ is
Hausdorff with respect to $\tau_0$.

        Let $X$ be a totally separated topological space.  Also let
$H$ be a compact subset of $X$, and let $y \in X \setminus H$ be
given.  If $x \in H$, then $x \ne y$, and so there is an open set
$U(x) \subseteq X$ such that $x \in U(x)$, $y \in X \setminus U(x)$,
and $U(x)$ is a closed set in $X$.  Because $H$ is compact, there are
finitely many elements $x_1, \ldots, x_n$ of $H$ such that
\begin{equation}
\label{H subseteq bigcup_{j = 1}^n U(x_j)}
        H \subseteq \bigcup_{j = 1}^n U(x_j).
\end{equation}
The right side of (\ref{H subseteq bigcup_{j = 1}^n U(x_j)}) is both
open and closed in $X$, and does not contain $y$.  Similarly, if $H$
and $K$ are disjoint compact subsets of $X$, then there is an open set
$U \subseteq X$ such that $H \subseteq U$, $K \subseteq X \setminus
U$, and $U$ is a closed set in $X$.  If $X$ is a topological space
with topological dimension $0$, $K \subseteq X$ is compact, $W
\subseteq X$ is open, and $K \subseteq W$, then an analogous argument
implies that there is an open set $U \subseteq X$ such that $K
\subseteq U \subseteq W$ and $U$ is a closed set in $X$.

        Let $X$ be a totally separated topological space again,
and suppose that $X$ is locally compact.  We would like to check that
$X$ has topological dimension $0$ under these conditions.  To do this,
let $x \in X$ and an open set $W \subseteq X$ that contains $x$ be
given.  Because $X$ is locally compact, there is an open set in
$X$ that contains $x$ and is contained in a compact subset of $X$.
We may as well suppose that $W$ is contained in a compact subset of
$X$, since otherwise we can replace $W$ with its intersection with
an open set that contains $x$ and is contained in a compact set.
Note that compact subsets of $X$ are closed sets, because $X$ is
Hausdorff, since it is totally separated.  It follows that the
closure $\overline{W}$ of $W$ in $X$ is contained in a compact set,
and hence that $\overline{W}$ is compact in $X$, because closed
subsets of compact sets are compact.  Similarly,
\begin{equation}
\label{partial W = overline{W} setminus W}
        \partial W = \overline{W} \setminus W
\end{equation}
is compact, and $x \not\in \partial W$, since $x \in W$.  As in the
previous paragraph, it follows that there is an open set $U_1
\subseteq X$ such that $\partial W \subseteq U_1$, $x \in X \setminus
U_1$, and $U_1$ is a closed set in $X$, because $X$ is totally
separated.  Put
\begin{equation}
\label{U_2 = U_1 cup (X setminus W) = U_1 cup (X setminus overline{W})}
        U_2 = U_1 \cup (X \setminus W) = U_1 \cup (X \setminus \overline{W}),
\end{equation}
where the second step uses the fact that $\partial W \subseteq U_1$.
It is easy to see that $U_2$ is both open and closed in $X$, because
of the analogous property of $U_1$, and the two expressions for $U_2$
in (\ref{U_2 = U_1 cup (X setminus W) = U_1 cup (X setminus
  overline{W})}).  Thus $X \setminus U_2$ is both open and closed in
$X$, and $X \setminus U_2 \subseteq W$, since $X \setminus W \subseteq
U_2$, by the definition of $U_2$.  We also have that $x \in X
\setminus U_2$, because $x \in X \setminus U_1$, by the way that $U_1$
was chosed, and $x \in W$, by hypothesis.  This shows that $X$ has
topological dimension $0$ at $x$, and hence that $X$ has topological
dimension $0$, since $x \in X$ is arbitrary.  It is well known that
any locally compact Hausdorff topological space that is totally
disconnected has topological dimension $0$.

\section{Dimension $0$, continued}
\label{dimension 0, continued}
\setcounter{equation}{0}

        Let $A$ be a commutative topological group.  If $A$ has
topogical dimension $0$ at any of its elements, then $A$ has
topological dimension $0$ at each point, because of continuity of
translations.  Suppose for the moment that the collection of open
subgroups of $A$ is a local base for the topology of $A$ at $0$.
Equivalently, this means that for open set $W \subseteq A$ that
contains $0$, there is an open subgroup $B$ of $A$ such that
\begin{equation}
\label{B subseteq W}
        B \subseteq W.
\end{equation}
This implies that $A$ has topological dimension $0$ at $0$, and hence
at every point in $A$, because open subgroups of $A$ are closed sets,
as in Section \ref{open subgroups}.  Note that the intersection of
finitely many open subgroups of $A$ is an open subgroup of $A$ too.
If there is a local sub-base of $A$ at $0$ consisting of open
subgroups of $A$, then it follows that the open subgroups of $A$ form
a local base for the topology of $A$ at $0$.  If the open subgroups of
$A$ form a local base for the topology of $A$ at $0$, then every
subgroup of $A$ has the same property with respect to the induced
topology.

        Suppose now that for each open set $W \subseteq A$ that contains
$0$ there are subsets $E$, $U$ of $A$ containing $0$ such that
$E \subseteq W$, $U$ is an open set in $A$, $U$ is symmetric about $0$,
and $E$ is $U$-separated from $A \setminus E$, as in Section
\ref{separated sets}.  In particular, this implies that $E$ and $A
\setminus E$ are separated in $A$, as before.  It follows that $E$ is
both open and closed in $A$, so that this condition implies that $A$
has topological dimension $0$ at $0$.  Let $W$, $E$, and $U$ be given
as in this condition, and let $B$ be obtained from $U$ as in (\ref{B =
  bigcup_{j = 1}^infty U_j}).  Thus $B$ is an open subgroup of $A$, as
in Section \ref{open subgroups}.  Because $E$ is $U$-separated from $A
\setminus E$, we also have that $E + B$ is contained in $E$, as in
(\ref{E + B subseteq E}).  This implies that
\begin{equation}
\label{B subseteq E}
        B \subseteq E,
\end{equation}
since $0 \in E$, and hence that (\ref{B subseteq W}) holds, because $E
\subseteq W$ by hypothesis.  This shows that the condition mentioned
at the beginning of the paragraph implies that the collection of open
subgroups of $A$ is a local base for the topology of $A$ at $0$.  Of
course, if $B$ is any open subgroup of $A$, then $B$ is $B$-separated
from $A \setminus B$, because $B + B$ is contained in $B$.  If the
collection of open subgroups of $A$ is a local base for the topology
of $A$ at $0$, then it follows that $A$ has the property mentioned at
the beginning of the paragraph.

        Suppose that $A$ is locally compact at $0$, which is to say
that there is an open set $W_0 \subseteq A$ that contains $0$ and is
contained in a compact set $K \subseteq A$.  Suppose also that $A$ has
topological dimension $0$ at $0$, and let $W \subseteq A$ be any open
set that contains $0$.  Thus $W \cap W_0$ is an open subset of $A$
that contains $0$, and so there is an open set $E \subseteq A$ such
that $0 \in E$,
\begin{equation}
\label{E subseteq W cap W_0}
        E \subseteq W \cap W_0,
\end{equation}
and $E$ is a closed set in $A$.  In particular,
\begin{equation}
\label{E subseteq W_0 subseteq K}
        E \subseteq W_0 \subseteq K,
\end{equation}
which implies that $E$ is compact in $A$, because $E$ is a closed set
and $K$ is compact.  Remember that $E$ is an open set in $A$ too, so
that $A \setminus E$ is a closed set.  Using an argument from Section
\ref{separated sets}, we get that there is an open set $U \subseteq A$
such that $0 \in U$, $U$ is symmetric about $0$, and $E$ is
$U$-separated from $A \setminus E$.  This shows that $A$ satisfies the
condition mentioned at the beginning of the previous paragraph under
these conditions, so that the collection of open subgroups of $A$ is a
local base for the topology of $A$ at $0$.

        If $A$ is any commutative group equipped with the discrete
topology, then $\{0\}$ is an open subgroup of $A$, which defines a
local base for the discrete topology on $A$.  Consider now the set
${\bf Q}$ of rational numbers as a commutative topological group with
respect to addition and the topology induced by the standard topology
on the real line.  It is easy to see that ${\bf Q}$ has topological
dimension $0$ as a topological space.  However, one can also check
that ${\bf Q}$ is the only open subgroup of itself.

\section{Translation-invariant semi-ultrametrics}
\label{translation-invariant semi-ultrametrics}
\setcounter{equation}{0}

        Let $A$ be a commutative group, and let $B$ be a subgroup
of $A$.  The collection of cosets of $B$ in $A$ defines a partition
$\mathcal{P}_B(A)$ of $A$.  Let
\begin{equation}
\label{d_{mathcal{P}_B(A)}(x, y)}
        d_{\mathcal{P}_B(A)}(x, y)
\end{equation}
be the discrete semi-ultrametric on $A$ corresponding to
$\mathcal{P}_B(A)$ as in Section \ref{semimetrics, partitions}.  In
this case, (\ref{d_{mathcal{P}_B(A)}(x, y)}) is equal to $0$ when
\begin{equation}
\label{x - y in B}
        x - y \in B,
\end{equation}
and otherwise (\ref{d_{mathcal{P}_B(A)}(x, y)}) is equal to $1$.  In
particular, (\ref{d_{mathcal{P}_B(A)}(x, y)}) is invariant under
translations on $A$.  Note that $\mathcal{P}_B(A)$ corresponds to the
standard quotient mapping from $A$ onto the quotient group $A / B$ as
in (\ref{{phi^{-1}(z) : z in phi(X)}}).  Thus
(\ref{d_{mathcal{P}_B(A)}(x, y)}) corresponds to the discrete metric
on $A / B$ and the standard quotient mapping from $A$ onto $A / B$
as in (\ref{d_Z(phi(x), phi(y))}).

        If $d(x, y)$ is a $q$-semimetric on $A$ for some $q > 0$,
then we have seen that (\ref{d(x, y) = 0}) defines an equivalence
relation on $A$, for which the corresponding equivalence classes are
the same as the closed balls in $A$ of radius $0$.  If $d(x, y)$ is
invariant under translations on $A$, then one can check that the
equivalence class containing $0$ is a subgroup of $A$, and the other
equivalence classes are cosets of this subgroup in $A$.  Similarly, if
$d(x, y)$ is a semi-ultrametric on $A$ and $r > 0$, then $A$ is
partitioned by the open balls of radius $r$ with respect to $d(x, y)$,
and by the closed balls of radius $r$.  If $d(x, y)$ is invariant
under translations, then one can verify that the open and closed balls
in $A$ with respect to $d(x, y)$ centered at $0$ with radius $r$ are
subgroups of $A$ for every $r > 0$.  Of course, the open and closed
balls in $A$ with respect to $d(x, y)$ centered at other points are
cosets of the corresponding balls centered at $0$.

        Let $\mathcal{M}$ be a nonempty collection of translation-invariant
semi-ultrametrics on $A$, and remember that $A$ is a commutative
topological group with respect to the topology determined by
$\mathcal{M}$, as in Section \ref{translation-invariant semimetrics}.
As in the preceding paragraph, the open and closed balls in $A$
centered at $0$ with respect to elements of $\mathcal{M}$ are
subgroups of $A$.  Of course, open balls in $A$ with respect to
elements of $\mathcal{M}$ are open sets with respect to the
corresponding topology.  In this situation, closed balls in $A$ of
positive radius with respect to elements of $\mathcal{M}$ are open
sets too, because the elements of $\mathcal{M}$ are semi-ultrametrics,
as in Section \ref{semi-ultrametrics}.  Remember that open balls in
$A$ centered at $0$ with respect to $\mathcal{M}$ form a local
sub-base for the topology on $A$ determined by $\mathcal{M}$, as in
Section \ref{collections of semimetrics}.  This is a local sub-base
for the topology of $A$ at $0$ consisting of open subgroups of $A$,
since open balls in $A$ centered at $0$ with respect to elements of
$\mathcal{M}$ are subgroups of $A$, as before.  It follows that the
open subgroups of $A$ form a local base for the topology of $A$ at $0$
under these conditions, as in the previous section.

        Let $\mathcal{B}$ be a nonempty collection of subgroups of $A$,
and let
\begin{equation}
\label{mathcal{M}(mathcal{B}) = {d_{mathcal{P}_B(A)} : B in mathcal{B}}}
 \mathcal{M}(\mathcal{B}) = \{d_{\mathcal{P}_B(A)} : B \in \mathcal{B}\}
\end{equation}
be the collection of discrete semi-ultrametrics on $A$ corresponding
to the partitions $\mathcal{P}_B(A)$ of $A$ associated to elements $B$
of $\mathcal{B}$ as in (\ref{d_{mathcal{P}_B(A)}(x, y)}).  Remember
that these semi-ultrametrics are invariant under translations on $A$.
This implies that $A$ is a commutative topological group with respect
to the topology determined by $\mathcal{M}(\mathcal{B})$, and the open
subgroups of $A$ form a local base for the topology of $A$ at $0$, as
in the preceding paragraph.  More precisely, each $B \in \mathcal{B}$
is the same as the open ball in $A$ centered at $0$ with radius $1$
with respect to the corresponding discrete semi-ultrametric
$d_{\mathcal{P}_B(A)}$, by construction.  Using this, it is easy to
see that $\mathcal{B}$ is a local sub-base for the topology of $A$ at
$0$ with respect to the topology determined by
$\mathcal{M}(\mathcal{B})$.

        Let $A$ be a commutative topological group, and put
\begin{equation}
\label{mathcal{B}(A) = {B subseteq A : B is an open subgroup of A}}
 \mathcal{B}(A) = \{B \subseteq A : B \hbox{ is an open subgroup of } A\}.
\end{equation}
Note that $A$ is automatically an element of $\mathcal{B}(A)$, so that
$\mathcal{B}(A) \ne \emptyset$.  Also let
$\mathcal{M}(\mathcal{B}(A))$ be the collection of discrete
semi-metrics corresponding to the partitions $\mathcal{P}_B(A)$
associated to $B \in \mathcal{B}(A)$, as in
(\ref{mathcal{M}(mathcal{B}) = {d_{mathcal{P}_B(A)} : B in
    mathcal{B}}}).  If $B \in \mathcal{B}(A)$, then the corresponding
discrete semi-metric $d_{\mathcal{P}_B(A)}$ is compatible with the
given topology on $A$, in the sense that every open set in $A$ with
respect to $d_{\mathcal{P}_B(A)}$ is also an open set with respect to
the given topology on $A$.  This implies that every open set in $A$
with respect to the topology determined by
$\mathcal{M}(\mathcal{B}(A))$ is an open set with respect to the given
topology on $A$.  Of course, every element of $\mathcal{B}(A)$ is an
open set in $A$ with respect to the topology determined by
$\mathcal{M}(\mathcal{B}(A))$.  If $\mathcal{B}(A)$ is a local base
for the given topology on $A$ at $0$, then the topology on $A$
determined by $\mathcal{M}(\mathcal{B}(A))$ is the same as the given
topology on $A$.

        Suppose that for each $y \in A$ with $y \ne 0$ there is an
open set $U \subseteq A$ such that $0 \in U$, $y \in A \setminus U$,
and $U$ is a closed set in $A$ too.  This is the same as the totally
separated property described in Section \ref{totally separated
  topological spaces} with $x = 0$.  In this situation, this property
implies that $A$ is totally separated as a topological space, because
of continuity of translations.

        Suppose now that for each $y \in A$ with $y \ne 0$ there is
an open subgroup $B$ of $A$ such that $y \not\in B$.  This implies the
condition mentioned in the previous paragraph, because open subgroups
of $A$ are closed sets in $A$, as in Section \ref{open subgroups}.
Equivalently, this new condition is the same as saying that the
intersection of all of the open subgroups of $A$ is equal to $\{0\}$.
This is also equivalent to asking that $\mathcal{M}(\mathcal{B}(A))$
be nondegenerate on $A$.  If $A$ has this property, then every
subgroup of $A$ has the same property with respect to the induced
topology.

\section{Cartesian products, continued}
\label{cartesian products, continued}
\setcounter{equation}{0}

        Let $I$ be a nonempty set, and let $Y_j$ be a topological space
for each $j \in I$.  Thus the Cartesian product
\begin{equation}
\label{Y = prod_{j in I} Y_j}
        Y = \prod_{j \in I} Y_j
\end{equation}
is also a topological space with respect to the corresponding product
and strong product topologies.  If $E_j \subseteq Y_j$ is a closed set
for each $j \in I$, then
\begin{equation}
\label{E = prod_{j in I} E_j, 3}
        E = \prod_{j \in I} E_j
\end{equation}
is a closed set in $Y$ with respect to the product topology, and hence
with respect to the strong product topology.  If $U_j \subseteq Y_j$
is an open set for each $j \in I$, and if $U_j = Y_j$ for all but
finitely many $j \in I$, then
\begin{equation}
\label{U = prod_{j in I} U_j, 5}
        U = \prod_{j \in I} U_j
\end{equation}
is an open set in $Y$ with respect to the product topology.  If $U_j$
is a closed set in $Y_j$ for each $j \in I$ too, then $U$ is a closed
set in $Y$ with respect to the product topology, as before.  Using
this, it is easy to see that if $Y_j$ has topological dimension $0$
for each $j$, then $Y$ has topological dimension $0$ with respect to
the product topology.  Similarly, if $Y_j$ is totally separated for
each $j \in I$, then $Y$ is totally separated with respect to the
product topology.

        If $U_j \subseteq Y_j$ is an open set for each $j \in I$,
then (\ref{U = prod_{j in I} U_j, 5}) is an open set in $Y$ with
respect to the strong product topology on $Y$.  If $U_j$ is also a
closed set in $Y_j$ for each $j \in I$, then $U$ is a closed set in
$Y$ with respect to the strong product topology, as in the preceding
paragraph.  If $Y_j$ has topological dimension $0$ for each $j \in I$,
then it follows that $Y$ has topological dimension $0$ with respect to
the strong product topology too, as before.  If $Y_j$ is totally
separated for each $j \in I$, then $Y$ is totally separated with
respect to the strong product topology as well.  However, this follows
from the analogous statement for the product topology on $Y$, since
open and closed subsets of $Y$ with respect to the product topology
have the same property with respect to the strong product topology.

        Let $X$ be another topological space, and let $\phi_j$ be a
continuous mapping from $X$ into $Y_j$ for each $j \in I$.  This leads
to a mapping
\begin{equation}
\label{phi : X to Y}
        \phi : X \to Y
\end{equation}
in a natural way, where the $j$th component of $\phi(x)$ in $Y$ is
equal to $\phi_j(x)$ for each $x \in X$ and $j \in I$.  It is easy to
see that $\phi$ is continuous with respect to the product topology on
$Y$, but this does not always work with respect to the strong product
topology on $Y$.  Suppose now that
\begin{equation}
\label{Y_j = {0, 1}}
        Y_j = \{0, 1\}
\end{equation}
equipped with the discrete topology for each $j \in I$, so that $Y$
has topological dimension $0$ with respect to the corresponding
product topology, as before.  In this case, $\phi_j$ is continuous if
and only if $\phi_j^{-1}(\{0\})$ and $\phi_j^{-1}(\{1\})$ are open
subsets of $X$, which is the same as saying that $\phi_j^{-1}(\{1\})$
is both open and closed in $X$.

        Let $W_j$ be a subset of $X$ that is both open and closed
for each $j \in I$, and let $\phi_j$ be the mapping from $X$ into
$\{0, 1\}$ such that
\begin{equation}
\label{phi_j = 1 on W_j and phi_j = 0 on X setminus W_j}
 \phi_j = 1 \hbox{ on } W_j \quad\hbox{and}\quad \phi_j = 0 \hbox{ on }
                                                    X \setminus W_j
\end{equation}
for each $j \in I$.  Thus $\phi_j$ is continuous for each $j \in I$,
so that the correpsonding mapping $\phi$ from $X$ into the Cartesian
product $Y$ is continuous with respect to the product topology on $Y$,
as in the previous paragraph.  Suppose that for each $x, x' \in X$
with $x \ne x'$ there is a $j \in I$ such that
\begin{equation}
\label{x in W_j and x' not in W_j, or x' in W_j and x not in W_j}
 x \in W_j \hbox{ and } x'\not\in W_j, \quad\hbox{or}\quad x' \in W_j
                                        \hbox{ and } x \not\in W_j.
\end{equation}
This is the same as saying that the corresponding mapping $\phi$ from
$X$ into $Y$ is injective.  Note that $X$ is totally separated if and
only if there is a family of subsets of $X$ that are both open and
closed, and which separates points in $X$ in this way.  If, in
addition to these conditions, $\{W_j\}_{j \in I}$ is a sub-base for
the topology of $X$, then $\phi$ is a homeomorphism from $X$ onto its
image in $Y$, with respect to the topology induced on $\phi(X)$ by the
product topology on $Y$.  Remember that $X$ has topological dimension
$0$ if and only if there is a base for the topology of $X$ consisting
of subsets of $X$ that are both open and closed, as in Section
\ref{dimension 0}.  If there is a sub-base for the topology of $X$
consisting of subsets of $X$ that are both open and closed, then the
corresponding base for the topology of $X$ consisting of finite
intersections of elements of the sub-base has the same property.  If
$X$ also satisfies the first or $0$th separation condition, then any
sub-base for the topology of $X$ separates points as well.

\section{Direct products}
\label{direct products}
\setcounter{equation}{0}

        Let $I$ be a nonempty set, and let $C_j$ be a commutative
topological group for each $j \in I$.  As in Section \ref{cartesian
products}, the Cartesian product
\begin{equation}
\label{C = prod_{j in I} C_j}
        C = \prod_{j \in I} C_j
\end{equation}
is a commutative group, where the group operations are defined
coordinatewise.  We have also seen that $C$ is a commutative
topological group with respect to the corresponding product topology.
Let $D_j$ be a subgroup of $C_j$ for each $j \in I$, so that
\begin{equation}
\label{D = prod_{j in I} D_j}
        D = \prod_{j \in I} D_j
\end{equation}
is a subgroup of $C$.  If $D_j$ is an open set in $C_j$ for each $j$,
and if $D_j = C_j$ for all but finitely many $j$, then $D$ is an open
set in $C$ with respect to the product topology.  If the open
subgroups of $C_j$ form a local base for the topology of $C_j$ at $0$
for every $j \in I$, then one can check that the open subgroups of $C$
form a local base for the product topology on $C$ at $0$.  If, for
each $j \in I$, the intersection of the open subgroups of $C_j$ is
equal to $\{0\}$, then $C$ has the analogous property with respect to
the product topology.

        Similarly, $C$ is a commutative topological group with respect
to the strong product topology, as in Section \ref{strong product topology}.
If $D_j$ is an open subgroup of $C_j$ for each $j \in I$, then
(\ref{D = prod_{j in I} D_j}) is an open subgroup of $C$ with respect
to the strong product topology.  As in the preceding paragraph, if,
for each $j \in I$, the intersection of the open subgroups of $C_j$
is equal to $\{0\}$, then $C$ has the same property with respect to
the strong product topology.  This follows from the analogous statement
for the product topology on $C$, since open subsets of $C$ with respect
to the product topology are open with respect to the strong product
topology as well.

        Let $A$ be another commutative topological group, and let
$\phi_j$ be a continuous homomorphism from $A$ into $C_j$ for each
$j \in I$.  As in the previous section, this leads to a mapping
\begin{equation}
\label{phi : A to C}
        \phi : A \to C,
\end{equation}
whose $j$th component is equal to $\phi_j$ for each $j \in I$.  Under
these conditions, $\phi$ is a continuous homomorphism from $A$ into
$C$, with respect to the product topology on $C$.  If $C_j$ is
equipped with the discrete topology for each $j \in I$, then the open
subgroups of $C$ with respect to the product topology form a local
base for $C$ at $0$, as mentioned earlier.  In this situation,
$\phi_j$ is continuous if and only if the kernel $\phi_j^{-1}(\{0\})$
of $\phi_j$ is an open subgroup of $A$.

        Let $B_j$ be an open subgroup of $A$ for each $j \in I$,
and let
\begin{equation}
\label{C_j = A / B_j}
        C_j = A / B_j
\end{equation}
be the quotient of $A$ by $B_j$, equipped with the discrete topology.
Also let $\phi_j$ be the corresponding quotient mapping from $A$ onto
$C_j$ for each $j \in I$, so that
\begin{equation}
\label{phi_j^{-1}({0}) = B_j}
        \phi_j^{-1}(\{0\}) = B_j
\end{equation}
for each $j \in I$.  Thus $\phi_j$ is continuous for each $j \in I$,
as in the preceding paragraph.  This implies that the corresponding
mapping $\phi$ from $A$ into the product $C$ of the $C_j$'s is a
continuous homomorphism with respect to the product topology on $C$,
as before.  By construction,
\begin{equation}
\label{phi^{-1}({0}) = bigcap_{j in I} phi_j^{-1}({0}) = bigcap_{j in I} B_j}
        \phi^{-1}(\{0\}) = \bigcap_{j \in I} \phi_j^{-1}(\{0\})
                         = \bigcap_{j \in I} B_j,
\end{equation}
so that $\phi$ is injective if and only if
\begin{equation}
\label{bigcap_{j in I} B_j = {0}}
        \bigcap_{j \in I} B_j = \{0\}.
\end{equation}
If (\ref{bigcap_{j in I} B_j = {0}}) holds, and if $\{B_j\}_{j \in I}$
is a local sub-base for the topology of $A$ at $0$, then $\phi$ is a
homeomorphism from $A$ onto its image in $C$, with respect to the
product topology on $C$.  Of course, if $\{0\}$ is a closed set in
$A$, then $\{x\}$ is a closed set in $A$ for every $x \in A$, which
implies that the intersection of all open subsets of $A$ that contain
$0$ is equal to $\{0\}$.  In this case, if $\{B_j\}_{j \in I}$ is a
local sub-base for the topology of $A$ at $0$, then (\ref{bigcap_{j in
    I} B_j = {0}}) holds automatically.

\section{Weak connectedness}
\label{weak connectedness}
\setcounter{equation}{0}

        Let us say that a commutative topological group $A$ is
\emph{weakly connected}\index{weak connectedness} if there are no
proper open subgroups of $A$.  If $A$ is connected as a topological
space in the usual sense, then $A$ is weakly connected, because open
subgroups of $A$ are closed sets in $A$, as in Section \ref{open
  subgroups}.  Remember that the set ${\bf Q}$ of rational numbers has
topological dimension $0$ with respect to the topology induced by the
standard topology on ${\bf R}$, and in particular ${\bf Q}$ is totally
disconnected.  We have also seen that ${\bf Q}$ is weakly connected
as a commutative topological group with respect to addition and this
topology, as in Section \ref{dimension 0, continued}.

        Suppose that $A_1$, $A_2$ are commutative topological groups,
and that $\phi$ is a continuous homomorphism from $A_1$ into $A_2$.
If $B_2$ is an open subgroup of $A_2$, then $\phi^{-1}(B_2)$ is an
open subgroup of $A_1$.  If $A_1$ is weakly connected, then
\begin{equation}
\label{phi^{-1}(B_2) = A_1}
        \phi^{-1}(B_2) = A_1,
\end{equation}
which is to say that
\begin{equation}
\label{phi(A_1) subseteq B_2}
        \phi(A_1) \subseteq B_2.
\end{equation}
As before, $B_2$ is a closed set in $A_2$, so that (\ref{phi(A_1)
  subseteq B_2}) implies that
\begin{equation}
\label{overline{phi(A_1)} subseteq B_2}
        \overline{\phi(A_1)} \subseteq B_2,
\end{equation}where $\overline{\phi(A_1)}$ is the closure of $\phi(A_1)$
in $A_2$.  If $\phi(A_1)$ is dense in $A_2$, then it follows that
$A_2$ is weakly connected as well.

        Let us say that a subgroup $A_0$ of a commutative topological
group $A$ is weakly connected\index{weak connectedness} if $A_0$ is
weakly connected as a commutative topological group, with respect to
the topology induced by the one on $A$.  If $A_0$ is a weakly
connected subgroup of $A$ and $B$ is an open subgroup of $A$, then
\begin{equation}
\label{A_0 subseteq B}
        A_0 \subseteq B.
\end{equation}
This follows from (\ref{phi(A_1) subseteq B_2}), applied to the
standard inclusion mapping of $A_0$ into $A$ as a continuous
homomorphism, or by observing that $A_0 \cap B$ is a relatively open
subgroup of $A_0$.  Since (\ref{A_0 subseteq B}) holds for every
weakly connected subgroup $A_0$ of $A$, we get that the subgroup of
$A$ generated by all of its weakly connected subgroups is contained in
$B$.  If $A$ is generated by its weakly connected subgroups, then $A$
is weakly connected too.

        Let $A_0$ be a subgroup of a commutative topological group
$A$ again, and let $B_0$ be a subgroup of $A_0$.  Also let
\begin{equation}
\label{B = overline{B_0}}
        B = \overline{B_0}
\end{equation}
be the closure of $B_0$ in $A$, which is a closed subgroup of $A$.
Suppose that $B_0$ is a relatively open subgroup of $A_0$, so that
$B_0$ is relatively closed in $A_0$ too, as in Section \ref{open
  subgroups}.  This implies that
\begin{equation}
\label{B cap A_0 = B_0}
        B \cap A_0 = B_0.
\end{equation}
Because $B_0$ is relatively open in $A$, there is an open set $U
\subseteq A$ such that
\begin{equation}
\label{B_0 = U cap A_0}
        B_0 = U \cap A_0.
\end{equation}
If $A_0$ is a dense subset of $A$, then we get that
\begin{equation}
\label{U subseteq overline{U cap A_0} = overline{B_0} = B}
        U \subseteq \overline{U \cap A_0} = \overline{B_0} = B.
\end{equation}
This implies that $B$ is an open set in $A$, using coninuity of
translations on $A$, and the fact that $0 \in B_0 \subseteq U$.  Note
that $B \ne A$ when $B_0 \ne A_0$, by (\ref{B cap A_0 = B_0}).  It
follows that $A_0$ is weakly connected when $A$ is weakly connected
and $A_0$ is dense in $A$.

\section{Weak connectedness, continued}
\label{weak connectedness, continued}
\setcounter{equation}{0}

        Let $V$ be a topological vector space over the field ${\bf Q}$
of rational numbers, with respect to the standard absolute value
function on ${\bf Q}$.  Suppose that $B$ is an open subset of $V$
which is also a subgroup of $V$ as a commutative group with respect to
addition.  As in Section \ref{topological vector spaces},
\begin{equation}
\label{phi_v(t) = t v, 2}
        \phi_v(t) = t \, v
\end{equation}
defines a continuous linear mapping from ${\bf Q}$ into $V$ for every
$v \in V$.  This implies that $\phi_v^{-1}(B)$ is an open subgroup of
${\bf Q}$ with respect to addition for every $v \in V$.  It follows
that
\begin{equation}
\label{phi_v^{-1}(B) = {bf Q}}
        \phi_v^{-1}(B) = {\bf Q}
\end{equation}
for every $v \in V$, because ${\bf Q}$ is weakly connected as a
commutative topological group.  In particular,
\begin{equation}
\label{v = phi_v(1) in B}
        v = \phi_v(1) \in B
\end{equation}
for every $v \in V$.  This shows that
\begin{equation}
\label{B = V}
        B = V,
\end{equation}
so that $V$ is weakly connected as a commutative topological group
with respect to addition under these conditions.

        Let $k$ be a field, and let $|\cdot|$ be a $q$-absolute
value function on $k$ for some $q > 0$.  Suppose that $|\cdot|$ is
archimedian on $k$, as in Sections \ref{absolute value functions} and
\ref{q-absolute value functions}.  This implies that $k$ has
characteristic $0$, so that there is a natural embedding of ${\bf Q}$
into $k$.  Thus $|\cdot|$ induces a $q$-absolute value function on
${\bf Q}$, which is archimedian as well.  Using Ostrowski's theorem,
as in Section \ref{q-absolute value functions}, we get that this
induced $q$-absolute value function on ${\bf Q}$ is equivalent to the
standard absolute value function, which means that the induced
absolute value function on ${\bf Q}$ is equal to a positive power of
the standard absolute value function on ${\bf Q}$.

        Let $V$ be a topological vector space over $k$.  We can also
think of $V$ as a topological vector space over ${\bf Q}$, using the
natural embedding of ${\bf Q}$ in $k$, and the $q$-absolute value
function on ${\bf Q}$ induced by $|\cdot|$ on $k$.  Because the
induced absolute value function on ${\bf Q}$ is equivalent to the
standard absolute value function on ${\bf Q}$, as in the preceding
paragraph,we can think of $V$ as a topological vector space over ${\bf
  Q}$ with respect to the standard absolute value function on ${\bf
  Q}$ too.  By the remarks at the beginning of the section, $V$ is
weakly connected as a commutative topological group with respect to
addition.

        Of course, if $k = {\bf R}$ or ${\bf C}$ with the standard
absolute value function, then $k$ is connected with respect to the
corresponding topology.  Similarly, if $V$ is a topological vector
space over ${\bf R}$ or ${\bf C}$, with respect to the standard
absolute value function, then $V$ is connected, and in fact pathwise
connected.  In particular, this implies that $V$ is weakly connected
as a commutative topological group with respect to addition, as in
the previous section.

\section{Open subgroups, continued}
\label{open subgroups, continued}
\setcounter{equation}{0}

        Let $k$ be a field with a $q$-absolute value function $|\cdot|$
for some $q > 0$, and let $V$ be a vector space over $k$.  If $E_1$,
$E_2$ are balanced subsets of $V$, then it is easy to see that their
sum
\begin{equation}
\label{E_1 + E_2}
        E_1 + E_2
\end{equation}
is balanced in $V$ too.  Suppose that $U \subseteq V$ contains $0$ and
is symmetric about $0$, and let $B \subseteq V$ be as in (\ref{B =
  bigcup_{j = 1}^infty U_j}).  Equivalently, if we consider $V$ as a
commutative group with respect to addition, then $B$ is the subgroup
of $V$ generated by $U$.  If $U$ is also balanced in $V$, then one can
check that $B$ is balanced in $V$ too.  More precisely, if $U_j$ is as
in Section \ref{open subgroups} for each $j \in {\bf Z}_+$, then one
can verify that $U_j$ is balanced for every $j$, using induction on
$j$.  This implies that $B$ is balanced, because $B$ is defined as the
union of the $U_j$'s.  Suppose now that $V$ is a topological vector
space over $k$, so that $V$ is a commutative topological group with
respect to addition in particular.  If $U \subseteq V$ is an open set
that contains $0$ and is symmetric about $0$, and if $B$ is as in
(\ref{B = bigcup_{j = 1}^infty U_j}) again, then $B$ is an open
subgroup of $V$, as in Section \ref{open subgroups}.  If $|\cdot|$ is
archimedian on $k$, then it follows that $B = V$, as in the previous
section.

        Let $k$ be a field with an ultrametric absolute value function
$|\cdot|$, and let $V$ be a topological vector space over $k$.  Suppose
that $B_0 \subseteq V$ is an open set which is also a subgroup of $V$
as a commutative group with respect to addition.  Suppose too that
there is a nonempty balanced open set $U \subseteq V$ that is
contained in $B_0$.  In particular, this means that $U$ is symmetric
about $0$.  Let $B$ be the subgroup of $V$ generated by $U$ as in
(\ref{B = bigcup_{j = 1}^infty U_j}).  Thus $B$ is a balanced open
subset of $V$ under these conditions, as in the previous paragraphs.
We also have that
\begin{equation}
\label{B subseteq B_0}
        B \subseteq B_0,
\end{equation}
because $U \subseteq B_0$ and $B_0$ is a subgroup of $V$ with respect
to addition.  If $|\cdot|$ is nontrivial on $k$, then there is always
a nonempty balanced open set $U \subseteq B_0$, as in Section
\ref{topological vector spaces}.  Similarly, if $|\cdot|$ is the
trivial absolute value function on $k$, and if the collection of
linear mappings on $V$ corresponding to multiplication by elements of
$k$ is equicontinuous at $0$, then there is always a nonempty balanced
open set $U \subseteq B_0$.

        Let $k$, $V$ be as in the preceding paragraph, and suppose
for the moment that $|\cdot|$ is nontrivial on $k$.  Also let $B$ be a
subgroup of $V$ with respect to addition that is balanced and
absorbing in $V$.  In particular, if $B$ is an open set in $V$, then
$B$ is absorbing, as in Section \ref{topological vector spaces}.  Note
that $B$ is $\infty$-convex in $V$, in the sense described in Section
\ref{balanced q-convexity}, because $B$ is a balanced subgroup of $V$.
If $N_B$ is the Minkowski functional on $V$ corresponding to $B$ as in
(\ref{N_A(v) = inf {|t| : t in k, v in t A}}), then $N_B$ is a
semi-ultranorm on $V$ under these conditions, as in Section
\ref{minkowski functionals}.

        Suppose now that $|\cdot|$ is the trivial absolute value function
on $k$.  In this case, if $B$ is a subgroup of $V$ with respect to addition
that is balanced, then $B$ is a linear subspace of $V$.  If we put
\begin{equation}
\label{N(v) = 0 when v in B and N(v) = 1 when v in V setminus B}
 N(v) = 0 \hbox{ when } v \in B \quad\hbox{and}\quad N(v) = 1 \hbox{ when }
                                                  v \in V \setminus B,
\end{equation}
then it is easy to see that $N$ is a semi-ultranorm on $V$.  This is
the same as the composition of the standard quotient mapping from $V$
onto the quotient vector space $V / B$ with the trivial ultranorm on
$V / B$.  The semi-ultrametric on $V$ associated to $N$ as in
(\ref{d(v, w) = N(v - w), 2}) is the same as the discrete semimetric
on $V$ corresponding to the partition of $V$ into cosets of $B$, as in
Sections \ref{semimetrics, partitions} and \ref{translation-invariant
  semi-ultrametrics}.

\newpage

%\addcontentsline{toc}{section}{Index}

\printindex

\end{document}